\input amstex
\documentstyle{amsppt}
\loadmsbm
% \magnification=\magstep1

\nologo

\TagsOnRight

\NoBlackBoxes

\define\acc{\operatorname{acc}}
\define\Real{\operatorname{Re}}

\define\esssup{\operatorname{ess\,sup}}

\define\supp{\operatorname{supp}}
\define\Lip{\operatorname{Lip}}
\define\diam{\operatorname{diam}}
\define\dist{\operatorname{dist}}

\def\floor{\mathbin{\hbox{\vrule height1.2ex width0.8pt depth0pt
        \kern-0.8pt \vrule height0.8pt width1.2ex depth0pt}}}

\font\letter=cmss10

\font\normal=cmss10 scaled 700

\font\bignormal=cmss10 scaled 1300

\define\Box{\text{\normal B}}

\define\Haus{\text{\normal H}}

\define\scon{\text{\normal con}}

\define\ssim{\text{\normal sim}}

\define\erg{\text{\normal erg}}

\define\con{\text{\normal c}}
\define\abs{\text{\normal ab}}

\define\distance{\text{\letter d}}

\define\Distance{\text{\letter D}}

\define\LDistance{\text{\letter L}}

\define\X{\text{\bignormal X}}

%
%\hsize = 5.5 true in
%%\hoffset = 0.5 true in
%\vsize = 8.5 true in
%

\hsize = 6 true in
%\hoffset = 0.5 true in
\vsize = 9.18 true in

\topmatter
\title
Multifractal spectra and
multifractal zeta-functions
\endtitle
\endtopmatter

\centerline{\smc V\. Mijovi\'c}
\centerline{Department of Mathematics}
\centerline{University of St\. Andrews}
\centerline{St\. Andrews, Fife KY16 9SS, Scotland}
\centerline{e-mail: {\tt vm27\@st-and.ac.uk}}

\medskip

\centerline{\smc L\. Olsen}
\centerline{Department of Mathematics}
\centerline{University of St\. Andrews}
\centerline{St\. Andrews, Fife KY16 9SS, Scotland}
\centerline{e-mail: {\tt lo\@st-and.ac.uk}}

\topmatter
\abstract
{
We introduce multifractal zeta-functions providing precise information 
of  a very general class of multifractal spectra, including, for example, the 
multifractal spectra of self-conformal 
 measures 
 and the multifractal spectra of ergodic Birkhoff averages of continuous functions.
 More precisely, we prove that these and more
 general
 multifractal spectra equal the 
abscissae of convergence of the associated zeta-functions.
}
\endabstract
\endtopmatter

\footnote""
{
\!\!\!\!\!\!\!\!
}

\leftheadtext{V\. Mijovi\'c \& L\. Olsen}

\rightheadtext{Multifractal spectra and multifractal zeta-functions}

\heading{1. Introduction.}\endheading

Measures with widely varying intensity
are called multifractals and have 
during the past 20 years
been
the focus of enormous 
 attention in the mathematical
 literature.
% and
%the purpose of this paper is to
%analyse the \lq\lq spectra" of
%(very general)
%multifractals using zeta-functions.
Loosely
speaking
there are two main ingredients
in multifractal analysis: the multifractal spectrum and the 
Renyi dimensions.
One of the main goals
in multifractal analysis is to understand these two ingredients and their 
relationship with each other.
It is generally believed by experts that
the multifractal spectrum 
% of a measure
and the 
Renyi dimensions of a measure
encode important
information about the measure,
and it is therefore of considerable importance
to find explicit 
formulas 
for these quantities.
In [Le-VeMe,Ol4,Ol5,Ol6]
the authors
used the zeta-function technique introduced and pioneered 
by M\. Lapidus et al
in the intriguing books
[Lap-vF1,Lap-vF2]
in order to 
find 
explicit formulas for the 
Renyi dimensions of
a self-similar measure. 
At this point we note that
% We now note that 
it is generally
%(and with substantial evidence)
believed that analysing
the multifractal spectrum of a measure is considerably more difficult
and challenging  than analysing 
its Renyi dimensions,
and
the main purpose of this paper is 
to address the 
substantially more difficult problem 
of
finding explicit formulas 
%performing 
%%a similar analysis 
%an analysis
 for the multifractal spectrum
of
a self-similar measure
similar to the 
explicit
formulas 
for its Renyi dimensions
found
in
[Le-VeMe,Ol4,Ol5,Ol6].
In particular, and as
a first step 
in this direction, 
%in this work 
we introduce multifractal zeta-functions providing precise information 
of very general classes of multifractal spectra, including, for example, the 
multifractal spectra of self-conformal 
 measures 
 and the multifractal spectra of ergodic Birkhoff averages of continuous functions.
 More precisely, we prove that these and more
 general
 multifractal spectra equal the 
abscissae of convergence of the associated zeta-functions.
%A more careful analysis 
%of the multifractal zeta-functions introduced in this work (see [MiOl])
%suggests
%that they can by used
%to
% provide explicit formulas for the multifractal spectra
%allowing 
%us to express 
%the
% multifractal spectra
%as sums involving the residues of these zeta-functions.

\bigskip

{\bf 1.1. The first ingredient in multifractal analysis: multifractal spectra.}
For a Borel measure $\mu$ on $\Bbb R^{d}$ 
with support equal to $K$
and
a positive number $\alpha$,
let us consider
the set $\Delta_{\mu}(\alpha)$ of
those points
$x$ in $\Bbb R^{d}$ for which the measure
$\mu(B(x,\delta))$ of the ball
$B(x,\delta)$ with center $x$ and radius $\delta$ behaves like
$\delta^{\alpha}$ for small $\delta$,
i\.e\. the set
 $$
   \Delta_{\mu}(\alpha)
   =
   \Bigg\{
   x\in K
   \,\Bigg|\,
   \lim_{r\searrow 0}
   \frac{\log\mu(B(x,r))}{\log r}
   =
   \alpha
   \Bigg\}\,.
 $$
If the intensity of the measure $\mu$ varies very widely, it may
happen that the sets $\Delta_{\mu}(\alpha)$
display a
fractal-like character for a range of values of $\alpha$. In this case
it is natural to study
the Hausdorff dimensions of the sets $\Delta_{\mu}(\alpha)$ as $\alpha$
varies.
We therefore define the the multifractal spectrum of $\mu$
by
 $$
 f_{\mu}(\alpha)
 =
 \dim_{\Haus}\Delta_{\mu}(\alpha)\,,
 \tag2.1
$$
where $\dim_{\Haus}$ denotes the Hausdorff dimension.
Here and below
we use the following convention,
namely, we define the Hausdorff
of the empty set to be $-\infty$, i\.e\. we put
 $$
 \dim_{\Haus}\varnothing
 =
 -
 \infty\,.
 $$
One of the main problems
in multifractal analysis is to study this and related functions.
The function $f_{\mu}(\alpha)$ was
first explicitly defined by the physicists Halsey et al. in 1986 in
their seminal paper [HaJeKaPrSh].

The multifractal spectrum
$f_{\mu}$ is defined using the Hausdorff dimension.
There is an alternative approach using
\lq\lq box-counting" arguments
leading to the coarse multifractal spectrum.
Namely,
for a   Borel probability measure $\mu$ on $\Bbb R^{d}$ 
with support equal to $K$
and a real number  $\alpha$,
the coarse multifractal spectrum is defined as follows.
For positive real numbers $r>0$ and $\delta>0$, we write
 $$
 \align
 N_{\mu,\delta}(\alpha;r)
 =
 \sup
 \Bigg\{
 |I|
 \,\Bigg|\,
&\text{$(B(x_{i},\delta))_{i\in I}$ is a finite family of balls such that:}\\
&\qquad\qquad\text{$x_{i}\in K$ for all $i$,}\\
&{}\\
&\qquad\qquad\text{$B(x_{i},\delta)\cap B(x_{j},\delta)=\varnothing$ for all $i\not=j$,}\\
&\qquad\qquad
 \alpha-r
 \le
 \frac{\log\mu(B(x_{i},\delta))}{\log\delta}
 \le
 \alpha+r
 \,\,\,\,
 \text{for all $i$}
 \Bigg\}\,,
 \tag1.3
 \endalign
 $$
and define the 
$r$-approximate
coarse multifractal spectrum $f_{\mu}^{\con}(\alpha;r)$ of $\mu$ by
 $$
 f_{\mu}^{\con}(\alpha;r)
 =
 \liminf_{\delta\searrow 0}
 \frac{\log N_{\mu,\delta}(\alpha;r)}{-\log \delta}\,.
 \tag1.4
 $$
Finally, the 
coarse multifractal spectrum $f_{\mu}^{\con}(\alpha)$ of $\mu$ is defined by
 $$
 f_{\mu}^{\con}(\alpha)
 =
 \lim_{r\searrow 0}
 f_{\mu}^{\con}(\alpha;r)
 \tag1.5
  $$
(it is clear that this limit exists since 
$f_{\mu}^{\con}(\alpha;r)$ 
is a monotone function of $r$).
We note that it is easily seen that
$$
 f_{\mu}(\alpha)
 \le
 f_{\mu}^{\con}(\alpha)\,,
$$
and that
this inequality 
may be strict, see, for example, [Fa1].

\bigskip

{\bf 1.2. The second ingredient in multifractal analysis: Renyi dimensions.}
Renyi dimensions
quantify
the 
varying intensity of a measure
by 
analyzing  its moments at different scales.
Formally,
Renyi dimensions
are defined as follows.
Let $\mu$ be a Borel measure on $\Bbb R^{d}$.
For $E\subseteq \Bbb R^{d}$,
$q\in\Bbb R$ and $\delta>0$, we 
define 
the $q$-moment $M_{\mu,\delta}(q;E)$ of 
 $\mu$ on $E$ at scale $\delta$ 
by
 $$
 \align
 M_{\mu,\delta}(q;E)
 =
 \sup
 \Bigg\{
 \sum_{i\in I}\mu(B(x_{i},\delta))^{q}
 \,\Bigg|\,
&\text{$(B(x_{i},\delta))_{i\in I}$ is a finite family of balls such that:}\\
&\qquad\qquad\text{$x_{i}\in K$ for all $i$,}\\
&\qquad\qquad\text{$B(x_{i},\delta)\cap B(x_{j},\delta)=\varnothing$ for all $i\not=j$}
 \Bigg\}\,,
 \endalign
 $$
We now
 define the 
lower and upper
Renyi spectra
$\underline\tau_{\mu}(\cdot;E),\overline\tau_{\mu}(\cdot;E):\Bbb R\to[-\infty,\infty]$
 of $\mu$ by
 $$
 \aligned
 \underline\tau_{\mu}(q;E)
&= 
 \,\liminf_{\delta\searrow 0}\,\frac{\log M_{\mu,\delta}(q;E)}{-\log \delta}
 \,,\\
  \overline\tau_{\mu}(q;E)
&=  
 \limsup_{\delta\searrow 0}\frac{\log M_{\mu,\delta}(q;E)}{-\log \delta}
 \,.
 \endaligned
 $$
If $E$ equals the support $\supp\mu$ of $\mu$, then we will use the following 
shorter notation
 $$
 M_{\mu,\delta}(q)
 =
 M_{\mu,\delta}(q;\supp\mu)\,,\,\,\,\,
 \underline\tau_{\mu}(q)
 =
 \underline\tau_{\mu}(q;\supp\mu)\,,\,\,\,\,
 \overline\tau_{\mu}(q)
 =
 \overline\tau_{\mu}(q;\supp\mu)\,.
 $$
We note that
the $q$-moment $M_{\mu,\delta}(q;E)$ 
is closely related to the box dimension $\dim_{\Box}E$ of $E$.
Indeed,
if  
we let 
$M_{\delta}(E)$
denote the greatest number of 
pairwise disjoint balls
of radii
$\delta$
with centers in $E$, then
it follows from the definition of the box dimension that
$\dim_{\Box}E=\lim_{\delta\to 0}\frac{\log M_{\delta}(E)}{-\log\delta}$
(provided the limit exists)
and we clearly have
 $$
 M_{\delta}(E)
 =
 M_{\mu,\delta}(0;E)\,.
 \tag1.6
 $$
%provided
%$E\subseteq\supp\mu$. 

 It is also possible to define an integral version of 
 the $q$-moments $M_{\mu,\delta}(q;E)$.
Namely, 
for $E\subseteq\Bbb R^{d}$,
$q\in\Bbb R$ and $\delta>0$, we 
define 
the integral $q$-moment $V_{\mu,\delta}(q)$ of 
 $\mu$ on $E$ at scale $\delta$ 
by
 $$
 V_{\mu,\delta}(q;E)
 =
 \,
 \int
 \limits_{B(E,\delta)}
 \mu(B(x,\delta))^{q}
 \,
 d\Cal L^{d}(x)
 $$  
where
$B(E,\delta)
=
\{x\in\Bbb R^{d}\,|\,
\dist(x,E)\le \delta\}$
and
$\Cal L^{d}$ denotes the Lebesgue measure in $\Bbb R^{d}$. 
We now
 define the 
lower and upper
integral
Renyi spectra
$\underline T_{\mu}(\cdot;E),
\overline T_{\mu}(\cdot;E):\Bbb R\to[-\infty,\infty]$
 of $\mu$ by
 $$
 \aligned
 \underline T_{\mu}(q;E)
&= 
 \,\liminf_{\delta\searrow 0}\,\frac{\log V_{\mu,\delta}(q;E)}{-\log \delta}
 \,,\\
  \overline T_{\mu}(q;E)
&=  
 \limsup_{\delta\searrow 0}\frac{\log V_{\mu,\delta}(q;E)}{-\log \delta}
 \,.
 \endaligned
 $$
As above,
if $E$ equals the support $\supp\mu$ of $\mu$, then we will use the following 
shorter notation
 $$
  V_{\mu,\delta}(q)
 =
 V_{\mu,\delta}(q;\supp\mu)\,,\,\,\,\,
 \underline T_{\mu}(q)
 =
 \underline T_{\mu}(q;\supp\mu)\,,\,\,\,\,
 \overline T_{\mu}(q)
 =
 \overline T_{\mu}(q;\supp\mu)\,.
 $$ 
As above, we note  that
the integral $q$-moment $V_{\mu,\delta}(q;E)$ 
is also closely related to the 
Minkowski volume of $E$
and the box dimension $\dim_{\Box}E$ of $E$.
Namely,
 if we let $V_{\delta}(E)$ denote the
$\delta$ approximate Minkowski volume of $E$, i\.e\.
$V_{\delta}(E)= \Cal L^{d}(\,B(E,\delta)\,)$, then
it is well-known that
 $\dim_{\Box}E
 =
\lim_{\delta\to 0}\frac{\log( \frac{1}{r^{d}}V_{\delta}(E))}{-\log\delta}$ 
(provided the limit exists)
and we clearly have
 $$
 V_{\delta}(E)
 =
 V_{\mu,\delta}(0;E)\,.
 \tag1.7
 $$
%provided
%$E\subseteq\supp\mu$. 

\bigskip

{\bf 1.3. The Multifractal Formalism.}
Based on a remarkable insight
together with a clever heuristic argument,
it was suggested by theoretical physicists
Halsey et al\. [HaJeKaPrSh]
that the multifractal spectra
$f_{\mu}$ and $f_{\mu}^{\con}$ can be computed 
using the
Renyi
dimensions.
This result is
known
as the
\lq\lq Multifractal Formalism" in the physics literature.
More precisely, the
\lq\lq Multifractal Formalism"
says
that
the
multifractal spectra
equal the 
Legendre transform of the 
Renyi dimensions.
Recall
that if
$\varphi:\Bbb R\to\Bbb R$ is a real valued function,
then the Legendre transform $\varphi^{*}:\Bbb R\to[-\infty,\infty]$
of $\varphi$ is defined by
 $$
 \varphi^{*}(x)
 =
 \inf_{y}(xy+\varphi(y))\,.
 \tag1.8
 $$
We can now state the 
\lq\lq Multifractal Formalism".

\bigskip

\proclaim{The Multifractal Formalism -- A Physics Folklore Theorem}
The multifractal spectrum
$f_{\mu}$ 
of $\mu$
and the coarse multifractal spectrum $f_{\mu}^{\con}$
of $\mu$
equal the 
Legendre transforms
$\underline\tau_{\mu}^{*}$, $\overline\tau_{\mu}^{*}$,
$(\underline T_{\mu})^{*}$ and 
$(\overline T_{\mu})^{*}$
of 
the Renyi dimensions, i\.e\.
 $$
 f_{\mu}(\alpha)
 =
 f_{\mu}^{\con}(\alpha)
 =
 \underline\tau_{\mu}^{*}(\alpha)
 =
 \overline\tau_{\mu}^{*}(\alpha)
 =
 \underline T_{\mu}^{*}(\alpha)
 =
 \overline T_{\mu}^{*}(\alpha)
 $$
 for all $\alpha$.
\endproclaim

\bigskip

\noindent
During the past 20 years
there has been an enormous interest 
in
verifying the
Multifractal Formalism 
and
computing the multifractal spectra of measures
in 
the mathematical literature.
In the mid 1990's
Cawley \& Mauldin [CaMa] and Arbeiter \& Patzschke [ArPa]
verified the Multifractal Formalism for self-similar measures 
satisfying the OSC,
and within the last
15 years the multifractal spectra of various classes of measures
in Euclidean space $\Bbb R^{d}$ 
exhibiting some degree of self-similarity have been computed 
rigorously, cf\. 
the textbooks [Fa2,Pe]
and the references therein.
Summarizing 
the previous paragraph
somewhat more succinctly,
previous work has almost entirely
% focused
concentrated on the following problem:

\bigskip

\block
{\it Previous work:}
% has concentrated on the following problem:}

\noindent
Previous work
has concentrated on
 finding the limiting behaviour of the following ratios, namely,
  $$
 \frac{\log M_{\mu,\delta}(q)}{-\log \delta}
  $$
 and
 $$
 \frac{\log N_{\mu,\delta}(\alpha;r)}{-\log\delta}\,.
 $$ 
Indeed,
computing the Renyi dimensions $ \underline\tau_{\mu}(q)$ and $ \overline\tau_{\mu}(q)$
involves
analysing the
limiting behaviour of
$\frac{\log I_{\mu,r}(q)}{-\log r}$,
and
computing the   coarse multifractal spectrum $f_{\mu}^{\con}(\alpha;r)$
involves
analysing the
limiting behaviour of 
$\frac{\log N_{\mu,\delta}(\alpha;r)}{-\log\delta}$.
\endblock

\bigskip

\noindent
Due to the importance of the 
quantities
$M_{\mu,\delta}(q)$
 and
$N_{\mu,\delta}(\alpha;r)$
%encode important  geometric information of
%the measure $\mu$,
%and for this reason 
it is  clearly 
desirable not only
to  find 
expressions 
for the 
limiting behaviour of   
$\frac{\log M_{\mu,\delta}(q)}{-\log \delta}$
and
$\frac{\log N_{\mu,\delta}(\alpha;r)}{-\log\delta}$,
but to find explicit expressions for the 
quantities
$M_{\mu,\delta}(q)$
 and
$N_{\mu,\delta}(\alpha;r)$
themselves.
The purpose of this work can bee seen as a first step in this direction.
%and
%the ultimate aim of the present works can be stated as follows.
Again,
summarizing 
this
somewhat more succinctly,
the present work
is
concentrated on the following problem:

\bigskip

\block
{\it Present work:}
% concentrates on the following problem:}

\noindent
This work explores methods of finding 
 explicit
 expressions for
  $$
 M_{\mu,\delta}(q)
  $$
 and
 $$
 N_{\mu,\delta}(\alpha;r)\,.
 $$ 
 \endblock

\bigskip

\noindent
It is clear that finding
explicit
 expressions for
 $M_{\mu,\delta}(q)$
 and
$N_{\mu,\delta}(\alpha;r)$
is a more challenging
undertaking
 than 
determining
the limiting behaviour of the ratios
$\frac{\log M_{\mu,\delta}(q)}{-\log \delta}$
 and
 $\frac{\log N_{\mu,\delta}(\alpha;r)}{-\log\delta}$;
 indeed,
 if  
 explicit
 expressions for
 $M_{\mu,\delta}(q)$
 and
$N_{\mu,\delta}(\alpha;r)$
are known, then
the limiting behaviour of the ratios
$\frac{\log M_{\mu,\delta}(q)}{-\log \delta}$
 and
 $\frac{\log N_{\mu,\delta}(\alpha;r)}{-\log\delta}$
 can be computed directly
 from these expressions.

We will now describe our  strategy for 
analysing the 
quantities
 $M_{\mu,\delta}(q)$
 and
$N_{\mu,\delta}(\alpha;r)$.
Very loosely speaking, the 
quantities
$M_{\mu,\delta}(q)$
 and
$N_{\mu,\delta}(\alpha;r)$
\lq\lq count"
the number of 
balls $B(x,\delta)$ 
satisfying certain conditions.
There are
two distinct and widely used
techniques for 
analysing the asymptotic behaviour
%or
%finding explicit expression for 
of
such (and similar)
\lq\lq counting functions", namely,
(1) using ideas from renewal theory 
or
(2) 
using the Mellin transform
and the residue theorem  to express 
the
\lq\lq counting functions"
 as sums
involving the residues of suitably defined
 zeta-functions.
Indeed,
renewal theory techniques were introduced and pioneered by Lalley
[La1,La2,La2]
in the 1980's,
and later investigated further by Gatzouras [Ga], Winter [Wi]
and most recently Kesseb\"ohmer \& Kombrink [KeKo],
in order
to analyse the asymptotic behaviour
of 
the
\lq\lq counting function"
$M_{\delta}(E)=M_{\mu,\delta}(0,E)=M_{\mu,\delta}(0)$ 
for self-similar sets $E$
% and self-similar measures $\mu$ with support equal to $E$, 
(see (1.6))
and similar \lq\lq counting functions" from fractal geometry.
%%to analyse a number 
%%of 
%%\lq\lq counting functions"
%%in fractal geometry.
%to find 
%expressions 
%for the
%Minkowski content
%of a self-similar set
%(in fact, 
%Lalley [La1,La2,La3]
%found an explicit
%expression
%for the smallest number of 
%balls needed to cover a self-similar set, but 
%it is clear that
%an identical
%argument 
%can be used to analyze the Minkowski content).
However,
while renewal theory
techniques are powerful tools
for analysing the asymptotic behaviour of
\lq\lq counting functions",
they do not yield 
% genuinely
\lq\lq explicit" 
formulas.
%expressions for the 
%\lq\lq counting functions"
%they are applied to.
%%instead, they provide 
%%expressions 
%%involving
%%complicated
%%(and
%%difficult to compute)
%%integrals of 
%%certain auxiliary
%%functions.
This is clearly
unsatisfactory and it would be desirable if \lq\lq explicit" expressions could be found. 
%In fact,
%even in the fractal case, the problem of nding explicit formulas for the Minkowski content
%is highly non-trivial. 
However, despite, or perhaps in spite, of the diculties, 
the problem of finding \lq\lq explicit" formulas
of
\lq\lq counting functions"
in fractal geometry
%this problem
has recently attracted considerable interest. 
In particular, Lapidus and collaborators
[LapPea1,LapPea2,LapPeaWi,Lap-vF1,Lap-vF2] have 
with
spectacular success
during the past 20 years 
 pioneered the use of 
applying the Mellin transform
to suitably defined
zeta-functions in order to obtain explicit formulas for the
Minkowski volume
$V_{\delta}(E)=V_{\mu,\delta}(0,E)=V_{\mu,\delta}(0)$
 of self-similar fractal subsets $E$ of the line
 (see (1.7)).

It would clearly be desirable if 
similar
formulas
 could be found
for the multifractal 
quantities
$M_{\mu,\delta}(q)$
 and
$N_{\mu,\delta}(\alpha;r)$ of 
self-similar (and more general) multifractal  measures $\mu$.
In multifractal 
analysis
it is generally
believed 
that 
 analysing
 the
 the $q$-moments $M_{\mu,\delta}(q)$
 and 
 the associated
 Renyi dimenions
 $\underline\tau_{\mu}^{*}(\alpha)$
 and
 $\overline\tau_{\mu}^{*}(\alpha)$
is 
less difficult than analysing the
\lq\lq counting function" $N_{\mu,\delta}(\alpha;r)$
and the
associated
multifractal spectra
$f_{\mu}$ and $f_{\mu}^{\con}$.
Indeed,
in [Le-VeMe,Ol4]
(see also the surveys [Ol5,Ol6])
the authors
 introduced 
 a one-parameter family
 of  multifractal
zeta-functions
and
established
explicit formulas for
the integral $q$-moments
$V_{\mu,\delta}(q)$
expressing $V_{\mu,\delta}(q)$
as a sum
involving the residues of these zeta-functions,
%using the Mellin transform
%and the residue theorem, we are able to express the multifractal tube formulas as sums
%involving the residues of the zeta-function.
and in [Ol1]
the asymptotic behaviour of the 
$q$-moments
$M_{\mu,\delta}(q)$
were analysed using
techniques from
renewal theory.
In addition,
we note that Lapidus and collaborators have introduced various intriguing multifractal 
zeta-functions
[LapRo,LapLe-VeRo]. However, the multifractal zeta-functions in 
[LapRo,LapLe-VeRo] serve very 
different
purposes and are 
significantly different from the multifractal
zeta-functions
introduced in  [Le-VeMe,Ol2,Ol4].
The purpose of this
paper is to
address the significantly more difficult and challenging
problem 
of 
performing
a similar analysis of the multifractal spectrum
\lq\lq counting function" $N_{\mu,\delta}(\alpha;r)$.
In particular,
the final aim is to introduce 
a class of multifractal zeta-functions
allowing
us
to derive
explicit formulas for
the \lq\lq counting function" $N_{\mu,\delta}(\alpha;r)$
expressing $N_{\mu,\delta}(\alpha;r)$
as a sum
involving the residues of these zeta-functions.
As a first step in this direction, 
in
this work we 
introduce multifractal
zeta-functions
providing
precise
information of very general classes
of
multifractal spectra, including, for example, the
spectra
$f_{\mu}$ and $f_{\mu}^{\con}$ of self-similar 
multifractal measures $\mu$.
More precisely, 
we prove that
the 
multifractal spectra
equal
the
abscissae
of convergence 
of the associated zeta-functions.
It is our hope 
that a more careful analysis
of these zeta-functions 
will provide
explicit formulas for
the \lq\lq counting function" $N_{\mu,\delta}(\alpha;r)$
allowing us to
express $N_{\mu,\delta}(\alpha;r)$
as a sum
involving the residues of these zeta-functions; this will be explored in [MiOl].
In order to illustrate the
ideas involved
we now consider a simple 
example.

\bigskip

{\bf 1.4. An example illustrating the ideas: 
self-similar measures.}
To illustrate the above ideas in a simple setting, 
we consider the following example
involving self-similar
measures.
Recall, that self-similar measures are defined as follows.
Let $(S_{1},\ldots,S_{N})$ be a list of contracting similarities 
$S_{i}:\Bbb R^{d}\to\Bbb R^{d}$
and let $r_{i}$ denote the similarity ratio of $S_{i}$.
Also,
let
$(p_{1},\ldots,p_{N})$ be a probability vector.
Then there is 
a unique Borel probability measure $\mu$ on $\Bbb R^{d}$ such that
 $$
 \mu
 =
 \sum_{i}p_{i}\mu\circ S_{i}^{-1}\,,
 \tag1.9
 $$
 see [Fa1,Hu].
The measure $\mu$ is called the
self-similar measure associated with the list
$(S_{1},\ldots,\allowmathbreak S_{N},\allowmathbreak p_{1},\allowmathbreak \ldots,p_{N})$.
If the so-called Open Set Condition (OSC) is satisfied,
then
the multifractal spectra
$f_{\mu}$ and $f_{\mu}^{\con}$
are given by the following formula.
Namely, 
if if the OSC is satisfied
and
if we define
$\beta:\Bbb R\to\Bbb R$ by
 $$
 \sum_{i}p_{i}^{q}r_{i}^{\beta(q)}=1\,,
 \tag1.10
 $$
then it follows from [,CaMa,Pa] that
 $$
 f_{\mu}(\alpha)
 =
 f_{\mu}^{\con}(\alpha)
 =
 \beta^{*}(\alpha)
 $$
for all $\alpha\in\Bbb R$
 where $\beta^{*}$ denotes the Legendre transform of $\beta$ 
 (recall, that the definition of the Legendre transform is given in (1.8)).

For $\alpha\in\Bbb R$,
we
are now attempting to introduce a 
\lq\lq natural" self-similar multifractal zeta-function $\zeta_{\alpha}^{\ssim}$
whose 
abscissa 
of convergence equals
$f_{\mu}(\alpha)$.
To do this we first introduce the following notation.
Write
$\Sigma^{*}=\{\bold i=i_{1}\ldots i_{n}\,|\,n\in\Bbb N\,,\, i_{j}\in\{1,\ldots,N\}\,\}$
i\.e\.
$\Sigma^{*}$
is the set of all finite strings
$\bold i=i_{1}\ldots i_{n}$ with $n\in\Bbb N$ and  $i_{j}\in\{1,\ldots,N\}$.
For a finite string
$\bold i=i_{1}\ldots i_{n}\in\Sigma^{*}$
of length 
$n$, we write
$|\bold i|=n$,
and we write
$r_{\bold i}=r_{i_{1}}\cdots r_{i_{n}}$
and
$p_{\bold i}=p_{i_{1}}\cdots p_{i_{n}}$.
With this notation, we can now motivate the 
introduction of a
\lq\lq natural"
multifractal zeta-function as follows.
Namely,
since 
$f_{\mu}(\alpha)$
measures 
the size of the 
set of points 
$x$ for which
$\lim_{\delta\searrow 0}\frac{\log\mu(B(x,\delta))}{\log \delta}=\alpha$
and since
$\frac{\log\mu(B(x,\delta))}{\log \delta}$
has the same form as
$\frac{\log p_{\bold i}}{\log r_{\bold i}}$,
it is natural to
define the 
self-similar multifractal zeta-function $\zeta_{\alpha}^{\ssim}$ by
 $$
 \zeta_{\alpha}^{\ssim}(s)
 =
 \sum
 \Sb
 \bold i\\
 {}\\
 \frac{\log p_{\bold i}}{\log r_{\bold i}}=\alpha
 \endSb
 r_{\bold i}^s
 \tag1.11
 $$
for those complex numbers $s$ for which the series converges absolutely.
An easy and straight forward calculation
(which we present below)
shows that
the abscissa of convergence
$\sigma_{\abs}(\,\zeta_{\alpha}^{\ssim}\,)$
 of $\zeta_{\mu}^{\ssim}$
is less than $f_{\mu}(\alpha)$,
i\.e\.
 $$ 
 \sigma_{\abs}(\,\zeta_{\alpha}^{\ssim}\,)\le f_{\mu}(\alpha)=f_{\mu}^{\con}(\alpha)\,.
 \tag1.12
 $$
Indeed, if
$\alpha\not\in[\min_{i}\frac{\log p_{i}}{\log r_{i}},\max_{i}\frac{\log p_{i}}{\log r_{i}}]$,
then it is easily seen that
that 
for all $\bold i\in\Sigma^{*}$, we have
$\frac{\log p_{\bold i}}{\log r_{\bold i}}\not=\alpha$, whence
$\sigma_{\abs}(\,\zeta_{\alpha}^{\ssim}\,)=-\infty$,
and inequality (1.12) is therefore trivially satisfied.
On the other hand, if
$\alpha\in[\min_{i}\frac{\log p_{i}}{\log r_{i}},\max_{i}\frac{\log p_{i}}{\log r_{i}}]$,
then it follows from [CaMa,Fa1,Pa]
that there is a (unique)
$q\in\Bbb R$ with 
$f_{\mu}(\alpha)
=
f_{\mu}^{\con}(\alpha)
=
\alpha q+\beta(q)$.
Hence, for each $\varepsilon>0$,
we have
(using the fact that
$\sum_{i}
 p_{i}^{q}
 r_{i}^{\beta(q)+\varepsilon}<1$)
 $$
 \align
 \zeta^{\ssim}_{\alpha}\big(\,f_{\mu}(\alpha)+\varepsilon\,\big)
&=
 \sum
 \Sb
 \bold i\\
 {}\\
 \frac{\log p_{\bold i}}{\log r_{\bold i}}=\alpha
 \endSb
 r_{\bold i}^{f_{\mu}(\alpha)+\varepsilon}\\
&{}\\ 
&=
 \sum
 \Sb
 \bold i\\
 {}\\
 \frac{\log p_{\bold i}}{\log r_{\bold i}}=\alpha
 \endSb
 r_{\bold i}^{\alpha q+\beta(q)+\varepsilon}\\
&{}\\ 
&=
 \sum
 \Sb
 \bold i\\
 {}\\
 \frac{\log p_{\bold i}}{\log r_{\bold i}}=\alpha
 \endSb
 p_{\bold i}^{q}
 r_{\bold i}^{\beta(q)+\varepsilon}\\
&{}\\ 
&\le
 \sum_{\bold i}
 p_{\bold i}^{q}
 r_{\bold i}^{\beta(q)+\varepsilon}\\
&=
\sum_{n}
 \sum_{|\bold i|=n}
 p_{\bold i}^{q}
 r_{\bold i}^{\beta(q)+\varepsilon}\\ 
&=
 \sum_{n}
 \Bigg(
 \sum_{i}
 p_{i}^{q}
 r_{i}^{\beta(q)+\varepsilon}
 \Bigg)^{n}\\
&<
 \infty\,.  
 \endalign
 $$
This shows that $\sigma_{\abs}(\,\zeta_{\alpha}^{\ssim}\,)\le f_{\mu}(\alpha)+\varepsilon$. 
Letting $\varepsilon\searrow 0$, now gives
$\sigma_{\abs}(\,\zeta_{\alpha}^{\ssim}\,)\le f_{\mu}(\alpha)$.
This proves (1.12).

However, it is also clear that we, in general, do not
have equality in (1.12).
Indeed, the set 
$\{\frac{\log p_{\bold i}}{\log r_{\bold i}}
 \,|\,
 \bold i\in\Sigma^{*}
 \}$ 
is clearly countable 
(because $\Sigma^{*}$ is countable)
and 
if 
$\alpha\in
 \Bbb R
 \setminus
 \{\frac{\log p_{\bold i}}{\log r_{\bold i}}
 \,|\,
 \bold i\in\Sigma^{*}
 \}$, then
$\sigma_{\abs}(\,\zeta_{\alpha}\,)
  =
  -\infty$
(because the series 
(1.11)
that defines $\zeta_{\alpha}^{\ssim}(s)$
is obtained by summing over the empty set).
Since it also
follows from [CaMa,Fa1,Pa]
 that
 $f_{\mu}(\alpha)=f_{\mu}^{\con}(\alpha)>0$
for all
 $\alpha\in  (\min_{i}\frac{\log p_{i}}{\log r_{i}},\max_{i}\frac{\log p_{i}}{\log r_{i}})$,
 we therefore conclude that:
 $$
 \align
 &{}\\
 \sigma_{\abs}(\,\zeta_{\alpha}^{\ssim}\,)
 &=
  -\infty
  <
  0
  <
  f_{\mu}(\alpha)=f_{\mu}^{\con}(\alpha)
  \,\,\,\,\\
  &{}\\
 &\qquad
  \text{
  for all 
  except at most countably many 
  $\alpha\in  (\min_{i}\tfrac{\log p_{i}}{\log r_{i}},\max_{i}\tfrac{\log p_{i}}{\log r_{i}})$.
  }\\
  &{}
 \tag1.13
   \endalign
  $$

It follows from the above discussion that
while 
the definition
of $\zeta_{\alpha}^{\ssim}(s)$
is \lq\lq natural",
it is not 
does not encode
sufficient
information
allowing us to recover the multifractal spectra
$f_{\mu}(\alpha)$ and $f_{\mu}^{\con}(\alpha)$.
The reason for the strict inequality in (1.13) 
is, of course, clear:
%Namely,
%the problem 
%for 
%$\alpha\in
% (\min_{i}\frac{\log p_{i}}{\log r_{i}},\max_{i}\frac{\log p_{i}}{\log r_{i}}) \setminus
% \{\frac{\log p_{\bold i}}{\log r_{\bold i}}
% \,|\,
% \bold i\in\Sigma^{*}
% \}$
%is
%% of course,
%that
even though
there are no
strings
$\bold i\in\Sigma^{*}$
for which
the ratio
$\frac{\log p_{\bold i}}{\log r_{\bold i}}$ equals $\alpha$
if
$\alpha\in
 (\min_{i}\frac{\log p_{i}}{\log r_{i}},\max_{i}\frac{\log p_{i}}{\log r_{i}}) \setminus
 \{\frac{\log p_{\bold i}}{\log r_{\bold i}}
 \,|\,
 \bold i\in\Sigma^{*}
 \}$,
there are
nevertheless
many sequences $(\bold i_{n})_{n}$ of strings
$\bold i_{n}\in\Sigma^{*}$
for which the
ratios 
$\frac{\log p_{\bold i_{n}}}{\log r_{\bold i_{n}}}$ converges to $\alpha$.
In order to capture this,
it is necessary
to 
ensure
that
those
 strings
$\bold i$
for which
the ratio
$\frac{\log p_{\bold i}}{\log r_{\bold i}}$ 
is 
\lq\lq close"
 to $\alpha$
are
also included in the series
defining 
 the multifractal zeta-function.
For this reason, we
modify the definition of
$\zeta_{\alpha}^{\ssim}$
and
introduce
a
self-similar multifractal zeta-function
obtained by
replacing the 
original 
small \lq\lq target" set
$\{\alpha\}$
by a larger
 \lq\lq target" set
$I$
(for example,
 we may
choose the enlarged
 \lq\lq target" set
$I$
to be a non-degenerate  interval).
In order to
make this idea precise we proceed as follows.
 For a closed interval
$I$,
we 
define the self-similar multifractal zeta-function $\zeta^{\ssim}_{I}$ by
 $$
 \zeta^{\ssim}_{I}(s)
 =
 \sum
 \Sb
 \bold i\\
 {}\\
 \frac{\log p_{\bold i}}{\log r_{\bold i}}\in I
 \endSb
 r_{\bold i}^s
 \tag1.14
 $$
for those complex numbers $s$ for which the series converges absolutely.
Observe that if $I=\{\alpha\}$, then
 $$
  \zeta^{\ssim}_{I}(s)
=
 \zeta^{\ssim}_{\alpha}(s)\,.
 $$
% Using the
% zeta-function
% $ \zeta^{\ssim}_{I}$
% associated with an enlarged
% \lq\lq target" set $I$,
% we can now proceed in two equally natural ways.
We can now 
proceed in two equally natural ways.
Either, we can 
%replace the 
%original 
%small \lq\lq target" set
%$\{\alpha\}$
%by 
consider
a family
of enlarged
\lq\lq target" sets 
%$I$
shrinking to the original 
main \lq\lq target" $\{\alpha\}$;
this approach will be referred to as the
shrinking target approach.
Or, alternatively,
we can consider
a
fixed enlarge 
\lq\lq target" set
% $I$
and regard this as our
original main \lq\lq target";
this approach will be referred to as the
fixed target approach.
%Firstly, one 
%can attempt to let the 
%enlarged 
%\lq\lq target" set $I$ shrink 
%to the originally 
%small
%\lq\lq target" set $\{\alpha\}$.
%Secondly,
%one could choose to keep the
%enlarged
% \lq\lq target"
%set $I$ fixd and unchanged.
%Both approaches
%yield natural results.
We now
discuss these approaches in more detail.

\bigskip

\noindent
{\it (1) The shrinking target approach.
%First approach: 
%%Letting 
%%the 
%%enlarged 
%%\lq\lq target" set $I$ shrink 
%%to  $\{\alpha\}$.}
%Shrinking
%the 
%enlarged 
%\lq\lq target" set $I$  
%to  $\{\alpha\}$.
}
For a given (small) \lq\lq target"
$\{\alpha\}$,
we
consider
the following family
$\big(\,[\alpha-r,\alpha+r]\,\big)_{r>0}$
of enlarged
\lq\lq target" sets $[\alpha-r,\alpha+r]$
shrinking  to the original 
main \lq\lq target" $\{\alpha\}$ as $r\searrow0$,
and 
%try 
attempt
to
relate the limiting behaviour of the 
abscissa convergence
of 
$\zeta_{[\alpha-r,\alpha+r]}^{\ssim}$
to the multifractal spectrum
$f_{\mu}(\alpha)$ at $\alpha$.
In order to make this idea formal we proceed as follows.
For 
each 
$\alpha\in\Bbb R$ and for each $r>0$,
we define the zeta-function
$\zeta_{\alpha}^{\ssim}(\cdot;r)$ by
 $$
 \align
 \zeta^{\ssim}_{\alpha}(s;r)
&=
 \zeta_{[\alpha-r,\alpha+r]}^{\ssim}(s)\\
 &{}\\
&= 
 \sum
 \Sb
 \bold i\\
 {}\\
 \frac{\log p_{\bold i}}{\log r_{\bold i}}\in [\alpha-r,\alpha+r]
 \endSb
 r_{\bold i}^s\\ 
&{}\\ 
&= 
 \sum
 \Sb
 \bold i\\
 {}\\
 \big|\frac{\log p_{\bold i}}{\log r_{\bold i}}-\alpha\big|\le r
 \endSb
 r_{\bold i}^s\,.
 \tag1.15
 \endalign
 $$
The next result,
which is an application of one of our main results (see Theorem 3.6),
shows that 
the multifractal zeta-functions
$ \zeta^{\ssim}_{\alpha}(\cdot;r)$
encode sufficient
information allowing us
to recover 
the multifractal spectra
$f_{\mu}(\alpha)$ and $f_{\mu}^{\con}(\alpha)$
by letting 
$r\searrow 0$.

\bigskip

\proclaim{Theorem 1.1. Shrinking targets}
Assume that the list $(S_{1},\ldots S_{N})$
satisfies the OSC
and
let $\mu$ be the self-similar measure defined by
(1.9).
For  $\alpha\in\Bbb R$ and $r>0$,
let $\zeta_{\alpha}^{\ssim}(\cdot;r)$ be defined by (1.15).
Then we have
 $$
 \lim_{r\searrow 0}
 \sigma_{\abs}\big(\,\zeta_{\alpha}^{\ssim}(\cdot;r)\,\big)
 =
 f_{\mu}(\alpha)
 =
 f_{\mu}^{\con}(\alpha)
 $$
where  
$\sigma_{\abs}\big(\,\zeta_{\alpha}^{\ssim}(\cdot;r)\,\big)$
denotes the abscissa of convergence of the 
zeta-function
$\zeta_{\alpha}^{\ssim}(\cdot;r)$.
\endproclaim

\bigskip

\noindent
{\it (2) The fixed target approach
%Second approach: 
%%Letting the 
%%enlarged 
%%\lq\lq target" set $I$ be fixed.}
%Fixing
% the 
%enlarged 
%\lq\lq target" set $I$.
}
Alternatively
we
can keep the
enlarged 
\lq\lq target" set $I$
fixed
and 
%try
attempt
to relate the
abscissa of convergence
of the
multifractal zeta-function $\zeta_{I}^{\ssim}$
associated with the enlarger \lq\lq target" set $I$
to the 
values
 of the multifractal spectrum
$f_{\mu}(\alpha)$ 
for $\alpha\in I$.
Of course, inequality(1.13) shows that
if the \lq\lq target" 
set $I$ is 
\lq\lq too small",
then this is not possible.
However,
if
the enlarger \lq\lq target" set $I$ 
satisfies 
a mild
non-degeneracy
condition, namely condition (1.16), 
guaranteeing
that $I$ is 
sufficiently
\lq\lq big",
then
the next result,
which is also 
an application of one of our main results (see Theorem 3.6),
shows that 
this is possible. 
More
precisely the result
shows that
if
the enlarger \lq\lq target" set $I$ 
satisfies condition 
 (1.16),
then
the multifractal zeta-function
$\zeta_{I}^{\ssim}$
associated with the enlarger \lq\lq target" set $I$
encode sufficient
information allowing us
to recover 
the suprema
$\sup_{\alpha\in I}
 f_{\mu}(\alpha)$
 and
$\sup_{\alpha\in I}
 f_{\mu}^{\con}(\alpha)$
 of the 
 multifractal spectra $ f_{\mu}(\alpha)$ and $ f_{\mu}^{\con}(\alpha)$
for $\alpha\in I$.

\bigskip

\proclaim{Theorem 1.2. Fixed targets}
Assume that the list $(S_{1},\ldots S_{N})$
satisfies the OSC
and
let $\mu$ be the self-similar measure defined by
(1.9).
For  a closed interval $I$,
let $\zeta_{I}^{\ssim}$ be defined by (1.14).
If 
 $$
 \overset{\,\circ}\to{I}
 \cap
 \Big(\min_{i}\tfrac{\log p_{i}}{\log r_{i}},\max_{i}\tfrac{\log p_{i}}{\log r_{i}}\Big)
 \not=
 \varnothing
 \tag1.16
 $$
(where $\overset{\,\circ}\to{I}$ denotes the interior of $I$),
then we have
 $$
 \sigma_{\abs}\big(\,\zeta_{I}^{\ssim}\,\big)
 =
 \sup_{\alpha\in I}
 f_{\mu}(\alpha)
 =
 \sup_{\alpha\in I}
 f_{\mu}^{\con}(\alpha)
 $$
 where  
$\sigma_{\abs}\big(\,\zeta_{I}^{\ssim}\,\big)$
denotes the abscissa of convergence of the 
zeta-function
$\zeta_{I}^{\ssim}$.
 \endproclaim

\bigskip

We emphasise that Theorem 1.1 and Theorem 1.2
are
presented
in order to motive this work
and
are special cases
of the
substantially
more general and abstract 
theory 
of multifractal zeta-function
developed in this paper.

The next section, i\.e\. Section 2, describes the 
general framework
developed in this paper and list our main results.
In Section 3 we will discuss a number of examples, including,
mixed and non-mixed
multifractal spectra 
of self-similar and self-conformal measures,
and multifractal spectra of Birkhoff ergodic averages.

\bigskip
\bigskip

\heading{2. Statements of main results.}\endheading

{\bf 2.1.
Main definitions:
the zeta-functions 
 $\zeta_{C}^{U,\Lambda}(\cdot)$
 and 
 $\zeta_{C}^{U,\Lambda}(\cdot;r)$.}
In this section we
describe
the
 framework
developed in this paper and list our main results.
We first recall and introduce some useful notation.
Fix a positive integer $N$.
 Let
$\Sigma
 =
 \{1,\ldots,N\}$ and 
 for a positive integer $n$,
 write
 $$
 \aligned
 \Sigma^{n}
&=\{1,\dots,N\}^{n}\,,\\
 \Sigma^{*}
&=\bigcup_{m}\Sigma^{m}\,,
 \endaligned
 $$
i\.e\. $\Sigma^{n}$ is the family of all
strings
$\bold i=i_{1}\ldots i_{n}$
of length $n$ 
with $i_{j}\in\{1,\ldots,N\}$
and
$\Sigma^{*}$ is the family of all finite strings
$\bold i=i_{1}\ldots i_{m}$
with $m\in\Bbb N$ and $i_{j}\in\{1,\ldots,N\}$.
Also write
$$
 \Sigma^{\Bbb N}
=\{1,\dots,N\}^{\Bbb N}\,,
 $$
i\.e\. $\Sigma^{\Bbb N}$ is the family of all
infinite
strings
$\bold i=i_{1}i_{2}\ldots $
with $i_{j}\in\{1,\ldots,N\}$.
%Next,
For an infinite string 
$\bold i=i_{1}i_{2}\ldots\in\Sigma^{\Bbb N} $
and a positive integer $n$, we will write
$\bold i|n
=
i_{1}\ldots i_{n}$.
In addition, for
a positive integer $n$
and
a finite string 
$\bold i=i_{1}\ldots i_{n}\in\Sigma^{n} $
with length equal to  $n$,
 we will write
$|\bold i|
=
n$, and we let $[\bold i]$ denote the cylinder 
generated by $\bold i$, i\.e\.
  $$
  [\bold i]
 =
 \Big\{
 \bold j\in\Sigma^{\Bbb N}
 \,\Big|\,
 \bold j|n=\bold i
 \Big\}\,.
 $$
Also, let $S:\Sigma^{\Bbb N}\to\Sigma^{\Bbb N}$ denote the shift map.
Finally, we  denote the family of Borel probability measures on 
$\Sigma^{\Bbb N}$
by $\Cal P(\Sigma^{\Bbb N})$
and 
we
equip
$\Cal P(\Sigma^{\Bbb N})$  with the weak topology.

The multifractal zeta-function
framework
developed 
 in this paper
depend on a space $X$ and two maps $U$ and $\Lambda$
satisfying various conditions.
We will now
introduce the space $X$ and the maps $U$ and $\Lambda$.

\medskip

\roster
\item"(1)"
First, we fix a metric space $X$.

\medskip

\item"(2)"
Next, we fix a continuous map
$U:\Cal P(\Sigma^{\Bbb N})\to X$.

\medskip

\item"(3)"
Finally,
we fix a function
$\Lambda:\Sigma^{\Bbb N}\to\Bbb R$
satisfying the 
following three conditions:
%\roster
%\item"(1)"
% {\it Continuity:}

\smallskip

(C1)
The function $\Lambda$ is continuous;
%\item"(2)"
% {\it Hyperbolicity:}

\smallskip

(C2)
There are constants $c_{\min}$ and $c_{\max}$ with
$-\infty<c_{\min}\le c_{\max}<0$ such that
 $c_{\min}
 \le
\Lambda
\le$

\quad\,\,\,\,
\,\,\,\,
$c_{\max}$;
%\item"(3)"
% {\it Bounded Distortion:}

\smallskip

(C3)
There is a constant
$c$ with $c\ge 1$
such that
for all positive integers $n$ and all $\bold i,\bold j\in\Sigma^{\Bbb N}$

\quad\,\,\,\,
\,\,\,\,
with
$\bold i|n=\bold j|n$, we have
 $$
 \frac{1}{c}
 \le
 \frac
 {\dsize\,\,\exp\,\,\sum_{k=0}^{n-1}\Lambda S^{k}\bold i\,\,}
 {\dsize\,\,\exp\,\,\sum_{k=0}^{n-1}\Lambda S^{k}\bold j\,\,}
 \le
 c\,.
 $$
%\endroster
Condition (C2) is clearly motivated by the hyperbolicity condition
from
dynamical systems,
and Condition (C3) is equally clearly motivated the bounded distortion property 
from
dynamical systems. 

\endroster

 \medskip

 \noindent
 Associated with the space $X$ and the maps $U$ and $\Lambda$, we now define the 
 following multifractal zeta-functions.

 \bigskip
 
 \proclaim{Definition.
 The zeta-functions 
 $\zeta_{C}^{U,\Lambda}$
 and 
 $\zeta_{C}^{U,\Lambda}(\cdot;r)$
associated with the space $X$ and the maps $U$ and $\Lambda$ }
 For a finite string $\bold i\in\Sigma^{n}$, let
 $$
 s_{\bold i}
 =
 \sup_{\bold k\in[\bold i]}
 \,\,
 \exp
 \,\,
 \sum_{k=0}^{n-1}\Lambda S^{k}\bold k
 \,,
 $$
and for a positive integer $n$ and an 
infinite string  $\bold i\in\Sigma^{\Bbb N}$,
let
$L_{n}:\Sigma^{\Bbb N}\to\Cal P(\Sigma^{\Bbb N})$ be defined 
by
 $$
 L_{n}\bold i
 =
 \frac{1}{n}\sum_{k=0}^{n-1}\delta_{S^{k}\bold i}\,.
 $$
 For $C\subseteq X$,
we define the zeta-function $\zeta_{C}^{U,\Lambda}$
associated with the space $X$ and the maps $U$ and $\Lambda$ by
 $$
  \zeta_{C}^{U,\Lambda}(s)
  =
   \sum
  \Sb
  \bold i\\
  {}\\
  UL_{|\bold i|}[\bold i]\subseteq C
  \endSb
 s_{\bold i}^{s}
 $$
for those complex numbers $s$ for which the series converges absolutely,
and
for $r>0$ and $C\subseteq X$,
we define the zeta-function $\zeta_{C}^{U,\Lambda}(\cdot;r)$
associated with the space $X$ and the maps $U$ and $\Lambda$ by
 $$
 \align
  \zeta_{C}^{U,\Lambda}(s;r)
 &=
 \zeta_{B(C,r)}^{U,\Lambda}(s)\\
 &=
   \sum
  \Sb
  \bold i\\
  {}\\
  UL_{|\bold i|}[\bold i]\subseteq B(C,r)
  \endSb
 s_{\bold i}^{s}
 \endalign
 $$
for those complex numbers $s$ for which the series converges absolutely
and
where
$B(C,r)
 =
 \{
 x\in X
 \,|\,
 \dist(x,C)\le r
 \}$
denotes the closed $r$ neighborhood of $C$.

\endproclaim

\bigskip

\noindent
Next, we formally define the 
abscissa of convergence
(of a zeta-function).

\bigskip

\proclaim{Definition.
Abscissa of convergence}
Let $(\,a_{\bold i}\,)_{\bold i\in\Sigma^{*}}$ be a family of positive numbers and define the
(zeta-)function $\zeta$ by
 $$
 \zeta(s)
 =
 \sum_{\bold i}a_{\bold i}^{s}
 $$
for those complex numbers $s$
for which the series converges.  
The abscissa of convergence
of $\zeta$ is defined by
 $$
 \sigma_{\abs}(\zeta)
 =
 \inf
 \Bigg\{
 t\in\Bbb R
 \,\Big|\,
 \text{the series}
 \,\,
  \sum_{\bold i}a_{\bold i}^{t}
  \,\,
 \text{converges absolutely}
 \Bigg\}\,. 
 $$
\endproclaim

\bigskip

Our main results, i\.e\. Theorem 2.1 and Theorem 2.2 below,
 relate the abscissa of converge of the zeta-functions
$ \zeta_{C}^{U,\Lambda}(\cdot;r)$
and
$ \zeta_{C}^{U,\Lambda}$
to various multifractal 
quantities, including,
the
coarse multifractal spectrum associated
with the space $X$ and the maps $U$ and $\Lambda$.
In order  to state Theorem 2.1 and Theorem 2.2
 we will
 now
 define the
coarse multifractal spectra.

\bigskip

  \proclaim{Definition.
 The coarse multifractal spectra
associated with the space $X$ and the maps $U$ and $\Lambda$ }
 For $\bold i=i_{1}\ldots i_{n}\in\Sigma^{*}$,
 we let
$\widehat{\bold i}=i_{1}\ldots i_{n-1}\in\Sigma^{*}$ 
denote the 
\lq\lq parent" of $\bold i$.
Next,
for $\bold i\in\Sigma^{*}$
 and $\delta>0$, we write
  $$
  s_{\bold i}
  \approx
  \delta
  $$
 if and only if
 $s_{\bold i}
  \le
  \delta
  <
  s_{\widehat{\bold i}}$.
 For $r>0$ and $C\subseteq X$, let
   $$
  \Pi_{\delta}^{U,\Lambda}(C,r)
  =
  \Big\{
  \bold i
  \,\Big|\,
  s_{\bold i}
  \approx
  \delta\,,\,
  UL_{|\bold i|}[\bold i]\subseteq B(C,r)
 \Big\}
  $$ 
 and
  $$
   N_{\delta}^{U,\Lambda}(C,r) 
 =
 \Big|
 \,
 \Pi_{\delta}^{U,\Lambda}(C,r)
 \,
 \Big|\,.
 $$
We define the lower and upper $r$-approximate
coarse multifractal spectrum associated
with the space $X$ and the maps $U$ and $\Lambda$
by
 $$
 \align
 \underline f^{U,\Lambda}(C,r)
&=
\,
\liminf_{\delta\searrow 0}
\,
\frac{\log  N_{\delta}^{U,\Lambda}(C,r)}{-\log \delta}\,,\\
\overline f^{U,\Lambda}(C,r)
&=
\limsup_{\delta\searrow 0}
\frac{\log  N_{\delta}^{U,\Lambda}(C,r)}{-\log \delta}\,,
\endalign
$$
and
we define the lower and upper 
coarse multifractal spectrum associated
with the space $X$ and the maps $U$ and $\Lambda$
by
 $$
 \align
 \underline f^{U,\Lambda}(C)
&=
\lim_{r\searrow 0}
\underline f^{U,\Lambda}(C,r)\,,\\
 \overline f^{U,\Lambda}(C)
&=
\lim_{r\searrow 0}
\overline f^{U,\Lambda}(C,r)\,.
\endalign
$$
\endproclaim

\bigskip

Below we state our main results. 
As suggested by the 
discussion in Section 1.4,
we will attempt to relate 
the abscissae of 
convergence of the
multifractal 
zeta-functions
$\zeta_{C}^{U,\Lambda}$
 and 
 $\zeta_{C}^{U,\Lambda}(\cdot;r)$
to
various multifractal spectra 
using two different but equally natural approaches:
the shrinking target approach or  the fixed target approach.
 The shrinking target approach is 
 discussed in Section 2.2 and the fixed target approach is discussed in Section 2.3.

\bigskip

{\bf 2.2.
First main result.
The shrinking target approach:
finding
 $\lim_{r\searrow 0}\sigma_{\abs}\big(\,\zeta_{C}^{U,\Lambda}(\cdot;r)\,\big)$.}
 For a given  \lq\lq target"
$C$,
we
consider
the following family
$\big(\,B(C,r)\,\big)_{r>0}$
of enlarged
\lq\lq target" sets $B(C,r)$
shrinking to the original 
main \lq\lq target" $C$ as $r\searrow0$,
and 
%try 
attempt
to
relate the limiting behaviour of the 
abscissa convergence
of 
the zeta-function
$\zeta_{C}^{U,\Lambda}(\cdot;r)=\zeta_{B(C,r)}^{U,\Lambda}$
to the coarse multifractal spectrum
$\underline f^{U,\Lambda}(C)$ 
and other multifractal quatities.
Our first main result,
i\.e\. Theorem 2.1 below,
shows that this is possible.
More precisely, Theorem 2.1 shows 
that
the 
abscissa of convergence
% $\sigma_{\abs}\big(\,\zeta_{C}^{U,\Lambda}(\cdot;r)\,\big)$
of the 
zeta-function
$ \zeta_{C}^{U,\Lambda}(\cdot;r)$
converges
as $r\searrow0$,
and that this limit 
%can be computed in
%the following  two ways.
%%Namely,
%Firstly, we show that the 
%limit
equals
the 
%natural
coarse multifractal spectrum  of $C$.
%Secondly,
%we show that
We also show that
the limit
can be obtained 
by
a variational principle 
involving the
supremum 
of the 
entropy
of all shift invariant Borel probability
measures
$\mu\in\Cal P(\Sigma^{\Bbb N})$
with $U\mu\in C$.
In Section 3
we show that
in many important 
cases
%examples 
the limit
$\lim_{r\searrow0}\sigma_{\abs}\big(\,\zeta_{C}^{U,\Lambda}(\cdot;r)\,\big)$ equals 
the traditional multifractal spectra.

 \bigskip

\proclaim{Theorem 2.1. Shrinking targets}
Let $X$ be a metric space and let $U:\Cal P(\Sigma^{\Bbb N})\to X$ be 
continuous with respect to the weak topology.
Let $C\subseteq X$ be a closed subset of $X$.
%
%We have
% $$
% \align
%  \lim_{r\searrow 0}
% \,\,
% \sigma_{\abs}\big(\,\zeta_{C}^{U,\Lambda}(\cdot;r)\,\big)
%&=
% \underline f^{U,\Lambda}(C)\,,\\ 
%&{}\\
% \lim_{r\searrow 0}
% \,\,
% \sigma_{\abs}\big(\,\zeta_{C}^{U,\Lambda}(\cdot;r)\,\big)
%&=
%\,\,\,\,
% \sup
%  \Sb
%  \mu\in\Cal P_{S}(\Sigma^{\Bbb N})\\
%  {}\\
%  U\mu\in C
%  \endSb
%  \,\,\,\,
% -\frac{h(\mu)}{\int \Lambda\,d\mu}\,.
% \endalign
% $$
%
\roster
\item"(1)"
The lower 
coarse multifractal spectrum associated
with the space $X$ and the maps $U$ and $\Lambda$:
we have
 $$
 \align
  \lim_{r\searrow 0}
 \,\,
 \sigma_{\abs}\big(\,\zeta_{C}^{U,\Lambda}(\cdot;r)\,\big)
&=
 \underline f^{U,\Lambda}(C)\,.
 \qquad\qquad
 \quad\,\,\\ 
%&{}\\
% \lim_{r\searrow 0}
% \,\,
% \sigma_{\abs}\big(\,\zeta_{C}^{U,\Lambda}(\cdot;r)\,\big)
%&=
%\,\,\,\,
% \sup
%  \Sb
%  \mu\in\Cal P_{S}(\Sigma^{\Bbb N})\\
%  {}\\
%  U\mu\in C
%  \endSb
%  \,\,\,\,
% -\frac{h(\mu)}{\int \Lambda\,d\mu}\,.
 \endalign
 $$
\item"(2)"
The variational principle:
we have 
 $$
 \align
%  \lim_{r\searrow 0}
% \,\,
% \sigma_{\abs}\big(\,\zeta_{C}^{U,\Lambda}(\cdot;r)\,\big)
%&=
% \underline f^{U,\Lambda}(C)\,,\\ 
%&{}\\
 \lim_{r\searrow 0}
 \,\,
 \sigma_{\abs}\big(\,\zeta_{C}^{U,\Lambda}(\cdot;r)\,\big)
&=
\,\,\,\,
 \sup
  \Sb
  \mu\in\Cal P_{S}(\Sigma^{\Bbb N})\\
  {}\\
  U\mu\in C
  \endSb
  \,\,\,\,
 -\frac{h(\mu)}{\int \Lambda\,d\mu}\,;
 \endalign
 $$
here $\Cal P_{S}(\Sigma^{\Bbb N})$ denotes the family of shift invariant Borel probability
measures on $\Sigma^{\Bbb N}$ 
and 
$h(\mu)$ denotes the entropy of $\mu\in \Cal P_{S}(\Sigma^{\Bbb N})$.
\endroster 
\endproclaim
 
 \bigskip

 \noindent
 In order to prove 
 Theorem 2.1 it suffices 
 to prove the following three inequalities:
  $$
  \align
   \limsup_{r\searrow 0}
 \,\,
 \sigma_{\abs}\big(\,\zeta_{C}^{U,\Lambda}(\cdot;r)\,\big)
&\le
\,\,\,\,
 \sup
  \Sb
  \mu\in\Cal P_{S}(\Sigma^{\Bbb N})\\
  {}\\
  U\mu\in C
  \endSb
  \,\,\,\,
 -\frac{h(\mu)}{\int \Lambda\,d\mu}\,,
 \tag2.1\\
&{}\\
  \sup
  \Sb
  \mu\in\Cal P_{S}(\Sigma^{\Bbb N})\\
  {}\\
  U\mu\in C
  \endSb
  \,\,\,\,
 -\frac{h(\mu)}{\int \Lambda\,d\mu}
 &\le
  \underline f^{U,\Lambda}(C)\,,
  \tag2.2\\
  &{}\\
  \underline f^{U,\Lambda}(C)
 &\le
  \liminf_{r\searrow 0}
 \,\,
 \sigma_{\abs}\big(\,\zeta_{C}^{U,\Lambda}(\cdot;r)\,\big)\,.
 \tag2.3
  \endalign
  $$ 
Inequality (2.1) is proven in 
Section 5 using techniques from 
the theory of large deviations.
Inequality (2.2) is proven in 
Section 6 using techniques from 
ergodic theory.
Finally,
inequality (2.3) follows 
directly
from the definitions and is proved in 
Section 7.

\bigskip

{\bf 2.3.
Second main result.
The fixed
target approach:
finding
 $\sigma_{\abs}\big(\,\zeta_{C}^{U,\Lambda}\,\big)$.}
 Alternatively,
instead of choosing a family of \lq\lq target" sets that shrinks to the given 
\lq\lq target" $C$,
we
can keep the
given
\lq\lq target" set $C$
fixed
and 
%try
attempt
to relate the
abscissa of convergence
of the
multifractal zeta-function $ \zeta_{C}^{U,\Lambda}$
associated with the \lq\lq target" set $C$
to the 
values
 of the multifractal spectrum
coarse multifractal spectrum
$\underline f^{U,\Lambda}(C)$.
Of course, the example in Section 1.4
 shows that
if the \lq\lq target" 
set $C$ is 
\lq\lq too small",
then this is not possible.
However,
if
the coarse multifractal spectrum
$\underline f^{U,\Lambda}$
satisfies a continuity condition at $C$
guaranteeing
that the interior of $C$ is 
\lq\lq sufficiently big",
then
our second main result,
i\.e\. Theorem 2.2 below,
shows that
this is possible. 
More precisely, Theorem 2.2 shows that
if
the coarse multifractal spectrum
$\underline f^{U,\Lambda}$
is inner continuous  at $C$
(the definition of inner continuity will be given below),
then
the 
abscissa of convergence
of the 
zeta-function
$ \zeta_{C}^{U,\Lambda}$
equals
the 
coarse multifractal spectrum of $C$.
In analogy with Theorem 2.1,
we also show that
the 
abscissa of convergence
of
$\zeta_{C}^{U,\Lambda}$
can be obtained 
by
a variational principle 
involving the
supremum 
of the 
entropy
of all shift invariant Borel probability
measures
$\mu\in\Cal P(\Sigma^{\Bbb N})$
with $U\mu\in C$.
However,
before stating Theorem 2.2,
we first define the
continuity condition that
the coarse multifractal spectrum
$\underline f^{U,\Lambda}$
is required to satisfy.

\bigskip

\proclaim{Definition. 
Inner continuity}
Let
$P(X)$ denote the family of subsets of $X$
and
for $C\subseteq X$ and $r>0$, write
 $$
 I(C,r)
% =
% \overline{
% C
% \setminus
% B(\partial C,r)}\,.
  =
 \Big\{
 x\in C
 \,\Big|\,
 \dist(x,\partial C)\ge r
 \Big\}\,.
 $$
% \roster
% \item"(1)"
We say that a function $\Phi:P(X)\to\Bbb R$ 
is inner continuous at $C\subseteq X$ if
 $$
 \Phi\big(\,I(C,r)\,\big)
 \to
 \Phi(C)
 \,\,\,\,
 \text{as $r\searrow0$}\,.
$$
%\item"(2)"
%We
%say that a function
%$f:X\to\Bbb R$
%is inner continuous at $C\subseteq X$ if
%the function
%  $$
%  \align
%  P(X)
% &\to\Bbb R\\
% E
% &\to
% \Phi\{f(x)\,|\, x\in E\}
% \endalign
% $$
%is inner continuous at $C\subseteq X$.
%\endroster
\endproclaim

\bigskip

\noindent
We can now state Theorem 2.2.

 \bigskip

\proclaim{Theorem 2.2. Fixed targets}
Fix a positive integer $M$.
Let $U:\Cal P(\Sigma^{\Bbb N})\to \Bbb R^{M}$ be 
continuous with respect to the weak topology.
Let $C\subseteq \Bbb R^{M}$ be a closed subset of $\Bbb R^{M}$
and assume that 
$ \underline f^{U,\Lambda}$
is inner continuous at $C$.

\roster
\item"(1)"
The lower 
coarse multifractal spectrum associated
with $\Bbb R^{M}$ and the maps $U$ and $\Lambda$:
we have
 $$
 \align
 \sigma_{\abs}\big(\,\zeta_{C}^{U,\Lambda}\,\big)
&=
 \underline f^{U,\Lambda}(C)\,.
 \qquad\qquad
 \quad\,\,\\ 
 \endalign
 $$
\item"(2)"
The variational principle:
we have 
 $$
 \align
 \sigma_{\abs}\big(\,\zeta_{C}^{U,\Lambda}\,\big)
&=
\,\,\,\,
 \sup
  \Sb
  \mu\in\Cal P_{S}(\Sigma^{\Bbb N})\\
  {}\\
  U\mu\in C
  \endSb
  \,\,\,\,
 -\frac{h(\mu)}{\int \Lambda\,d\mu}\,;
 \endalign
 $$
here $\Cal P_{S}(\Sigma^{\Bbb N})$ denotes the family of shift invariant Borel probability
measures on $\Sigma^{\Bbb N}$ 
and 
$h(\mu)$ denotes the entropy of $\mu\in \Cal P_{S}(\Sigma^{\Bbb N})$. 
\endroster 
\endproclaim

 \bigskip
 
 \noindent
Theorem 2.2 follows easily from Theorem 2.1
and is proved in Section 8.

\bigskip

{\bf 2.4. Euler product.}
We will now prove that the 
multifractal zeta-function $\zeta_{C}^{U,\Lambda}$ has
a natural Euler product.
We begin with a definition.

\bigskip

\proclaim{Definition. Composite and prime}
A finite string $\bold i\in\Sigma^{*}$ is called composite (or peiodic)
if there is $\bold u\in\Sigma^{*}$ and a positive integer $n>1$ such that
$\bold i
=
%\underbrace{\bold u\ldots\bold u}_{\text{$n$ times}}
\bold u\ldots\bold u$
where $\bold u$ is repeated $n$ times.
A finite string $\bold i\in\Sigma^{*}$ is called prime if it is not composite.
\endproclaim

\bigskip

\noindent
Theorem 2.3 shows that
 $\zeta_{C}^{U,\Lambda}$ has
an Euler product.
In Theorem 2.3 we use the following notation, namely,
if $f$ is a holomorphic function that does not attain the value $0$, then we
let $Lf$ denote the logarithmic derivative of $f$, i\.e\.
$Lf=\frac{f'}{f}$.
%In addition, for complex numbers $u,w\in \Bbb C$ with $u\not=0$
%and
%$|\frac{u-1}{u}|<1$, we
%define $u^{w}$ by
%$u^{w}
%=
%\exp(w\sum_{n}\frac{1}{n}(\frac{u-1}{u})^{n})$;
%note that this definition is consistent with
%the usual definition of $u^{w}$ for 
%$u\in \Bbb R$ with $u\not=0$
%and
%$|\frac{u-1}{u}|<1$ --
%indeed,
%for 
%$u\in \Bbb R$ with $u\not=0$
%and
%$|\frac{u-1}{u}|<1$, we have
%$\exp(w\sum_{n}\frac{1}{n}(\frac{u-1}{u})^{n})
%=
%\exp(w\log(1-(\frac{u-1}{u}))^{-1})
%=
%\exp(w\log u)
%=
%u^{w}$.
We can now state Theorem 2.3.

\bigskip

\proclaim{Theorem 2.3. Euler product}
Let $X$ be a metric space and let 
$U:\Cal P(\Sigma^{\Bbb N})\to X$ be continuous
with respect to the weak topology. 
Assume that
 $$
 s_{\bold i\bold j}
 =
 s_{\bold i}s_{\bold j}
 $$
for all $\bold i,\bold j\in\Sigma^{*}$.
Let 
$C\subseteq X$ be a closed subset of $X$.
 \roster
 \item"(1)"
For complex numbers $s$
with
$\Real(s)>\sigma_{\abs}(\,\zeta_{C}^{U,\Lambda}\,)$,
the product
$$
Q_{C}^{U,\Lambda}(s)
=
\prod
\Sb
\bold i\\
{}\\
\text{$\bold i$ is prime}\\
{}\\
UL_{|\bold i|}[\bold i]\subseteq C
\endSb
\Bigg(
\frac{1}{1-s_{\bold i}^{s}}
\Bigg)^{\frac{1}{\log s_{\bold i}}}
$$
converges
and 
$Q_{C}^{U,\Lambda}(s)\not=0$.
The product
$Q_{C}^{U,\Lambda}(s)$ is called the
Euler product of
$\zeta_{C}^{U,\Lambda}$.
\item"(2)"
For all
complex numbers $s$
with
$\Real(s)>\sigma_{\abs}(\,\zeta_{C}^{U,\Lambda}\,)$, we have
  $$
  \zeta_{C}^{U,\Lambda}(s)
  =
   L\, Q_{C}^{U,\Lambda}(s) \,.
 $$
 \endroster
\endproclaim

\bigskip

\noindent
Theorem 2.3 is proved in Section 9.

  \bigskip
  \bigskip

 %%%%%%%%%%%%%%%%%%%%%%%%%%%%%%%%%%%%%%%%%
 %%%%%%%%%%%%%%%%%%%%%%%%%%%%%%%%%%%%%%%%%
 %%%%%%%%%%%%%%%%%%%%%%%%%%%%%%%%%%%%%%%%%
 %%%%%%%%%%%%%%%%%%%%%%%%%%%%%%%%%%%%%%%%%
 %%%%%%%%%%%%%%%%%%%%%%%%%%%%%%%%%%%%%%%%% 

%  
%  
%\heading{3. Applications: 
%multifractal spectra of measures
%and 
%multifractal spectra of ergodic Birkhoff averages}
%\endheading
% 
% 

\centerline{\smc 3. Applications: }

\centerline{\smc multifractal spectra of measures}

\centerline{\smc  and}

\centerline{\smc  
multifractal spectra of ergodic Birkhoff averages}

 \medskip

We will now consider several of applications of Theorem 2.1 and Theorem 2.2
to multifractal spectra of measures and ergodic averages.
In particular, we consider the following examples:

\bigskip

$\bullet$
 Section 3.1: Multifractal spectra of self-conformal measures.

\bigskip

$\bullet$
 Section 3.2: Mixed multifractal spectra of self-conformal measures.

\bigskip

$\bullet$
Section 3.3: Multifractal spectra of self-similar  measures.

\bigskip

$\bullet$
Section 3.4: Multifractal spectra of ergodic Birkhoff averages.

\bigskip

{\bf 3.1.
Multifractal 
spectra of 
self-conformal measures.}
Since
our examples are formulated in the setting of self-conformal 
(or self-similar) measures
we begin be recalling the definition of self-conformal 
(and self-similar) measures. 
A conformal iterated function system with probabilities 
is a list
 $$
 \big(
 \,
 V
 \,,\,
 X
 \,,\,
 (S_{i})_{i=1,\ldots,N}
% \,,\,
% (p_{i})_{i=1,\ldots,N}
 \,
 \big)
 $$
where
\roster
\item"(1)"
$V$ is an open, connected subset of $\Bbb R^{d}$.
\item"(2)" 
$X$ is a compact set with $X\subseteq V$
and
$X^{\circ\,-}=X$.
\item"(3)" 
$S_{i}:V\to V$ is a contractive
$C^{1+\gamma}$ diffeomorphism with
$0<\gamma<1$ such that
 $S_{i}X
 \subseteq
 X$
for all $i$. 
\item"(4)"
The Conformality Condition:
For each $x\in V$, we have that
$(DS_{i})(x)$ is a contractive similarity map, i\.e\.
there exists
$r_{i}(x)\in(0,1)$ such that
$|(DS_{i})(x)u-(DS_{i})(x)v|
 =
 r_{i}(x)|u-v|$
for all $u,v\in\Bbb R^{d}$;
here $(DS_{i})(x)$ denotes the derivative of $S_{i}$ at $x$.
%\item"(5)" 
%$(p_{i})_{i}=1,\ldots,N$ is a probability vector.
\endroster

\bigskip

\noindent
It follows from [Hu] that there exists a unique
non-empty compact set $K$ with $K\subseteq X$
 such that
 $$
 K
 =
 \bigcup_{i}\,
 S_{i}K\,.
 \tag3.1
 $$
The set $K$
is called the self-conformal set
associated with the list 
$\big(
 \,
 V
 \,,\,
 X
 \,,\,
 (S_{i})_{i=1,\ldots,N}
% \,,\,
% (p_{i})_{i=1,\ldots,N}
 \,
 \big)$;
in particular,
 if each map $S_{i}$ is a contracting similarity, then the 
 set $K$ is called the 
 self-similar set
 associated with the list 
$\big(
 \,
 V
 \,,\,
 X
 \,,\,
 (S_{i})_{i=1,\ldots,N}
% \,,\,
% (p_{i})_{i=1,\ldots,N}
 \,
 \big)$.
In addition,
if
$(p_{i})_{i=1,\ldots,N}$ is a probability vector
then it follows from [Hu]
that there is a unique
probability measure $\mu$ with
$\supp\mu=K$ such that
 $$
 \mu
 =
 \sum_{i}\,
 p_{i}\,\mu\circ S_{i}^{-1}\,.
 \tag3.2
 $$
The measure $\mu$
is called the self-conformal measure
associated with the list
$\big(
 \,
 V
 \,,\,
 X
 \,,\,
 (S_{i})_{i=1,\ldots,N}
  \,,\,\allowmathbreak
 (p_{i})_{i=1,\ldots,N}
 \,
 \big)$;
 if each map $S_{i}$ is a contracting similarity,
 then the
 measure $\mu$
is called the self-similar measure
associated with the list
$\big(
 \,
 V
 \,,\,
 X
 \,,\,
 (S_{i})_{i=1,\ldots,N}
  \,,\,\allowmathbreak
 (p_{i})_{i=1,\ldots,N}
 \,
 \big)$.
We will frequently assume that the list
$\big(
 \,
 V
 \,,\,
 X
 \,,\,
 (S_{i})_{i=1,\ldots,N}\allowmathbreak
% \,,\,
% (p_{i})_{i=1,\ldots,N}
 \,
 \big)$
 satisfies
  the Open Set Condition defined below.
Namely,
the list 
$\big(
 \,
 V
 \,,\,
 X
 \,,\,
 (S_{i})_{i=1,\ldots,N}
 \,
 \big)$
satisfies the Open Set Condition (OSC) if
there exists an
open, non-empty and bounded set $O$
with $O\subseteq X$
and
$S_{i}O
 \subseteq
 O$
for all $i$ such that
$S_{i}O
 \,\cap\,
 S_{j}O
 =
 \varnothing$
for all $i,j$ with 
$i\not=j$.

 %\noindent
Next, we 
define 
the natural projection map $\pi:\Sigma^{\Bbb N}\to K$.
However, we first make the follwing definitions.
Namely,
for $\bold i=i_{1}\ldots i_{n}\in\Sigma^{*}$, write
 $$
\align
S_{\bold i}
&=
S_{i_{1}}\cdots S_{i_{n}}\,,\\
K_{\bold i}
&=
S_{\bold i}K\,.
\endalign
$$
The natural projection map $\pi:\Sigma^{\Bbb N}\to K$
is now defined 
 by
 $$
 \Big\{\,\pi(\bold i)\,\Big\}
 =
 \bigcap_{n}S_{\bold i|n}K
 $$
for $\bold i\in\Sigma^{\Bbb N}$.

%

%\noindent
%Finally, we will
%define the natural projection map $\pi:\Sigma^{\Bbb N}\to K$.
%For $\bold i=i_{1}\ldots i_{n}\in\Sigma^{*}$, we write
%$S_{\bold i}=S_{i_{1}}\cdots S_{i_{n}}$.
%Also, for 
%$\bold i=i_{1}i_{2}\ldots\in\Sigma^{\Bbb N}$
%and a positive integer $n$,
%we write
%$\bold i|n=i_{1}\ldots i_{n}$.
%We now
%define the natural projection map $\pi:\Sigma^{\Bbb N}\to K$
% by
% $$
% \Big\{\,\pi(\bold i)\,\Big\}
% =
% \bigcap_{n}S_{\bold i|n}K
% $$
%for $\bold i\in\Sigma^{\Bbb N}$. 
%Also, recall that $\Sigma^{\Bbb N}\to\Sigma^{\Bbb N}$ denotes the shift map.

Finally,  we collect the definitions and results from
multifractal analysis of self-conformal measures
that we need
in order to state our main results.
We first recall, that
the Hausdorff multifractal spectrum
$f_{\mu}$ of $\mu$ 
is defined
by
 $$
 \align
 f_{\mu}(\alpha)
&=
   \,\dim_{\Haus}
   \left\{x\in K
    \,\left|\,
     \lim_{r\searrow0}
     \frac
     {\log\mu B(x,r)}{\log r}
     =
     \alpha
      \right.
   \right\}\,,
  \endalign
 $$
for $\alpha\in\Bbb R$
where $\dim_{\Haus}$ 
denotes the Hausdorff dimension.
In the late 1990's
Patzschke [Pa],
building on works by
Cawley \& Mauldin [CaMa]  and Arbeiter \& Patzschke [ArPa],
succeeded in computing the multifractal spectra 
$f_{\mu}(\alpha)$ assuming the OSC.
In order to state 
Patzschke's result we make the following definitions.
Define
$\Phi,\Lambda:\Sigma^{\Bbb N}\to\Bbb R$ by
$\Phi(\bold i)=\log p_{i_{1}}$ and 
$\Lambda(\bold i)=\log |DS_{i_{1}}(\pi S\bold i)|$
for
$\bold i=i_{1}i_{2}\ldots\in\Sigma^{\Bbb N}$,
and for $q\in\Bbb R$,
let 
$\beta(q)$
be the unique real number such that
 $$
 0
 =
 P\big(
 \,
 \beta(q)\Lambda
 +
 q\Phi
 \,
 \big)\,;
 $$
 here, and below, we use the following standard notation, namely
 if $\varphi:\Sigma^{\Bbb N}\to\Bbb R$ is a H\"older continuous function,
 then 
 $P(\varphi)$ denotes the pressure
 of  $\varphi$.
Also, recall that the Legendre transform is defined in (1.8).
We can now state 
Patzschke's result.

\bigskip

\proclaim{Theorem A [P]}
Let $\mu$ be defined by (3.2) and $\alpha\in\Bbb R$.
If the OSC is satisfied, then we have 
 $$
 f_{\mu}(\alpha)
 =
 \beta^{*}(\alpha)\,.
 $$
\endproclaim

 \bigskip

Of course, in general, the limit
 $\lim_{r\searrow0}
     \frac
     {\log\mu B(x,r)}{\log r}$
     may not exist.
     Indeed, recently 
     Barreira \& Schmeling [BaSc]
     (see also
     Olsen \& Winter [OlWi1,OlWi2],
     Xiao, Wu \& Gao [XiWuGa]
     and 
Moran [Mo])
have shown that 
the set of divergence points, 
i\.e\. the set 
 $$
 \align
 \Delta_{\mu}
&=
 \Bigg\{
 x\in K
     \,\Bigg|\,
     \text{the expression}
    \,\,
     \frac
     {\log \mu B(x,r)}{\log r}
     \,\,
     \text{diverges as $r\searrow0$}
     \,\,
      \Bigg\}
 \endalign  
 $$
of points $x$ for which the 
limit
$\lim_{r\searrow0}
     \frac
     {\log \mu B(x,r)}{\log r}$
does not exist, typically is  highly
\lq\lq visible" and
\lq\lq observable", namely it has full Hausdorff dimension.
More precisely, it follows from 
[BaSc]
% XiWuGa]
that if
the OSC is satisfied and $t$ denotes the Hausdorff 
dimension of $K$, then
 $$
  \Bigg\{
 x\in K
     \,\Bigg|\,
     \text{the expression}
    \,\,
     \frac
     {\log \mu B(x,r)}{\log r}
     \,\,
     \text{diverges as $r\searrow0$}
     \,\,
      \Bigg\}
 =
 \varnothing
 $$
provided $\mu$
is proportional to the $t$-dimensional Hausdorff measure restricted to $K$,
and
 $$
 \dim_{\Haus}
  \Bigg\{
 x\in K
     \,\Bigg|\,
     \text{the expression}
    \,\,
     \frac
     {\log \mu B(x,r)}{\log r}
     \,\,
     \text{diverges as $r\searrow0$}
     \,\,
      \Bigg\}
 =
 \dim_{\Haus} K
 $$ 
provided $\mu$
is 
not
proportional to the $t$-dimensional Hausdorff measure restricted to $K$.
This
 suggests that the set $\Delta_{\mu}$
 has a surprising rich and complex  
 fractal structure,
 and in order to explore this more
 carefully
  Olsen \& Winter [OlWi1,OlWi2]
 introduced various
 generalised multifractal spectra functions designed to 
 \lq\lq see"
 different sets of divergence points.
 In order to define these spectra 
 we introduce the following notation.
 If $M$ 
 is a  metric space
 and
 $\varphi:(0,\infty)\to M$ is a function, then we write
 $\acc_{r\searrow 0}f(r)$
 for the set of accumulation 
 points of $f$ as $r\searrow 0$, i\.e\.
  $$
  \underset {r\searrow0}\to\acc\,\,\varphi(r)
  =
  \Big\{x\in M
  \,\Big|\,
  \text{$x$ is an accumulation point of $f$ as $r\searrow 0$}
  \Big\}\,.
  $$
  In [OlWi1]
Olsen \& Winter
introduced and investigated
 the generalised
 Hausdorff multifractal spectrum
 $F_{\mu}$ of $\mu$ 
 defined
by
 $$
 \align
F_{\mu}(C)
&=
   \,\dim_{\Haus}
   \left\{x\in K
    \,\left|\,
     \,\,
     \underset {r\searrow0}\to\acc
 \,
     \frac
     {\log\mu B(x,r)}{\log r}
  \subseteq
 C
      \right.
      \,
   \right\}
  \endalign
 $$
 for 
 $C\subseteq \Bbb R$.
Note that the generalised spectrum
is a genuine extension of the
traditional multifractal spectrum 
$f_{\mu}(\alpha)$, namely if
$C=\{\alpha\}$ is a singleton consisting of the point $\alpha$, 
then clearly
$F_{\mu}(C)=f_{\mu}(\alpha)$.
There 
is a natural
divergence point analogue of Theorem A.
Indeed,
the following divergence point analogue of Theorem A
was first obtained by
Moran [Mo] and
Olsen \& Winter [OlWi1],
and later in a less restrictive setting by
Li, Wu \& Xiong [LiWuXi]
(see also [Ca,Vo]
for earlier
but related results in 
a
slightly different setting).

 \bigskip

\proclaim{Theorem B [LiWuXi,Mo,OlWi1]}
Let $\mu$ be defined by (3.2)
and
let $C$ be a closed subset of $\Bbb R$.
If the
 OSC is satisfied, then we have 
 $$
 F_{\mu}(C)
 =
 \sup_{\alpha\in C}\beta^{*}(\alpha)\,.
 $$
\endproclaim

As a first application of 
Theorem 2.1 and Theorem 2.2
we obtain a 
zeta-function
whose
abscissa of convergence equals the 
generalised 
multifractal spectrum $F_{\mu}(C)$
of a self-conformal measure $\mu$.
The is the content of the next theorem.

\bigskip

\proclaim{Theorem 3.1.
Multifractal zeta-functinons for
multifractal spectra of self-conformal measures}
Let
$(p_{1},\ldots,p_{N})$ be a probability vector, and let
$\mu$ denote the 
self-conformal measure associated with the list
$\big(
 \,
 V
 \,,\,
 X
 \,,\,
 (S_{i})_{i=1,\ldots,N}
  \,,\,\allowmathbreak
 (p_{i})_{i=1,\ldots,N}
 \,
 \big)$, i\.e\.
$\mu$ is the unique probability measure such that
$\mu
 =
 \sum_{i}p_{l,i}\mu\circ S_{i}^{-1}$.

For $\bold i\in\Sigma^{*}$, let
 $$
 s_{\bold i}
 =
 \sup_{\bold u\in\Sigma^{\Bbb N}}
 |DS_{\bold i}(\pi\bold u)|\,.
 $$
For a closed set $C\subseteq\Bbb R$, 
we define the self-conformal multifractal zeta-function by 
$$
 \zeta_{C}^{\scon}(s)
 =
 \sum
 \Sb
 \bold i\\
 {}\\
 \frac{\log p_{\bold i}}{\log\diam K_{\bold i}}
\in
 C
 \endSb
 s_{\bold i}^{s}\,,
 \qquad
 \,\,\,\,\,
 $$
For a closed set $C\subseteq\Bbb R$
and $r>0$, 
we define the self-conformal multifractal zeta-function by 
$$
\align
 \zeta_{C}^{\scon}(s;r)
&=
 \zeta_{B(C,r)}^{\scon}(s)\\
&= 
 \sum
 \Sb
 \bold i\\
 {}\\
 \dist
 \big(
 \,
 \frac{\log p_{\bold i}}{\log\diam K_{\bold i}}
 \,,\,
 C
 \,
 \big)
 \,
 \le 
 \,
 r
 \endSb
 s_{\bold i}^{s}\,,
 \endalign
 $$
and if $\alpha\in\Bbb R$
and
$C=\{\alpha\}$ is the singleton consisting of $\alpha$,
then
we
 write $ \zeta_{C}^{\scon}(s;r)=\zeta_{\alpha}^{\scon}(s;r)$, 
 i\.e\. 
 we
write
 $$
 \zeta_{\alpha}^{\scon}(s;r)
 =
 \sum
 \Sb
 \bold i\\
 {}\\
 \big|
 \,
 \frac{\log p_{\bold i}}{\log\diam K_{\bold i}}
 \,-\,
 \alpha
 \,
 \big|
 \,
 \le 
 \,
 r
 \endSb
 s_{\bold i}^{s}\,.
 \quad
 \,\,\,
 $$

Define $\Phi,\Lambda:\Sigma^{\Bbb N}\to\Bbb R$ by
$\Phi(\bold i)=\log p_{i_{1}}$ and 
$\Lambda(\bold i)=\log |DS_{i_{1}}(\pi S\bold i)|$
for
$\bold i=i_{1}i_{2}\ldots\in\Sigma^{\Bbb N}$.
Define 
$\beta:\Bbb R^{M}\to\Bbb R$ by
 $$
 0
 =
 P\big(
 \,
 \beta(q)\Lambda
 +
 q\Phi
 \,
 \big)
 $$
 for $q\in\Bbb R$.
 Let $C$ be a closed subset of $\Bbb R$.
Then the following hold:
\roster
\item"(1.1)"
We have
$$
 \lim_{r\searrow0}
 \sigma_{\abs}\big(\, \zeta_{C} ^{\scon}(\cdot;r)\,\big)
 =
 \sup_{\alpha\in C}
 \beta^{*}(\alpha)\,.
 $$
In particular, if $\alpha\in\Bbb R$, then we have
 $$
 \lim_{r\searrow0}
 \sigma_{\abs}\big(\, \zeta_{\alpha} ^{\scon}(\cdot;r)\,\big)
 =
 \beta^{*}(\alpha)\,.
 \qquad
 $$

\item"(1.2)"
If the OSC is satisfied, then we have
$$
\align
\quad
 \lim_{r\searrow0}
 \sigma_{\abs}\big(\, \zeta_{C} ^{\scon}(\cdot;r)\,\big)
 &=
 \sup_{\alpha\in C}
 \dim_{\Haus}
 \Bigg\{
 x\in K
 \,\Bigg|\,
 \lim_{r\searrow 0}
 \frac{\log\mu(B(x,r))}{\log r}
 =
 \alpha
 \Bigg\}\\ 
 &=
  \quad\,\,\,\,\,
 \dim_{\Haus}
 \Bigg\{
 x\in K
 \,\Bigg|\,
 \,\underset{r\searrow 0}\to\acc
 \frac{\log\mu(B(x,r))}{\log r}
\subseteq
C
 \Bigg\}\,.\\
 \endalign
 $$
In particular, if the OSC is satisfied and $\alpha\in\Bbb R$, then we have 
$$
\align
\quad
 \lim_{r\searrow0}
 \sigma_{\abs}\big(\, \zeta_{\alpha} ^{\scon}(\cdot;r)\,\big)
 &=
  \quad\,\,\,\,\,
 \dim_{\Haus}
 \Bigg\{
 x\in K
 \,\Bigg|\,
 \lim_{r\searrow 0}
 \frac{\log\mu(B(x,r))}{\log r}
 =
 \alpha
 \Bigg\}\,.
 \endalign
 $$

\item"(2.1)"
If $C$ is an interval and
$\overset{\circ}\to{C}\cap\,\big( -\beta'(\Bbb R)\big)\not=\varnothing$,
 then we have
$$
 \sigma_{\abs}\big(\, \zeta_{C} ^{\scon}\,\big)
 =
 \sup_{\alpha\in C}
 \beta^{*}(\alpha)\,.
 $$
 
\item"(2.2)"
If $C$ is an interval and
$\overset{\circ}\to{C}\cap \,\big( -\beta'(\Bbb R)\big)\not=\varnothing$
and
 the OSC is satisfied, then we have
$$
\align
\quad
 \sigma_{\abs}\big(\, \zeta_{C} ^{\scon}\,\big)
 &=
 \sup_{\alpha\in C}
 \dim_{\Haus}
 \Bigg\{
 x\in K
 \,\Bigg|\,
 \lim_{r\searrow 0}
 \frac{\log\mu(B(x,r))}{\log r}
 =
 \alpha
 \Bigg\}\\ 
 &=
  \quad\,\,\,\,\,
 \dim_{\Haus}
 \Bigg\{
 x\in K
 \,\Bigg|\,
 \,\underset{r\searrow 0}\to\acc
 \frac{\log\mu(B(x,r))}{\log r}
\subseteq
C
 \Bigg\}\,.\\
 \endalign
 $$

\endroster
\endproclaim
\noindent{\it Proof}\newline
This follows immediately from 
the more general Theorem 3.2 in Section 3.2 by putting $M=1$.
\hfill$\square$

\bigskip

{\bf 3.2.
Mixed multifractal 
spectra of 
self-conformal measures.}
Recently
mixed
(or simultaneous)
multifractal spectra
have generated an enormous interest in
the mathematical literature, see
[BaSa,Mo,Ol2,Ol3].
Indeed, previous
 result (Theorem A and Theorem B) only considered
the scaling behaviour  
of a single measure.
Mixed
multifractal analysis investigates
the {\it simultaneous} scaling behaviour
of finitely many 
measures.
Mixed
multifractal analysis
thus
combines local 
characteristics which depend simultaneously on various different 
aspects of the underlying dynamical system,
and provides the basis for a significantly
better understanding of 
the underlying dynamics.
We will now make these ideas precise.
For $m=1,\ldots, M$, let
$(p_{m,1},\ldots,p_{m,N})$ be a probability vector, and let
$\mu_{m}$ denote the 
self-conformal measure associated with the list
$\big(
 \,
 V
 \,,\,
 X
 \,,\,
 (S_{i})_{i=1,\ldots,N}
  \,,\,\allowmathbreak
 (p_{m,i})_{i=1,\ldots,N}
 \,
 \big)$, i\.e\.
$\mu_{m}$ is the unique probability measure such that
 $$
 \mu_{m}
 =
 \sum_{i}p_{m,i}\mu_{m}\circ S_{i}^{-1}\,.
 \tag3.3
 $$
The mixed multifractal spectrum 
$f_{\pmb\mu}$
of the list $\pmb\mu=(\mu_{1},\ldots,\mu_{M})$ is defined by
 $$
f_{\pmb\mu}(\pmb\alpha)
=
\dim_{\Haus}
 \Bigg\{
 x\in K
 \,\Bigg|\,
 \lim_{r\searrow 0}
 \Bigg(
 \frac{\log\mu_{1}(B(x,r))}{\log r}
 ,\ldots,
 \frac{\log\mu_{M}(B(x,r))}{\log r}
 \Bigg)
 =
 \pmb\alpha
 \Bigg\}
 $$
for $\pmb\alpha\in\Bbb R^{M}$.  
Of course, it is also possible to define 
generalised
mixed 
multifractal spectra
designed to
\lq\lq see"
different sets of divergence points.
Namely, we
define
the generalised  mixed
 Hausdorff multifractal spectrum
$F_{\pmb\mu}$
of the 
list $\pmb\mu=(\mu_{1},\ldots,\mu_{M})$
by
  $$
  F_{\pmb\mu}(C)
  =
\dim_{\Haus}
 \Bigg\{
 x\in K
 \,\Bigg|\,
 \,\underset{r\searrow 0}\to\acc
 \Bigg(
 \frac{\log\mu_{1}(B(x,r))}{\log r}
 ,\ldots,
 \frac{\log\mu_{M}(B(x,r))}{\log r}
 \Bigg)
 \subseteq
 C
 \Bigg\}
 $$
for  $C\subseteq\Bbb R^{M}$. 
Again we  
note that the generalised mixed multifractal spectrum
is a genuine extensions of the
traditional mixed multifractal spectrum
$F_{\mu}(\pmb\alpha)$, namely, if
$C=\{\pmb\alpha\}$ is a singleton consisting of the point $\pmb\alpha$, 
then clearly
$F_{\mu}(C)=f_{\mu}(\pmb\alpha)$.
Assuming the OSC,
 the  generalised mixed multifractal spectrum
$F_{\mu}(C)$ can be computed [Mo,Ol2].
 In order to state the result from [Mo,Ol2], 
 we introduce the following 
 definitions.
Define $\Lambda,\Phi_{m}:\Sigma^{\Bbb N}\to\Bbb R$ for $m=1,\ldots,M$ by
 $\Lambda(\bold i)=\log |DS_{i_{1}}(\pi S\bold i)|$
and
$\Phi_{m}(\bold i)=\log p_{m,i_{1}}$ 
for
$\bold i=i_{1}i_{2}\ldots\in\Sigma^{\Bbb N}$,
and
write $\pmb\Phi=(\Phi_{1},\ldots,\Phi_{M})$.
Define 
$\beta:\Bbb R^{M}\to\Bbb R$ by
 $$
 0
 =
 P\big(
 \,
 \beta(\bold q)\Lambda
 +
 \langle\bold q|\pmb\Phi\rangle
 \,
 \big)
 $$
for
$\bold q\in\Bbb R^{M}$;
recall that if $\varphi:\Sigma^{\Bbb N}\to\Bbb R$ is a 
H\"older continuous map, then $P(\varphi)$
denotes the pressure of $\varphi$.
Also, for $\bold x,\bold y\in\Bbb R^{M}$, we let
$\langle\bold x|\bold y\rangle$ denote the 
usual inner product of $\bold x$ and $\bold y$,
and 
if $\varphi:\Bbb R^{M}\to\Bbb R$ is a function, we define the Legendre transform
$\varphi^{*}:\Bbb R^{M}\to[-\infty,\infty]$ of $\varphi$ by
 $$
 \varphi^{*}(\bold x)
 =
 \inf_{\bold y}(\langle\bold x|\bold y\rangle+\varphi(\bold y))\,.
 $$
The  generalised mixed multifractal spectra
$f_{\pmb\mu}$
and
$F_{\pmb\mu}$ 
are now given by the following theorem.

\bigskip

\proclaim{Theorem C. [Mo,Ol2]}
Let $\mu_{1},\ldots,\mu_{M}$ be defined by (3.3)
and
let $C\subseteq\Bbb R^{M}$ be a
closed  set.
Put
$\pmb\mu
=
(\mu_{1},\ldots,\mu_{M})$.
If the OSC is satisfied,
then we have 
 $$
 F_{\pmb\mu}(C)
 =
 \sup_{\pmb\alpha\in C}\beta^{*}(\pmb\alpha)\,.
 $$
In particular,
if the OSC is satisfied
and $\pmb\alpha\in\Bbb R^{M}$, then we have
 $$
 f_{\pmb\mu}(\pmb\alpha)
 =
 \beta^{*}(\pmb\alpha)\,.
 $$
\endproclaim

\bigskip

As a second
 application of 
Theorem 2.1 and Theorem 2.2
we obtain a 
zeta-function
whose
abscissa of convergence equals the 
generalised 
mixed 
multifractal spectrum $F_{\pmb\mu}(C)$
of a list $\pmb\mu$ of self-conformal measures.
The is the content of the next theorem.

\bigskip

\proclaim{Theorem 3.2. 
Multifractal zeta-functinons for mixed
multifractal spectra of self-conformal measures}
For $m=1,\ldots, M$, let
$(p_{m,1},\ldots,p_{m,N})$ be a probability vector, and let
$\mu_{m}$ denote the 
self-conformal measure associated with the list
$\big(
 \,
 V
 \,,\,
 X
 \,,\,
 (S_{i})_{i=1,\ldots,N}
  \,,\,\allowmathbreak
 (p_{m,i})_{i=1,\ldots,N}
 \,
 \big)$, i\.e\.
$\mu_{m}$ is the unique probability measure such that
$\mu_{m}
 =
 \sum_{i}p_{m,i}\mu_{m}\circ S_{i}^{-1}$.

For $\bold i\in\Sigma^{*}$, let
 $$
 s_{\bold i}
 =
 \sup_{\bold u\in\Sigma^{\Bbb N}}
 |DS_{\bold i}(\pi\bold u)|\,.
 $$
For a closed set $C\subseteq\Bbb R^{M}$, 
we define the 
self-conformal
multifractal
zeta-function by
 $$
 \zeta_{C}^{\scon}(s)
 =
 \sum
 \Sb
 \bold i\\
 {}\\
 \big(
 \frac{\log p_{1,\bold i}}{\log \diam K_{\bold i}},
 \ldots,
 \frac{\log p_{M,\bold i}}{\log \diam K_{\bold i}}
 \big)
\in
 C
 \endSb
 s_{\bold i}^{s}
 $$
For a closed set $C\subseteq\Bbb R^{M}$  and $r>0$, 
we define the 
self-conformal
multifractal
zeta-function by
 $$
 \align
 \zeta_{C}^{\scon}(s;r)
&=
 \zeta_{B(C,r)}^{\scon}(s;r)\\
&=
 \sum
 \Sb
 \bold i\\
 {}\\
 \dist
 \big(
 \,
 \big(
 \frac{\log p_{1,\bold i}}{\log \diam K_{\bold i}},
 \ldots,
 \frac{\log p_{M,\bold i}}{\log \diam K_{\bold i}}
 \big)
 \,,\,
 C
 \,
 \big)
 \,
 \le 
 \,
 r
 \endSb
 s_{\bold i}^{s}
 \endalign
 $$

Define $\Lambda,\Phi_{m}:\Sigma^{\Bbb N}\to\Bbb R$ for $m=1,\ldots,M$ by
$\Lambda(\bold i)=\log |DS_{i_{1}}(\pi S\bold i)|$
and
$\Phi_{m}(\bold i)=\log p_{m,i_{1}}$ 
for
$\bold i=i_{1}i_{2}\ldots\in\Sigma^{\Bbb N}$,
and
write $\pmb\Phi=(\Phi_{1},\ldots,\Phi_{M})$.
Define 
$\beta:\Bbb R^{M}\to\Bbb R$ by
 $$
 0
 =
 P\big(
 \,
 \beta(\bold q)\Lambda
 +
 \langle\bold q|\pmb\Phi\rangle
 \,
 \big)
 $$
for
$\bold q\in\Bbb R^{M}$.
Let $C$ be a closed subset of $\Bbb R^{M}$.
Then the following hold:
\roster
\item"(1.1)"
We have
$$
 \lim_{r\searrow0}
 \sigma_{\abs}\big(\, \zeta_{C}^{\scon}(\cdot;r)\,\big)
 =
 \sup_{\pmb\alpha\in C}
 \beta^{*}(\pmb\alpha)\,.
 $$

\item"(1.2)"
If the OSC is satisfied, then we have
$$
\align
\quad
 \lim_{r\searrow0}
&\sigma_{\abs}\big(\, \zeta_{C}^{\scon}(\cdot;r)\,\big)\\
 &=
 \sup_{\pmb\alpha\in C}
 \dim_{\Haus}
 \Bigg\{
 x\in K
 \,\Bigg|\,
 \lim_{r\searrow 0}
 \Bigg(
 \frac{\log\mu_{1}(B(x,r))}{\log r}
 ,\ldots,
 \frac{\log\mu_{M}(B(x,r))}{\log r}
 \Bigg)
 =
 \pmb\alpha
 \Bigg\}\\
 &=
 \quad\,\,\,\,\,
 \dim_{\Haus}
 \Bigg\{
 x\in K
 \,\Bigg|\,
 \,\underset{r\searrow 0}\to\acc
 \Bigg(
 \frac{\log\mu_{1}(B(x,r))}{\log r}
 ,\ldots,
 \frac{\log\mu_{M}(B(x,r))}{\log r}
 \Bigg)
 \subseteq
 C
 \Bigg\}\,.\\
 \endalign
 $$

\item"(2.1)"
If $C$ is convex and
$\overset{\circ}\to{C}\cap \,\big( -\nabla\beta(\Bbb R^{M})\big)\not=\varnothing$,
then we have
$$
 \sigma_{\abs}\big(\, \zeta_{C}^{\scon}\,\big)
 =
 \sup_{\pmb\alpha\in C}
 \beta^{*}(\pmb\alpha)\,.
 $$
 
 \item"(2.2)"
 If $C$ is convex and
$\overset{\circ}\to{C}\cap \,\big( -\nabla\beta(\Bbb R^{M})\big)\not=\varnothing$
and
the OSC is satisfied, then we have
$$
\align
\quad
\sigma_{\abs}
&\big(\, \zeta_{C}^{\scon}\,\big)\\
 &=
 \sup_{\pmb\alpha\in C}
 \dim_{\Haus}
 \Bigg\{
 x\in K
 \,\Bigg|\,
 \lim_{r\searrow 0}
 \Bigg(
 \frac{\log\mu_{1}(B(x,r))}{\log r}
 ,\ldots,
 \frac{\log\mu_{M}(B(x,r))}{\log r}
 \Bigg)
 =
 \pmb\alpha
 \Bigg\}\\
 &=
 \quad\,\,\,\,\,
 \dim_{\Haus}
 \Bigg\{
 x\in K
 \,\Bigg|\,
 \,\underset{r\searrow 0}\to\acc
 \Bigg(
 \frac{\log\mu_{1}(B(x,r))}{\log r}
 ,\ldots,
 \frac{\log\mu_{M}(B(x,r))}{\log r}
 \Bigg)
 \subseteq
 C
 \Bigg\}\,.\\
 \endalign
 $$

\endroster

\endproclaim

\bigskip

We  will now prove Theorem 3.2.
Recall  that the function $\Lambda:\Sigma^{\Bbb N}\to\Bbb R$ is defined by
 $$
 \Lambda(\bold i)=\log|DS_{i_{1}}(\pi S\bold i)|
 \tag3.4
 $$
 for $\bold i=i_{1}i_{2}\ldots\in\Sigma^{\Bbb N}$.
It is well-known that $\Lambda$
satisfies Conditions (C1)--(C3) in Section 2.1.
Also, a straight forward calculation shows that
$\sup_{\bold k\in[\bold i]}\exp\sum_{k=0}^{|\bold i|-1}\Lambda S^{k}\bold k
=
\sup_{\bold u\in\Sigma^{\Bbb N}}|DS_{\bold i}(\pi\bold u)|=s_{\bold i}$
for
$\bold i\in\Sigma^{*}$.
Next,
recall that
$\pmb\Phi=(\Phi_{1},\ldots,\Phi_{M})$
where
$\Phi_{m}:\Sigma^{\Bbb N}\to\Bbb R$ is defined by
$\Phi_{m}(\bold i)=\log p_{m,i_{1}}$ 
for
$\bold i=i_{1}i_{2}\ldots\in\Sigma^{\Bbb N}$.
For $\mu\in\Cal P(\Sigma^{\Bbb N})$, we will write
$\int\pmb\Phi\,d\mu
=
(\int\Phi_{1}\,d\mu,\ldots,\int\Phi_{M}\,d\mu)$.
Finally, 
define
$U:\Cal P(\Sigma^{\Bbb N})\to\Bbb R^{M}$ by
 $$
 U\mu
 =
 \frac{\int\pmb\Phi\,d\mu}{\int\Lambda\,d\mu}\,,
 \tag3.5
 $$
and note that if $\bold i\in\Sigma^{*}$, then
 $$
 UL_{|\bold i|}[\bold i]
 =
 \Bigg\{
  \Bigg(
 \,
 \frac{\log p_{1,\bold i}}{\log |DS_{\bold i}(\pi\bold u)|}
 \,
 ,
 \ldots
 ,
 \,
 \frac{\log p_{M,\bold i}}{\log |DS_{\bold i}(\pi\bold u)|} 
 \,
 \Bigg)
 \,\Bigg|\,
 \bold u\in\Sigma^{\Bbb N}
 \Bigg\}\,.
 $$
Hence, for $C\subseteq\Bbb R^{M}$ we have
 $$
 \align
 \zeta_{C}^{U,\Lambda}(s;r)
&=
 \sum
 \Sb
 \bold i\\
 {}\\
 UL_{|\bold i|}[\bold i]
  \subseteq
 B(C,r)
 \endSb
 s_{\bold i}^{s}\\ 
&{}\\ 
&=
 \sum
 \Sb
 \bold i\\
 {}\\
  \big\{
  \big(
 \,
 \frac{\log p_{1,\bold i}}{\log |DS_{\bold i}(\pi\bold u)|}
 \,
 ,
 \ldots
 ,
 \,
 \frac{\log p_{M,\bold i}}{\log |DS_{\bold i}(\pi\bold u)|} 
 \,
 \big)
 \,\big|\,
 \bold u\in\Sigma^{\Bbb N}
 \big\}
 \subseteq
 B(C,r)
 \endSb
 s_{\bold i}^{s}\\
&{}\\ 
&=
 \sum
 \Sb
 \bold i\\
 {}\\
 \forall
 \bold u\in\Sigma^{\Bbb N}
 \,\,:\,\,
 \dist
 \big(
 \,
  \big(
 \,
 \frac{\log p_{1,\bold i}}{\log |DS_{\bold i}(\pi\bold u)|}
 \,
 ,
 \ldots
 ,
 \,
 \frac{\log p_{M,\bold i}}{\log |DS_{\bold i}(\pi\bold u)|} 
 \,
 \big)
 \,,\,
 C
 \,
 \big)
 \,
 \le
 \,
 r
 \endSb
 s_{\bold i}^{s}\,.
 \tag3.6
 \endalign 
 $$
 In order to prove Theorem 3.2, we first
 prove  the following three auxiliary results, namely,
 Propositions 3.3--3.5.

\bigskip

\proclaim{Proposition 3.3}
Let $U$ and $\Lambda$
be defined by (3.5) and (3.4), respectively.
For $\pmb\alpha\in\Bbb R^{M}$,
we have
 $$
 \sup
 \Sb
 \mu\in\Cal P_{S}(\Sigma^{\Bbb N})\\
 {}\\
 U\mu=\pmb\alpha
 \endSb
 -
 \frac{h(\mu)}{\int\Lambda\,d\mu}
 =
 \beta^{*}(\pmb\alpha)\,.
$$
\endproclaim
\noindent{\it Proof}\newline
This result is folklore for $M=1$.
%and various versions of the result for $M=1$
%can be found in [???].
The proof of Proposition 3.3
for an arbitrary positive integer
can (with some modifications)
be modelled
on the argument for $M=1$.
However, for the sake of brevity we have decided to omit the proof.
\hfill$\square$

\bigskip

\proclaim{Proposition 3.4}
Let $U$ and $\Lambda$
be defined by (3.5) and (3.4), respectively.
Let $C$ be a closed subset of $\Bbb R^{M}$.
 If $C$ is convex and
$\overset{\circ}\to{C}\cap \,\big( -\nabla\beta(\Bbb R^{M})\big)\not=\varnothing$, 
then 
$\underline f^{U,\Lambda}$ is inner continuous at $C$.
\endproclaim
\noindent{\it Proof}\newline
\noindent
Note that it follows from 
Theorem 2.1 and Proposition 3.3
that if $W$ is closed subset of $\Bbb R^{M}$, then
 $$
 \align
 \underline f^{U,\Lambda}(W)
&=
  \sup
  \Sb
 \mu\in\Cal P_{S}(\Sigma^{\Bbb N})\\
 {}\\
 U\mu\in W
 \endSb
 -
 \frac{h(\mu)}{\int\Lambda\,d\mu}\\
&=
 \sup_{\pmb\alpha\in W}
 \sup
 \Sb
 \mu\in\Cal P_{S}(\Sigma^{\Bbb N})\\
 {}\\
 U\mu=\pmb\alpha
 \endSb
 -
 \frac{h(\mu)}{\int\Lambda\,d\mu}\\
&=
 \sup_{\pmb\alpha\in W} 
 \beta^{*}(\pmb\alpha)\,.
 \tag3.8
 \endalign
 $$
Also, since
the function
$\beta^{*}$ with
 $\{\pmb\alpha\in\Bbb R^{M}
 \,|\,
 \beta^{*}(\pmb\alpha)>-\infty
 \}
 =
 \nabla\beta(\Bbb R^{M})$
 (see [Ro, Corollary 26.4.1])
and
the set
$C$ is convex 
with
$\overset{\circ}\to{C}\cap \,\big( -\nabla\beta(\Bbb R^{M})\big)\not=\varnothing$, 
we conclude immediately from (3.8)
that 
$\underline f^{U,\Lambda}$
is inner continuous at $C$.
\hfill$\square$

\bigskip

\proclaim{Proposition 3.5}
Let $U$ and $\Lambda$
be defined by (3.5) and (3.4), respectively.
\roster
\item"(1)"
There is a sequence $(\Delta_{n})_{n}$ with $\Delta_{n}>0$ 
% for all $n$ 
and 
$\Delta_{n}\to 0$ such that
for all
closed subsets $C$ of $\Bbb R^{M}$
and
for
all
$n\in\Bbb N$, $\bold i\in\Sigma^{n}$ and $\bold u\in\Sigma^{\Bbb N}$, we have
 $$
 \align
 \quad\quad\quad\,\,\,
 \dist
 \Bigg(
 \,
  \Bigg(
 \,
&\frac{\log p_{1,\bold i}}{\log |DS_{\bold i}(\pi\bold u)|}
 \,
,
 \ldots
 ,
 \,
 \frac{\log p_{M,\bold i}}{\log |DS_{\bold i}(\pi\bold u)|} 
 \,
 \Bigg)
 \,,\,
 C
 \,
 \Bigg)\\
&\le
\dist
\Bigg(
 \,
  \Bigg(
 \,
 \frac{\log p_{1,\bold i}}{\log \diam K_{\bold i}}
 \,
,
 \ldots
 ,
 \,
 \frac{\log p_{M,\bold i}}{\log \diam K_{\bold i}} 
 \,
 \Bigg)
 \,,\,
 C
 \,
 \Bigg)
\,\,
\,\,\,\,\,\,\,\,\,\,
+
\,\,
\Delta_{n}\,,
\tag3.9\\
&{}\\
\dist
\Bigg(
 \,
  \Bigg(
 \,
&\frac{\log p_{1,\bold i}}{\log \diam K_{\bold i}}
 \,
,
 \ldots
 ,
 \,
 \frac{\log p_{M,\bold i}}{\log \diam K_{\bold i}} 
 \,
 \Bigg)
 \,,\,
 C
 \,
 \Bigg)
\\
&\le
\dist
 \Bigg(
 \,
  \Bigg(
 \,
 \frac{\log p_{1,\bold i}}{\log |DS_{\bold i}(\pi\bold u)|}
 \,
,
 \ldots
 ,
 \,
 \frac{\log p_{M,\bold i}}{\log |DS_{\bold i}(\pi\bold u)|} 
 \,
 \Bigg)
 \,,\,
 C
 \,
 \Bigg)
 \,\,
+
\,\,
\Delta_{n}\,.
\tag3.10
\endalign
$$

\item"(2)"
For all closed subsets $W$ of $\Bbb R^{M}$ and all $r>0$, we have
 $$
 \align
 \sigma_{\abs}
 \big(
 \,
 \zeta_{W}^{U,\Lambda}(\cdot;r)
 \,
 \big)
&\le
 \sigma_{\abs}
 \big(
 \,
 \zeta_{B(W,2r)}^{\scon}
 \,
 \big)\,,
 \tag3.11\\
 \sigma_{\abs}
 \big(
 \,
 \zeta_{B(W,r)}^{\scon}
 \,
 \big) 
&\le
 \sigma_{\abs}
 \big(
 \,
 \zeta_{W}^{U,\Lambda}(\cdot;2r)
 \,
 \big)\,.
 \tag3.12
 \endalign
 $$

\item"(3)"
Let $C$ be a closed subset of $\Bbb R^{M}$.
We have
 $$
 \lim_{r\searrow0}
 \sigma_{\abs}
 \big(
 \,
 \zeta_{C}^{\scon}(\cdot;r)
 \,
 \big)
 =
 \lim_{r\searrow0}
 \sigma_{\abs}
 \big(
 \,
 \zeta_{C}^{U,\Lambda}(\cdot;r)
 \,
 \big)\,.
$$

\item"(4)"
Let $C$ be a closed subset of $\Bbb R^{M}$.
 If $C$ is convex and
$\overset{\circ}\to{C}\cap \,\big( -\nabla\beta(\Bbb R^{M})\big)\not=\varnothing$, 
then we have
 $$
 \sigma_{\abs}
 \big(
 \,
 \zeta_{C}^{\scon}
 \,
 \big)
 =
 \sigma_{\abs}
 \big(
 \,
 \zeta_{C}^{U,\Lambda}
 \,
 \big)\,.
$$

\endroster
\endproclaim
\noindent{\it Proof}\newline
\noindent
(1)
It is well-known that there is a constant $c_{0}>0$ such that
for all $\bold i\in\Sigma^{*}$ and all $\bold u\in\Sigma^{\Bbb N}$,
 we have
$\frac{1}{c_{0}}
\le
\frac{\diam K_{\bold i}}{|DS_{\bold i}(\pi\bold u)|}
\le
c_{0}$, see, for example, [Fa2] or [Pa].
It is not difficult to see that the desired result follows from this
and the fact 
that the function $\Lambda:\Sigma^{\Bbb N}\to\Bbb R$
defined
by
$\Lambda(\bold i)=\log|DS_{i_{1}}(\pi S\bold i)|$ for $\bold i=i_{1}i_{2}\ldots\in\Sigma^{\Bbb N}$
satisfies Conditions (C1)--(C3) in Section 2.1.

\noindent
(2)
Fix $r>0$.
Let $(\Delta_{n})_{n}$ be the sequence from (1).
Since $\Delta_{n}\to 0$, we can find
a positive integer
$N_{r}$
such that if $n\ge N_{r}$, then $\Delta_{n}\le r$.
Consequently, using (3.10) in Part (1),
for $s\in\Bbb R$, we have
$$
 \align
 \zeta_{W}^{U,\Lambda}(s;r)
&=
 \sum
 \Sb
 \bold i\\
 {}\\
 \forall
 \bold u\in\Sigma^{\Bbb N}
 \,\,:\,\,
 \dist
 \big(
 \,
  \big(
 \,
 \frac{\log p_{1,\bold i}}{\log |DS_{\bold i}(\pi\bold u)|}
 \,
 ,
 \ldots
 ,
 \,
 \frac{\log p_{M,\bold i}}{\log |DS_{\bold i}(\pi\bold u)|} 
 \,
 \big)
 \,,\,
 W
 \,
 \big)
 \,
 \le
 \,
 r
 \endSb
 s_{\bold i}^{s}\\
&{}\\ 
&\le
  \sum
 \Sb
 \bold i\\
 {}\\
 |\bold i|< N_{r}
 \endSb
 s_{\bold i}^{s}
 \,\,
 +
 \,\,
   \sum
 \Sb
 \bold i\\
 {}\\
 |\bold i|\ge N_{r}\\
 {}\\
 \forall
 \bold u\in\Sigma^{\Bbb N}
 \,\,:\,\,
 \dist
 \big(
 \,
  \big(
 \,
 \frac{\log p_{1,\bold i}}{\log |DS_{\bold i}(\pi\bold u)|}
 \,
 ,
 \ldots
 ,
 \,
 \frac{\log p_{M,\bold i}}{\log |DS_{\bold i}(\pi\bold u)|} 
 \,
 \big)
 \,,\,
 W
 \,
 \big)
 \,
 \le
 \,
 r
 \endSb
 s_{\bold i}^{s}\\
&{}\\ 
&\le
  \sum
 \Sb
 \bold i\\
 {}\\
 |\bold i|< N_{r}
 \endSb
 s_{\bold i}^{s}
 \,\,
 +
 \,\,
   \sum
 \Sb
 \bold i\\
 {}\\
 |\bold i|\ge N_{r}\\
 {}\\
 \dist
 \big(
 \,
  \big(
 \,
 \frac{\log p_{1,\bold i}}{\log \diam K_{\bold i}}
 \,
 ,
 \ldots
 ,
 \,
 \frac{\log p_{M,\bold i}}{\log \diam K_{\bold i}} 
 \,
 \big)
 \,,\,
 W
 \,
 \big)
 \,
 \le
 \,
 r+\Delta_{|\bold i|}
 \endSb
 s_{\bold i}^{s}\\
&{}\\ 
&\le
  \sum
 \Sb
 \bold i\\
 {}\\
 |\bold i|< N_{r}
 \endSb
 s_{\bold i}^{s}
 \,\,
 +
 \,\,
   \sum
 \Sb
 \bold i\\
 {}\\
 |\bold i|\ge N_{r}\\
 {}\\
 \dist
 \big(
 \,
  \big(
 \,
 \frac{\log p_{1,\bold i}}{\log \diam K_{\bold i}}
 \,
 ,
 \ldots
 ,
 \,
 \frac{\log p_{M,\bold i}}{\log \diam K_{\bold i}} 
 \,
 \big)
 \,,\,
 W
 \,
 \big)
 \,
 \le
 \,
 2r
 \endSb
 s_{\bold i}^{s}  \\
&{}\\ 
&\le
  \sum
 \Sb
 \bold i\\
 {}\\
 |\bold i|< N_{r}
 \endSb
 s_{\bold i}^{s}
 \,\,
 +
 \,\,
   \sum
 \Sb
 \bold i\\
 {}\\
 \dist
 \big(
 \,
  \big(
 \,
 \frac{\log p_{1,\bold i}}{\log \diam K_{\bold i}}
 \,
 ,
 \ldots
 ,
 \,
 \frac{\log p_{M,\bold i}}{\log \diam K_{\bold i}} 
 \,
 \big)
 \,,\,
 W
 \,
 \big)
 \,
 \le
 \,
 2r
 \endSb
 s_{\bold i}^{s}   \\
&{}\\ 
&=
  \sum
 \Sb
 \bold i\\
 {}\\
 |\bold i|< N_{r}
 \endSb
 s_{\bold i}^{s}
 \,\,
 +
 \,\,
 \zeta_{B(W,2r)}^{\scon}(s)\,. 
 \tag3.13
 \endalign 
 $$
A similar argument using (3.1) in Part 1 shows that
 $$
 \zeta_{B(W,r)}^{\scon}(s)
 \le 
   \sum
 \Sb
 \bold i\\
 {}\\
 |\bold i|< N_{r}
 \endSb
 s_{\bold i}^{s}
 \,\,
 +
 \,\,
   \zeta_{W}^{U,\Lambda}(s;2r)\,.
 \tag3.14
 $$
The desired results follow immediately from inequalities
(3.13) and (3.14).

\noindent
(3)
This result follows from Part (2) by letting $r\searrow0$.

\noindent
(4)
\lq\lq$\le$"
It follows from (3.12) and Theorem 2.1  that
 $$
 \align
 \sigma_{\abs}
 \big(
 \,
 \zeta_{C}^{\scon}
 \,
 \big)
&\le
\liminf_{r\searrow 0}
  \sigma_{\abs}
 \big(
 \,
 \zeta_{B(C,r)}^{\scon}
 \,
 \big)
 \qquad\qquad
 \quad\,\,\,
 \text{[since $C\subseteq B(C,r)$]}\\
&\le
\liminf_{r\searrow 0}
  \sigma_{\abs}
 \big(
 \,
 \zeta_{C}^{U,\Lambda}(\cdot;2r)
 \,
 \big)
 \qquad\qquad
 \text{[by (3.12)]}\\ 
&=
 \underline f^{U,\Lambda}(C)
  \qquad\qquad
   \qquad\qquad
    \qquad\,\,
 \text{[by Theorem 2.1]}\,.
 \tag3.15
 \endalign
 $$
 Next, since 
$C$ is convex and
$\overset{\circ}\to{C}\cap \,\big( -\nabla\beta(\Bbb R^{M})\big)\not=\varnothing$, 
we conclude from Proposition 3.4
that 
$\underline f^{U,\Lambda}$
is inner continuous at $C$, and it therefore follows from Theorem 2.2 that
$\underline f^{U,\Lambda}(C)
 =
  \sigma_{\abs}
 \big(
 \,
 \zeta_{C}^{U,\Lambda}
 \,
 \big)$.
The desired result follows from this and (3.15).

\noindent
\lq\lq$\ge$"
Let $\varepsilon>0$.
For all $r>0$ with $2r<\varepsilon$, it follows from 
(3.11) applied to $W=I(C,\varepsilon)$
(recall that 
$I(C,\varepsilon)
=
 \{
 x\in C
 \,|\,
 \dist(x,\partial C)\ge \varepsilon
 \}$, see Section 2.3)
that
 $$
\sigma_{\abs}
 \big(
 \,
 \zeta_{I(C,\varepsilon)}^{U,\Lambda}(\cdot;r)
 \,
 \big)
 \le
 \sigma_{\abs}
 \big(
 \,
 \zeta_{B(I(C,\varepsilon),2r)}^{\scon}
 \,
 \big)\,.
 \tag3.16
$$
However, for $2r<\varepsilon$
it is not difficult to see that
$B(I(C,\varepsilon),2r)\subseteq C$
(see, for example, the proof of Lemma 8.2), whence
$ \sigma_{\abs}
 \big(
 \,
 \zeta_{B(I(C,\varepsilon),2r)}^{\scon}
 \,
 \big)
 \le
  \sigma_{\abs}
 \big(
 \,
 \zeta_{C}^{\scon}
 \,
 \big)$,
 and we therefore conclude from
(3.16) that
if $2r<\varepsilon$, then
  $$
\sigma_{\abs}
 \big(
 \,
 \zeta_{I(C,\varepsilon)}^{U,\Lambda}(\cdot;r)
 \,
 \big)
 \le
 \sigma_{\abs}
 \big(
 \,
 \zeta_{C}^{\scon}
 \,
 \big)\,.
 \tag3.17
$$
Letting $r\searrow 0$ in (3.17) we now deduce that
$$
\limsup_{r\searrow0}
\sigma_{\abs}
 \big(
 \,
 \zeta_{I(C,\varepsilon)}^{U,\Lambda}(\cdot;r)
 \,
 \big)
 \le
 \sigma_{\abs}
 \big(
 \,
 \zeta_{C}^{\scon}
 \,
 \big)\,.
 \tag3.18
$$
Next, 
since 
$I(C,\varepsilon)$ is closed, we deduce from Theorem 2.1
that
$\limsup_{r\searrow0}
\sigma_{\abs}
 \big(
 \,
 \zeta_{I(C,\varepsilon)}^{U,\Lambda}(\cdot;r)
 \,
 \big)
 =
 \underline f^{U,\Lambda}(\,I(C,\varepsilon)\,)$,
 and (3.18) therefore implies that
 $$
\underline f^{U,\Lambda}(\,I(C,\varepsilon)\,)
 \le
 \sigma_{\abs}
 \big(
 \,
 \zeta_{C}^{\scon}
 \,
 \big)\,.
 \tag3.19
$$
Finally, it follows from Proposition 3.4 that 
$\underline f^{U,\Lambda}$ is 
is inner continuous at $C$, whence
$\lim_{\varepsilon\searrow 0}\underline f^{U,\Lambda}(\,I(C,\varepsilon)\,)
=
\underline f^{U,\Lambda}(C)$.
The desired result follows from this and (3.19).
\hfill$\square$

\bigskip

We can now prove Theorem 3.2.

\bigskip

\noindent{\it Proof of Theorem 3.2}\newline
(1.1) and (2.1)
The statements in Part (1.1) and Part (2.1)
of
Theorem 3.2 follow immediately from 
Theorem 2.1, 
Proposition 3.3
and
Proposition 3.5.

\noindent
(1.2) and (2.2)
The statements in Part (1.2) and Part (2.2)
of
Theorem 3.2 follow immediately from 
Part (1.1) and Part (2.1)
using
Theorem 2.2
and
Theorem C.
\hfill$\square$

\bigskip

{\bf 3.3. Multifractal spectra of self-similar measures.}
Due to important role 
 self-similar measures 
play in fractal geometry,
it is instructive to note the following special case of Theorem 3.1.

\bigskip

\proclaim{Theorem 3.6.
Multifractal zeta-functinons for
multifractal spectra of self-similar measures}
Assume that the maps $S_{1},\ldots, S_{N}$ are contracting similarities and let $r_{i}$ denote the contraction ratio of $S_{i}$.
For $\bold i=i_{1}\ldots i_{n}\in\Sigma^{*}$, let
 $$
 r_{\bold i}
 =
 r_{i_{1}}\cdots r_{i_{n}}\,.
 $$
Let
$(p_{1},\ldots,p_{N})$ be a probability vector, and let
$\mu$ denote the 
self-conformal measure associated with the list
$\big(
 \,
 V
 \,,\,
 X
 \,,\,
 (S_{i})_{i=1,\ldots,N}
  \,,\,\allowmathbreak
 (p_{i})_{i=1,\ldots,N}
 \,
 \big)$, i\.e\.
$\mu$ is the unique probability measure such that
$\mu
 =
 \sum_{i}p_{l,i}\mu\circ S_{i}^{-1}$.

For a closed set $C\subseteq\Bbb R$,
we define the self-similar multifractal zeta-function 
by
 $$
 \zeta_{C}^{\ssim}(s)
 =
 \sum
 \Sb
 \bold i\\
 {}\\
 \frac{\log p_{\bold i}}{\log\diam K_{\bold i}}
\in
C
 \endSb
 r_{\bold i}^{s}\,,
 $$ 
For a closed set $C\subseteq\Bbb R$ and $r>0$,
we define the self-similar multifractal zeta-function 
by
 $$
 \zeta_{C}^{\ssim}(s;r)
 =
 \sum
 \Sb
 \bold i\\
 {}\\
 \dist
 \big(
 \,
 \frac{\log p_{\bold i}}{\log\diam K_{\bold i}}
 \,,\,
 C
 \,
 \big)
 \,
 \le 
 \,
 r
 \endSb
 r_{\bold i}^{s}\,,
 $$
and if $\alpha\in\Bbb R$
and $C=\{\alpha\}$
is the singleton consisting of $\alpha$, then
we
 write $ \zeta_{C}(s;r)=\zeta_{\alpha}(s;r)$, 
 i\.e\. 
 we
write
 $$
 \zeta_{\alpha}^{\ssim}(s;r)
 =
 \sum
 \Sb
 \bold i\\
 {}\\
 \big|
 \,
 \frac{\log p_{\bold i}}{\log\diam K_{\bold i}}
 \,-\,
 \alpha
 \,
 \big|
 \,
 \le 
 \,
 r
 \endSb
 r_{\bold i}^{s}\,.
 $$

Define 
$\beta:\Bbb R^{M}\to\Bbb R$ by
 $$
 \sum_{i}p_{i}^{q}r_{i}^{\beta(q)}=1
 $$
for $q\in\Bbb R$. 
Let $C$ be a closed subset of $\Bbb R$.
Then the following hold:
\roster
\item"(1.1)"
We have
$$
 \lim_{r\searrow0}
 \sigma_{\abs}\big(\, \zeta_{C}^{\ssim}(\cdot;r)\,\big)
 =
 \sup_{\alpha\in C}
 \beta^{*}(\alpha)\,.
 $$
In particular, if $\alpha\in\Bbb R$, then we have
 $$
 \lim_{r\searrow0}
 \sigma_{\abs}\big(\, \zeta_{\alpha}^{\ssim}(\cdot;r)\,\big)
 =
 \beta^{*}(\alpha)\,.
 \qquad
 $$

\item"(1.2)"
If the OSC is satisfied, then we have
$$
\align
\quad
 \lim_{r\searrow0}
 \sigma_{\abs}\big(\, \zeta_{C} ^{\ssim}(\cdot;r)\,\big)
 &=
 \sup_{\alpha\in C}
 \dim_{\Haus}
 \Bigg\{
 x\in K
 \,\Bigg|\,
 \lim_{r\searrow 0}
 \frac{\log\mu(B(x,r))}{\log r}
 =
 \alpha
 \Bigg\}\\ 
 &=
  \quad\,\,\,\,\,
 \dim_{\Haus}
 \Bigg\{
 x\in K
 \,\Bigg|\,
 \,\underset{r\searrow 0}\to\acc
 \frac{\log\mu(B(x,r))}{\log r}
\subseteq
C
 \Bigg\}\,.\\
 \endalign
 $$
In particular, if the OSC is satisfied and $\alpha\in\Bbb R$, then we have 
$$
\align
\quad
 \lim_{r\searrow0}
 \sigma_{\abs}\big(\, \zeta_{\alpha} ^{\ssim}(\cdot;r)\,\big)
 &=
  \quad\,\,\,\,\,
 \dim_{\Haus}
 \Bigg\{
 x\in K
 \,\Bigg|\,
 \lim_{r\searrow 0}
 \frac{\log\mu(B(x,r))}{\log r}
 =
 \alpha
 \Bigg\}\,.
 \endalign
 $$
 
 \item"(2.1)"
 If $C$ is an interval
 and
 $\overset{\circ}\to{C}
 \cap
 \big(\min_{i}\frac{\log p_{i}}{\log r_{i}},\max_{i}\frac{\log p_{i}}{\log r_{i}}\big)
 \not=\varnothing$, then we have
$$
 \sigma_{\abs}\big(\, \zeta_{C}^{\ssim}\,\big)
 =
 \sup_{\alpha\in C}
 \beta^{*}(\alpha)\,.
 $$

\item"(2.2)"
If $C$ is an interval
and
$\overset{\circ}\to{C}
 \cap
 \big(\min_{i}\frac{\log p_{i}}{\log r_{i}},\max_{i}\frac{\log p_{i}}{\log r_{i}}\big)
 \not=\varnothing$
 and the OSC is satisfied, then we have
$$
\align
\quad
 \sigma_{\abs}\big(\, \zeta_{C} ^{\ssim}\,\big)
 &=
 \sup_{\alpha\in C}
 \dim_{\Haus}
 \Bigg\{
 x\in K
 \,\Bigg|\,
 \lim_{r\searrow 0}
 \frac{\log\mu(B(x,r))}{\log r}
 =
 \alpha
 \Bigg\}\\ 
 &=
  \quad\,\,\,\,\,
 \dim_{\Haus}
 \Bigg\{
 x\in K
 \,\Bigg|\,
 \,\underset{r\searrow 0}\to\acc
 \frac{\log\mu(B(x,r))}{\log r}
\subseteq
C
 \Bigg\}\,.\\
 \endalign
 $$

\endroster
\endproclaim
\noindent{\it Proof}\newline
\noindent
Theorem 3.6 follows immediately from Theorem 3.1.
\hfill$\square$

\bigskip

\noindent
It is, of course, also possible to formulate a 
version of 
Theorem 3.2 for 
a finite list self-similar measures.
However, for sake of brevity we have decided not to do this.

\bigskip

{\bf 3.4. Multifractal spectra of ergodic Birkhoff averages.}
We first fix $\gamma\in(0,1)$ and define  the metric
 $\distance_{\gamma}$ on $\Sigma^{\Bbb N}$ by
$\distance_{\gamma}(\bold i,\bold j)=\gamma^{\max\{n\,|\,\bold i|n=\bold j|n\}}$;
throughout this section, we equip
 $\Sigma^{\Bbb N}$
with the metric  $\distance_{\gamma}$ 
and continuity and Lipschitz properties of functions $f:\Sigma^{\Bbb N}\to\Bbb R$
from $\Sigma^{\Bbb N}$ to $\Bbb R$
will  always
refer to
the metric  $\distance_{\gamma}$.
Multifractal analysis of Birkhoff
averages has received significant interest
during
the past 10 years, see, for example,
[BaMe,FaFe,FaFeWu,FeLaWu,Oli,Ol3,OlWi2].
The 
multifractal
spectrum
$F_{f}^{\erg}$
of ergodic Birkhoff averages of a continuous function
$f:\Sigma^{\Bbb N}\to\Bbb R$ is defined by
 $$
F_{f}^{\erg}(\alpha)
=
 \dim_{\Haus}
 \pi
 \Bigg\{
 \bold i\in\Sigma^{\Bbb N}
 \,\Bigg|\,
 \lim_{n}
 \frac{1}{n}\sum_{k=0}^{n-1}f(S^{k}\bold i)
 =
 \alpha
 \Bigg\}
$$
for $\alpha\in\Bbb R$;
recall that the projection map $\pi:\Sigma^{\Bbb N}\to\Bbb R^{d}$
is defined in Section 3.1
and that
$S:\Sigma^{\Bbb N}\to\Sigma^{\Bbb N}$
denotes the shift map.
One of the main problems
in
multifractal analysis of Birkhoff
averages
is the detailed study of the multifractal
spectrum
$F_{f}^{\erg}$.
For example,
Theorem D below is
proved
in different settings and at various levels of generality
in [FaFe,FaFeWu,FeLaWu,Oli,Ol3,OlWi2].
Before we can state
we introduce the following 
notation.
If $(x_{n})_{n}$ is a sequence of real numbers, then we write
$\acc_{n} x_{n}$ for the set of accumulation points 
of $(x_{n})_{n}$, i\.e\.
 $$
 \underset{n}\to{\acc} \,x_{n}
=
\Big\{
x\in\Bbb R
\,\Big|\,
\text{
$x$ is an accumulation point of $(x_{n})_{n}$
}
\Big\}\,.
$$
Also, recall that
$\Cal P_{S}(\Sigma^{\Bbb N})$ denotes the family of shift invariant Borel probability
measures on $\Sigma^{\Bbb N}$ 
and 
that
$h(\mu)$ denotes the entropy of $\mu\in \Cal P_{S}(\Sigma^{\Bbb N})$.
We can now state Theorem D.

\bigskip

\proclaim{Theorem D. [FaFe,FaFeWu,FeLaWu,Oli,Ol3,OlWi2]}
Let $f:\Sigma^{\Bbb N}\to\Bbb R$ be a Lipschitz function.
Define $\Lambda:\Sigma^{\Bbb N}\to\Bbb R$ by
$\Lambda(\bold i)=\log |DS_{i_{1}}(\pi S\bold i)|$
for
$\bold i=i_{1}i_{2}\ldots\in\Sigma^{\Bbb N}$.
Let $C$ be a closed subset of $\Bbb R$.
If the OSC is satisfied, then
 $$
 \dim_{\Haus}
 \pi
 \Bigg\{
 \bold i\in\Sigma^{\Bbb N}
 \,\Bigg|\,
 \,\underset{n}\to\acc
 \frac{1}{n}\sum_{k=0}^{n-1}f(S^{k}\bold i)
\subseteq
C
 \Bigg\}
 =
 \sup_{\alpha\in C}
 \,
 \sup
 \Sb
 \mu\in\Cal P_{S}(\Sigma^{\Bbb N})\\
 \int f\,d\mu=\alpha
 \endSb
 -
 \frac{h(\mu)}{\int\Lambda\,d\mu}\,. 
 $$
In particular, if the OSC is satisfied and $\alpha\in\Bbb R$, then we have 
 $$
  \dim_{\Haus}
 \pi
 \Bigg\{
 \bold i\in\Sigma^{\Bbb N}
 \,\Bigg|\,
 \lim_{n}
 \frac{1}{n}\sum_{k=0}^{n-1}f(S^{k}\bold i)
 =
 \alpha
 \Bigg\}
 =
 \sup
 \Sb
 \mu\in\Cal P_{S}(\Sigma^{\Bbb N})\\
 \int f\,d\mu=\alpha
 \endSb
 -
 \frac{h(\mu)}{\int\Lambda\,d\mu}\,.
 \qquad
 $$ 
\endproclaim

\bigskip

\noindent
As a third 
 application of 
Theorem 2.1
we obtain a 
zeta-function
whose
abscissa of convergence equals the 
multifractal
spectrum
$F_{f}^{\erg}$
of ergodic Birkhoff averages of a Lipschitz function
$f$.
This is the content of the next theorem.

\bigskip

\proclaim{Theorem 3.7.
Multifractal zeta-functinons for
multifractal spectra of of ergodic Birkhoff averages}
%Let $f_{1},\ldots,f_{M}:\Sigma^{\Bbb N}\to\Bbb R$ be Lipschitz functions.
Let $f:\Sigma^{\Bbb N}\to\Bbb R$ be a Lipschitz function.
 
For $\bold i\in\Sigma^{*}$, let
 $$
 s_{\bold i}
 =
 \sup_{\bold u\in\Sigma^{\Bbb N}}
 |DS_{\bold i}(\pi\bold u)|\,.
 $$
 and write
 $\overline{\bold i}=\bold i\bold i\bold i\ldots\in\Sigma^{\Bbb N}$.
For a closed set $C\subseteq\Bbb R^{M}$,
we define the self-similar multifractal zeta-function of $f$
by
 $$
 \zeta_{C}^{\erg}(s;r)
 =
 \sum
 \Sb
 \bold i\\
 {}\\
 \dist
 \big(
 \,
 \frac{1}{|\bold i|}\sum_{k=0}^{|\bold i|-1}f(S^{k}\overline{\bold i})
% \big(
% \frac{1}{|\bold i|}\sum_{k=0}^{|\bold i|-1}f_{1}(S^{k}\overline{\bold i})
%,\ldots,
%\frac{1}{|\bold i|}\sum_{k=0}^{|\bold i|-1}f_{M}(S^{k}\overline{\bold i})
%\big)
 \,,\,
 C
 \,
 \big)
 \,
 \le 
 \,
 r
 \endSb
 s_{\bold i}^{s}\,,
 $$
and if $\alpha\in\Bbb R$
and $C=\{\alpha\}$
is the singleton consisting of $\alpha$, then
we
 write $ \zeta_{C}(s;r)=\zeta_{\alpha}(s;r)$, 
 i\.e\. 
 we
write
 $$
 \zeta_{\alpha}^{\erg}(s;r)
 =
 \sum
 \Sb
 \bold i\\
 {}\\
 \big|
 \,
 \frac{1}{|\bold i|}\sum_{k=0}^{|\bold i|-1}f(S^{k}\overline{\bold i})
 \,-\,
 \alpha
 \,
 \big|
 \,
 \le 
 \,
 r
 \endSb
 s_{\bold i}^{s}\,.
 $$

Define $\Lambda:\Sigma^{\Bbb N}\to\Bbb R$ by
$\Lambda(\bold i)=\log |DS_{i_{1}}(\pi S\bold i)|$
for
$\bold i=i_{1}i_{2}\ldots\in\Sigma^{\Bbb N}$.
Then the following hold:
\roster
\item"(1)"
We have
$$
 \lim_{r\searrow0}
 \sigma_{\abs}\big(\, \zeta_{C}^{\erg}(\cdot;r)\,\big)
 =
 \sup_{\alpha\in C}
 \,
 \sup
 \Sb
 \mu\in\Cal P_{S}(\Sigma^{\Bbb N})\\
 \int f\,d\mu=\alpha
 \endSb
 -
 \frac{h(\mu)}{\int\Lambda\,d\mu}\,.
 $$
In particular, if $\alpha\in\Bbb R$, then we have
 $$
 \lim_{r\searrow0}
 \sigma_{\abs}\big(\, \zeta_{\alpha}^{\erg}(\cdot;r)\,\big)
 =
  \,
 \sup
 \Sb
 \mu\in\Cal P_{S}(\Sigma^{\Bbb N})\\
 \int f\,d\mu=\alpha
 \endSb
 -
 \frac{h(\mu)}{\int\Lambda\,d\mu}\,.
 \qquad
 $$

\item"(2)"
If the OSC is satisfied, then we have
$$
\align
\quad
 \lim_{r\searrow0}
 \sigma_{\abs}\big(\, \zeta_{C} ^{\erg}(\cdot;r)\,\big)
 &=
 \sup_{\alpha\in C}
 \dim_{\Haus}
 \pi
 \Bigg\{
 \bold i\in\Sigma^{\Bbb N}
 \,\Bigg|\,
 \lim_{n}
 \frac{1}{n}\sum_{k=0}^{n-1}f(S^{k}\bold i)
 =
 \alpha
 \Bigg\}\\ 
 &=
  \quad\,\,\,\,\,
 \dim_{\Haus}
 \pi
 \Bigg\{
 \bold i\in\Sigma^{\Bbb N}
 \,\Bigg|\,
 \,\underset{n}\to\acc
 \frac{1}{n}\sum_{k=0}^{n-1}f(S^{k}\bold i)
\subseteq
C
 \Bigg\}\,.\\
 \endalign
 $$
In particular, if the OSC is satisfied and $\alpha\in\Bbb R$, then we have 
$$
\align
\quad
 \lim_{r\searrow0}
 \sigma_{\abs}\big(\, \zeta_{\alpha} ^{\erg}(\cdot;r)\,\big)
 &=
  \quad\,\,\,\,\,
 \dim_{\Haus}
 \pi
 \Bigg\{
 \bold i\in\Sigma^{\Bbb N}
 \,\Bigg|\,
 \lim_{n}
 \frac{1}{n}\sum_{k=0}^{n-1}f(S^{k}\bold i)
 =
 \alpha
 \Bigg\}\,.
 \endalign
 $$

\endroster
\endproclaim

\bigskip

We will now prove  Theorem 3.7.
Recall  that the function $\Lambda:\Sigma^{\Bbb N}\to\Bbb R$ is defined by
 $$
 \Lambda(\bold i)=\log|DS_{i_{1}}(\pi S\bold i)|
 \tag3.20
 $$
  for $\bold i=i_{1}i_{2}\ldots\in\Sigma^{\Bbb N}$.
It is well-known that $\Lambda$
satisfies Conditions (C1)--(C3) in Section 2.1.
Also, a straight forward calculation shows that
$\sup_{\bold k\in[\bold i]}\exp\sum_{k=0}^{|\bold i|-1}\Lambda S^{k}\bold k
=
\sup_{\bold u\in\Sigma^{\Bbb N}}|DS_{\bold i}(\pi\bold u)|=s_{\bold i}$
for
$\bold i\in\Sigma^{*}$.
Finally, 
define
$U:\Cal P(\Sigma^{\Bbb N})\to\Bbb R^{M}$ by
 $$
 U\mu
 =
\int f\,d\mu\,.
\tag3.21
 $$
and note that if $\bold i\in\Sigma^{*}$, then
 $$
 UL_{|\bold i|}[\bold i]
 =
 \Bigg\{
 \frac{1}{|\bold i|}
 \sum_{k=0}^{|\bold i|-1}f(S^{k}(\bold i\bold u))
  \,\Bigg|\,
 \bold u\in\Sigma^{\Bbb N}
 \Bigg\}\,.
 $$
Hence, for $C\subseteq\Bbb R$ we have
 $$
 \align
 \zeta_{C}^{U,\Lambda}(s;r)
&=
 \sum
 \Sb
 \bold i\\
 {}\\
 UL_{|\bold i|}[\bold i]
  \subseteq
 B(C,r)
 \endSb
 s_{\bold i}^{s}\\ 
&{}\\ 
&=
 \sum
 \Sb
 \bold i\\
 {}\\
  \big\{
  \frac{1}{|\bold i|}
 \sum_{k=0}^{|\bold i|-1}f(S^{k}(\bold i\bold u))
  \,\big|\,
 \bold u\in\Sigma^{\Bbb N}
 \big\}
 \subseteq
 B(C,r)
 \endSb
 s_{\bold i}^{s}\\
&{}\\ 
&=
 \sum
 \Sb
 \bold i\\
 {}\\
 \forall
 \bold u\in\Sigma^{\Bbb N}
 \,\,:\,\,
 \dist
 \big(
 \,
 \frac{1}{|\bold i|}
 \sum_{k=0}^{|\bold i|-1}f(S^{k}(\bold i\bold u))
  \,,\,
 C
 \,
 \big)
 \,
 \le
 \,
 r
 \endSb
 s_{\bold i}^{s}\,.
 \tag3.22
 \endalign 
 $$
 In order to prove Theorem 3.7, we first
 prove  the following  auxiliary result, namely,
 Proposition 3.8.

\bigskip

\proclaim{Proposition 3.8}
Let $U$ and $\Lambda$
be defined by (3.21) and (3.20), respectively.
\roster
\item"(1)"
There is a sequence $(\Delta_{n})_{n}$ with $\Delta_{n}>0$ 
for all $n$ and 
$\Delta_{n}\to 0$ such that
for all closed subsets $C$ of $\Bbb R$
and
for
all
$n\in\Bbb N$, $\bold i\in\Sigma^{n}$ and $\bold u\in\Sigma^{\Bbb N}$, we have
 $$
 \align
 \dist
 \Bigg(
 \,
 \frac{1}{|\bold i|}
 \sum_{k=0}^{|\bold i|-1}f(S^{k}(\bold i\bold u))
 \,,\,
 C
 \,
 \Bigg)
&\le
\dist
\Bigg(
 \,
 \frac{1}{|\bold i|}
 \sum_{k=0}^{|\bold i|-1}f(S^{k}(\overline{\bold i})) \,,\,
 C
 \,
 \Bigg)
\,\,\,\,\,\,
+
\,\,
\Delta_{n}\,,\\
\dist
\Bigg(
 \,
 \frac{1}{|\bold i|}
 \sum_{k=0}^{|\bold i|-1}f(S^{k}(\overline{\bold i})) \,,\,
 C
 \,
 \Bigg)
&\le
 \dist
 \Bigg(
 \,
 \frac{1}{|\bold i|}
 \sum_{k=0}^{|\bold i|-1}f(S^{k}(\bold i\bold u))
 \,,\,
 C
 \,
 \Bigg)
 \,\,
+
\,\,
\Delta_{n}\,.
\endalign
$$
\item"(2)"
We have
 $$
 \lim_{r\searrow0}
 \sigma_{\abs}
 \big(
 \,
 \zeta_{C}^{\erg}(\cdot;r)
 \,
 \big)
 =
 \lim_{r\searrow0}
 \sigma_{\abs}
 \big(
 \,
 \zeta_{C}^{U,\Lambda}(\cdot;r)
 \,
 \big)\,.
$$
\endroster
\endproclaim
\noindent{\it Proof}\newline
\noindent
(1)
Let $\Lip(f)$ denote the Lipschitz constant of $f$. 
It is clear that
for
all $n\in\Bbb N$,
$\bold i\in\Sigma^{n}$ and $\bold u\in\Sigma^{\Bbb N}$, we have
 $$
 \align
 \Bigg|
  \frac{1}{n}
 \sum_{k=0}^{n-1}f(S^{k}(\overline{\bold i})) 
 -
  \frac{1}{n}
 \sum_{k=0}^{n-1}f(S^{k}(\bold i\bold u)) 
 \Bigg|
 &\le
 \frac{1}{n}
 \sum_{k=0}^{n-1}
 |f(S^{k}(\overline{\bold i}))  - f(S^{k}(\bold i\bold u)) |\\
&\le
\Lip(f)
 \frac{1}{n}
 \sum_{k=0}^{n-1}
 \distance_{\gamma}
 \big(
 \,
 S^{k}(\overline{\bold i}),S^{k}(\bold i\bold u)
 \,
 \big) \\ 
&\le
\Lip(f)
 \frac{1}{n}
 \sum_{k=0}^{n-1}\gamma^{k}\\
&\le
\Lip(f)
 \frac{1}{n(1-\gamma)}\,.
 \tag3.23
  \endalign
 $$
It is not difficult to see that the desired result follows from (3.23).

\noindent
(2)
This statement follows from Part (1)
by an argument very similar to the proof of Part (2) and Part (3)
in
Proposition 3.5, and the proof is therefore omitted.
\hfill$\square$

\bigskip

We can now prove Theorem 3.7.

\bigskip

\noindent{\it Proof of Theorem 3.7}\newline
(1) This statement
follows immediately from 
Theorem 2.1 and Proposition 3.8.

\noindent
(2)  This statement
follows immediately from 
Part (1)
using
Theorem 2.2 and Theorem D.
\hfill$\square$

  \bigskip
  \bigskip

 %%%%%%%%%%%%%%%%%%%%%%%%%%%%%%%%%%%%%%%%%
 %%%%%%%%%%%%%%%%%%%%%%%%%%%%%%%%%%%%%%%%%
 %%%%%%%%%%%%%%%%%%%%%%%%%%%%%%%%%%%%%%%%%
 %%%%%%%%%%%%%%%%%%%%%%%%%%%%%%%%%%%%%%%%%
 %%%%%%%%%%%%%%%%%%%%%%%%%%%%%%%%%%%%%%%%% 

\heading
{
4. Preliminary results}
\endheading

The purpose of this short section is to prove Proposition 4.1 establishing various
auxiliary 
results 
needed for the proof of Theorem 2.1.
Let $c_{\min}$ and $c_{\max}$ be the constants from the 
Condition (C2) in Section 2.1
and
 write
  $$
  \aligned  
  s_{\min}
 &=
  e^{c_{\min}}\,,\\
 s_{\max}
 &=
  e^{c_{\max}}\,.
  \endaligned
  \tag4.1
  $$
we can now state and prove Proposition 4.1.
Recall, that for $\bold i\in\Sigma^{n}$, the number $s_{\bold i}$ is defined 
by
$s_{\bold i}
=
\sup_{\bold k\in[\bold i]}\exp\sum_{k=0}^{n-1}\Lambda S^{k}\bold k$,
see Section 2.1.

\bigskip   
  
 \proclaim{Proposition 4.1} 
 Let $c$ be the constant from Condition (C3) in Section 2.1.
 Let $\bold i,\bold j\in\Sigma^{*}$.
 \roster
 \item"(1)" $0<s_{\min}^{|\bold i|}\le s_{\bold i}\le s_{\max}^{|\bold i|}<1$.
 \item"(2)" $s_{\bold i\bold j}\le s_{\bold i}s_{\bold j}\le c s_{\bold i\bold j}$.
 \item"(3)" $s_{\bold i}< s_{\hat{\bold i}}$.
 \item"(4)"
 For $\bold k\in\Sigma^{\Bbb N}$ and a positive integer $n$, we have
 $\exp
  \,\,
  \sum_{k=0}^{n-1}\Lambda S^{k}\bold k
 \le
 s_{\bold k|n}
 \le
 c
 \,
 \exp
 \,\,
 \sum_{k=0}^{n-1}\Lambda S^{k}\bold k$.
 \item"(5)"
 For $\bold k\in\Sigma^{\Bbb N}$ and a real number $\alpha$, 
 the following two statements are equivalent:

${}$\,\,{\rm (i)}
 \,\,
 $\dsize
 \frac{1}{n}\sum_{k=0}^{n-1}\Lambda S^{k}\bold k\to\alpha$.

{\rm (ii)}
 \,\,
 $\dsize
 \frac{1}{n}\log s_{\bold k|n}\to\alpha$.

 \endroster
 \endproclaim 
 \noindent{\it Proof}\newline
 \noindent
 Statements (1), (2) and (4) follow easily from the definitions.
 Statement (3) follows from (1) and (2),
 and statement (5) follows from (4).
 \hfill$\square$

  \bigskip
  \bigskip

 %%%%%%%%%%%%%%%%%%%%%%%%%%%%%%%%%%%%%%%%%
 %%%%%%%%%%%%%%%%%%%%%%%%%%%%%%%%%%%%%%%%%
 %%%%%%%%%%%%%%%%%%%%%%%%%%%%%%%%%%%%%%%%%
 %%%%%%%%%%%%%%%%%%%%%%%%%%%%%%%%%%%%%%%%%
 %%%%%%%%%%%%%%%%%%%%%%%%%%%%%%%%%%%%%%%%% 

\heading
{
5. Proof of inequality (2.1)
% $$
%    \limsup_{r\searrow 0}
% \,\,
% \sigma_{\abs}\big(\,\zeta_{C}^{U,\Lambda}(\cdot;r)\,\big)
%\le
%\,\,\,\,
% \sup
%  \Sb
%  \mu\in\Cal P_{S}(\Sigma^{\Bbb N})\\
%  {}\\
%  U\mu\in C
%  \endSb
%  \,\,\,\,
% -\frac{h(\mu)}{\int \Lambda\,d\mu}\
% $$
}
\endheading

The purpose of this section
is to prove Theorem 5.5
providing a proof of inequality (2.1).
The proof of (2.1) is based on
results from large deviation theory. 
In particular, 
we need Varadhan's [Va]
large deviation theorem (Theorem 5.1.(i) below),
and a non-trivial application of this (namely Theorem 5.1.(ii) below) 
providing first order asymptotics of certain
\lq\lq Boltzmann distributions".

\bigskip

\proclaim{Definition}
Let $X$
be a complete separable metric space
and
let $(P_{n})_{n}$ be a sequence of probability measures on $X$. 
Let
$(a_{n})_{n}$ be a sequence of positive numbers with
$a_{n}\to\infty$ 
and let
$I:X\to[0,\infty]$ be a lower semicontinuous function
with compact level sets.
The sequence $(P_{n})_{n}$ is said to have the large deviation property
with constants $(a_{n})_{n}$ and rate function $I$ if the following 
two condistions
hold:
\roster
\item"(i)" For each closed subset $K$ of $X$, we have
 $$
 \limsup_{n}\frac{1}{a_{n}}\log P_{n}(K)\le-\inf_{x\in K}I(x)\,;
 $$
\item"(ii)" For each open subset $G$ of $X$, we have
 $$
 \liminf_{n}\frac{1}{a_{n}}\log P_{n}(G)\ge-\inf_{x\in G}I(x)\,.
 $$
\endroster
\endproclaim

\bigskip

\proclaim{Theorem 5.1}
Let $X$ be a complete separable metric space
and
let $(P_{n})_{n}$ be a sequence of probability measures on $X$.
Assume that the sequence $(P_{n})_{n}$
has the large deviation property
with constants $(a_{n})_{n}$ and rate function $I$.
Let $F:X\to\Bbb R$ be a continuous function
satisfying the following two conditions:
\roster
\item"(i)"
For all $n$, we have
 $$
 \int\exp(a_{n}F)\,dP_{n}
 <
 \infty\,.
 $$
\item"(ii)"
We have
 $$
 \lim_{M\to\infty}\,\,
 \limsup_{n}\,\,
 \frac{1}{a_{n}}
 \log
 \int_{\{M\le F\}}
 \exp(a_{n}F)\,dP_{n}
 =
 -\infty\,.
 $$
 \endroster
(Observe that the  Conditions (i)--(ii)
 are satisfied if $F$ is bounded.)
Then the following statements hold.
\roster
\item"(1)" We have
 $$
 \lim_{n}\,\,
 \frac{1}{a_{n}}
 \log
 \int
 \exp(a_{n}F)\,dP_{n}
 =
 -\inf_{x\in X}(I(x)-F(x))\,.
 $$
\item"(2)"
For each $n$ define a probability measure $Q_{n}$ on $X$ by
 $$
 Q_{n}(E)
 =
 \,\,
 \frac
 {\int_{E}\exp(a_{n}F)\,dP_{n}}
 {\int\exp(a_{n}F)\,dP_{n}}\,.
 $$
Then the sequence $(Q_{n})_{n}$
has the large deviation property with constants $(a_{n})_{n}$
and rate function
$(I-F)-\inf_{x\in X}(I(x)-F(x))$.
\endroster
\endproclaim 
\noindent{\it Proof}\newline
\noindent 
Statement (1) follows from [El, Theorem II.7.1] or [DeZe, Theorem 4.3.1], and
statement (2) follows from [El, Theorem II.7.2].
\hfill$\square$

 \bigskip
 
 \noindent
 Using Theorem 5.1
 we first establish
 the following auxiliary result.
 
  \bigskip

\proclaim{Theorem 5.2}
Let $X$ be a metric space and let $U:\Cal P(\Sigma^{\Bbb N})\to X$ be 
continuous with respect to the weak topology.
Let $C\subseteq X$ be a closed subset of $X$ and $r>0$.

If $t\in\Bbb R$, then
 $$
 \limsup_{n}
 \frac{1}{n}
 \log
 \sum
   \Sb
   |\bold i|=n\\
   {}\\
   UL_{n}[\bold i]\subseteq B(C,r)
	 \endSb
 s_{\bold i}^{t} 
 \le
 \sup
  \Sb
	\mu\in\Cal P_{S}(\Sigma^{\Bbb N})\\
	{}\\
	U\mu\in B(C,r)
	\endSb
 \left(t\int\Lambda\,d\mu+h(\mu)\right)\,.
 $$
\endproclaim

\noindent{\it Proof}\newline
\noindent
We start by introducing some notation.
If $\bold i\in\Sigma^{*}$, then we define 
$\overline{\bold i}\in\Sigma^{\Bbb N}$ by
$\overline{\bold i}=\bold i\bold i\ldots$.
We also
define $M_{n}:\Sigma^{\Bbb N}\to\Cal P_{S}(\Sigma^{\Bbb N})$ by
 $$
 \align
 M_{n}\bold i
&=
 L_{n}\left(\,\overline{\bold i|n}\,\right)\\
&= 
 \frac{1}{n}\sum_{k=0}^{n-1}\delta_{S^{k}(\,\overline{\bold i|n}\,)}
 \endalign
 $$
for $\bold i\in\Sigma^{\Bbb N}$;
recall, that
the map
$L_{n}:\Sigma^{\Bbb N}\to\Cal P(\Sigma^{\Bbb N})$ is defined in Section 2.
Furthermore,
note that
if
$\bold i\in\Sigma^{\Bbb N}$,
then
$M_{n}\bold i$ is shift invariant,
i\.e\.
$M_{n}$ 
maps 
$\Sigma^{\Bbb N}$ into $\Cal P_{S}(\Sigma^{\Bbb N})$
as claimed.
Next, let $P$ denote the probability measure on $\Sigma^{\Bbb N}$ given by
%$P={\X^{}}^{\infty}_{n=1}(\sum_{i=1}^{N}\frac{1}{N}\delta_{i})$.
 $$
 P={\underset{\Bbb N}\to \X}\,\,\sum_{i=1}^{N}\frac{1}{N}\delta_{i}\,.
 $$
Finally,
we define $F:\Cal P_{S}(\Sigma^{\Bbb N})\to\Bbb R$ by
 $$
 F(\mu)=t\int\Lambda\,d\mu\,.
 $$
Observe that since $\Lambda$ is bounded, 
i\.e\. $\|\Lambda\|_{\infty}<\infty$,
we conclude that
$\|F\|_{\infty}
=
|t|\,\|\Lambda\|_{\infty}<\infty$.
 Also, for a positive integer $n$, define
probability measures
$P_{n},Q_{n}\in\Cal P(\Cal P_{S}(\Sigma^{\Bbb N}))$ by
 $$
 \align
 P_{n}
&=
 P\circ M_{n}^{-1}\,,\\
 Q_{n}(E)
&=
 \frac
 {\int_{E}\exp(nF)\,dP_{n}}
 {\int\exp(nF)\,dP_{n}}
 \quad
 \text{for $E\subseteq\Cal P_{S}(\Sigma^{\Bbb N})$.}
 \endalign
 $$
We now prove the following two claims.

 \bigskip

{\it Claim 1. We have
 $$
 \sum
   \Sb
  |\bold k|=n\\
  {}\\
  UL_{n}[\bold k]\subseteq B(C,r)
	 \endSb
 s_{\bold k}^{t}
\le
 \sum
   \Sb
  |\bold k|=n\\
  {}\\
  UM_{n}[\bold k]\subseteq B(C,r)
	 \endSb
 s_{\bold k}^{t}\,.
 $$
} 
\noindent{\it Proof of Claim 1.} 
Observe that
if $|\bold k|=n$, then
$M_{n}[\bold k]
=
\{M_{n}(\bold k\bold l)\,|\,\bold l\in\Sigma^{\Bbb N}\}
=
\{L_{n}(\,\overline{(\bold k\bold l)|n}\,)\,|\,\bold l\in\Sigma^{\Bbb N}\}
=
\{L_{n}\overline{\bold k}\,|\,\bold l\in\Sigma^{\Bbb N}\}
=
\{L_{n}\overline{\bold k}\}
\subseteq
L_{n}[\bold k]$.
The desired result follows immediately from this inclusion.
This proves Claim 1.

\bigskip

 {\it Claim 2. We have
 $$
 \sum
   \Sb
  |\bold k|=n\\
  {}\\
  UM_{n}[\bold k]\subseteq B(C,r)
	 \endSb
 s_{\bold k}^{t}
 \le
 N^{n}
 \int
 \limits_{
 \big\{
 \bold j\in\Sigma^{\Bbb N}
 \,\big|\,
 U M_{n}[\bold j|n]\subseteq B(C, r)
 \big\}
 }
s_{\bold i|n}^{t}
 \,dP(\bold i)\,.
 $$
} 
\noindent{\it Proof of Claim 2.} 
It  follows that
 $$
 \align
 \int
 \limits_{
 \big\{
 \bold j\in\Sigma^{\Bbb N}
 \,\big|\,
 U M_{n}[\bold j|n]\subseteq B(C, r)
 \big\}
 }
&s_{\bold i|n}^{t}
 \,dP(\bold i)\\
&=
\qquad\!
 \sum_{|\bold k|=n}
 \,\,\,\,
 \int\limits_{
 [\bold k]
 \,\,\cap\,\,
 \big\{
 \bold j\in\Sigma^{\Bbb N}
 \,\big|\,
 U M_{n}[\bold j|n]\subseteq B(C, r)
  \big\}
 }
 s_{\bold i|n}^{t}
 \,dP(\bold i)\\
\allowdisplaybreak 
&=
\qquad\!
 \sum_{|\bold k|=n}
 \qquad
 s_{\bold k}^{t}
 \,\,
 P
 \Big(
 \,
  [\bold k]
 \,\,\cap\,\,
 \Big\{
 \bold j\in\Sigma^{\Bbb N}
 \,\Big|\,
 U M_{n}[\bold j|n]\subseteq B(C, r)
 \Big\}
 \,
 \Big)\\
\allowdisplaybreak 
&\ge
 \sum
   \Sb
   |\bold k|=n\\
   {}\\
   U M_{n}[\bold k]\subseteq B(C, r)
	 \endSb
 s_{\bold k}^{t}
 \,\,
 P
 \Big(
 \,
  [\bold k]
 \,\,\cap\,\,
 \Big\{
 \bold j\in\Sigma^{\Bbb N}
 \,\Big|\,
 U M_{n}[\bold j|n]\subseteq B(C, r)
  \Big\}
 \,
 \Big)\,.\\
& \tag5.1
 \endalign
 $$
However,
for 
$\bold k$ with
 $|\bold k|=n$
 and
  $ U M_{n}[\bold k]\subseteq B(C, r)$, it is clear that
  $[\bold k]\subseteq
 \{
 \bold j\in\Sigma^{\Bbb N}
 \,|\,
  U M_{n}[\bold j|n]\subseteq B(C, r)
 \}$,
 whence
 $[\bold k]\cap
 \{
 \bold j\in\Sigma^{\Bbb N}
 \,|\,
  U M_{n}[\bold j|n]\subseteq B(C, r)
 \}
 =
 [\bold k]$.
 This and (5.1) now imply that
 $$
 \align
  \int
 \limits_{
 \big\{
 \bold j\in\Sigma^{\Bbb N}
 \,\big|\,
 U M_{n}[\bold j|n]\subseteq B(C, r)
 \big\}
 }
&s_{\bold i|n}^{t}
 \,dP(\bold i)\\
&\ge
 \sum
   \Sb
   |\bold k|=n\\
   {}\\
   U M_{n}[\bold k]\subseteq B(C, r)
	 \endSb
 s_{\bold k}^{t}
 \,\,
 P
 \Big(
 \,
  [\bold k]
 \,\,\cap\,\,
 \Big\{
 \bold j\in\Sigma^{\Bbb N}
 \,\Big|\,
 U M_{n}[\bold j|n]\subseteq B(C, r)
  \Big\}
 \,
 \Big)\\ 
&=
 \sum
   \Sb
   |\bold k|=n\\
   {}\\
 U M_{n}[\bold k]\subseteq B(C, r)
	 \endSb
 s_{\bold k}^{t}
 \,\,
 P\big(\,
 [\bold k]
 \,\big)\\
\allowdisplaybreak 
&=
 \sum
   \Sb
  |\bold k|=n\\
  {}\\
   \dist(UL\bold k,C\,)\le r
	 \endSb
 s_{\bold k}^{t}\,\frac{1}{N^{n}}\,.
 \endalign
 $$

 \noindent
Hence
 $$
 \align
 \sum
   \Sb
  |\bold k|=n\\
  {}\\
   U M_{n}[\bold k]\subseteq B(C, r)
	 \endSb
 s_{\bold k}^{t}
&\le
 N^{n}
 \int
 \limits_{
 \big\{
 \bold j\in\Sigma^{\Bbb N}
 \,\big|\,
 U M_{n}[\bold j|n]\subseteq B(C, r)
 \big\}
 }
s_{\bold i|n}^{t}
 \,dP(\bold i)\,.
 \endalign
 $$
This completes the proof of Claim 2.

\bigskip

Combining Claim 1 and Claim 2 shows that
  $$
 \align
 \sum
   \Sb
  |\bold k|=n\\
  {}\\
  UL_{n}[\bold k]\subseteq B(C,r)
	 \endSb
 s_{\bold k}^{t}
 &\le
 \sum
   \Sb
  |\bold k|=n\\
  {}\\
  UM_{n}[\bold k]\subseteq B(C,r)
	 \endSb
 s_{\bold k}^{t}\\
 &\le
 N^{n}
 \int
 \limits_{
 \big\{
 \bold j\in\Sigma^{\Bbb N}
 \,\big|\,
 U M_{n}[\bold j|n]\subseteq B(C, r)
 \big\}
 }
s_{\bold i|n}^{t}
 \,dP(\bold i)\,.
 \tag5.2
 \endalign
 $$ 
 Let $c$ be the constant from Condition (C3) in Section 2.1,
 and notice that it follows from Proposition 4.1
that
if $\bold i\in\Sigma^{\Bbb N}$ and $n$ is a positive integer, then
we have
 $s_{\bold i|n}^{t}
 \le
c^{|t|}
\exp(\,t\sum_{k=0}^{n-1}\Lambda S^{k}(\,\overline{\bold i|n}\,)\,)$.
We  conclude from this and (5.2) that
 $$
 \align
  \sum
   \Sb
  |\bold k|=n\\
  {}\\
   U L_{n}[\bold k]\subseteq B(C, r)
	 \endSb
 s_{\bold k}^{t}
&\le
 N^{n}
 \int
 \limits_{
 \big\{
 \bold j\in\Sigma^{\Bbb N}
 \,\big|\,
 U M_{n}[\bold j|n]\subseteq B(C, r)
 \big\}
 }
s_{\bold i|n}^{t}
 \,dP(\bold i)\\
&\le 
 c^{|t|}
 \,
 N^{n}
 \int\limits_{
 \big\{
 \bold j\in\Sigma^{\Bbb N}
 \,\big|\,
U M_{n}[\bold j|n]\subseteq B(C, r)
 \big\}
 }
 \exp
 \Bigg(
 t\sum_{k=0}^{n-1}\Lambda S^{k}\left(\,\overline{\bold i|n}\,\right)
 \Bigg)
 \,dP(\bold i)\\
\allowdisplaybreak 
&=
 c^{|t|}
 \,
 N^{n}
 \int\limits_{
 \big\{
 \bold j\in\Sigma^{\Bbb N}
 \,\big|\,
U M_{n}[\bold j|n]\subseteq B(C, r)
 \big\}
  }
 \exp
 \Bigg(
 tn\int\Lambda\,d(M_{n}\bold i)
 \Bigg)
 \,dP(\bold i)\\
\allowdisplaybreak 
&=
 c^{|t|}
 \,
 N^{n}
 \int\limits_{
 \big\{
 \bold j\in\Sigma^{\Bbb N}
 \,\big|\,
 U M_{n}[\bold j|n]\subseteq B(C, r)
 \big\}
 }
 \exp\left(n F(M_{n}\bold i)\right)
 \,dP(\bold i)\,.
 \tag5.3
 \endalign
 $$
Noticing that
$\{\bold j\in\Sigma^{\Bbb N}
 \,|\,
 U M_{n}[\bold j|n]\subseteq B(C, r)
 \}
 \subseteq
 \{\bold j\in\Sigma^{\Bbb N}
 \,|\,
 U M_{n}\bold j\subseteq B(C, r)
 \}
 =
 \{UM_{n}\in B(C,r)\}$, we now deduce from (5.3) that
 $$
 \align
  \sum
   \Sb
  |\bold k|=n\\
  {}\\
   U L_{n}[\bold k]\subseteq B(C, r)
	 \endSb
 s_{\bold k}^{t}
&\le
 c^{|t|}
 \,
 N^{n}
 \int\limits_{
 \big\{
 \bold j\in\Sigma^{\Bbb N}
 \,\big|\,
 U M_{n}[\bold j|n]\subseteq B(C, r)
 \big\}
 }
 \exp\left(n F(M_{n}\bold i)\right)
 \,dP(\bold i)\\ 
&\le
 c^{|t|}
 \,
 N^{n}
 \int\limits_{
 \big\{U M_{n}\in B(C,r)\big\}
 }
 \exp\left(nF(M_{n}\bold i)\right)
 \,dP(\bold i)\\
\allowdisplaybreak 
&=
c^{|t|}
 \,
 N^{n}
 \int\limits_{
 \big\{U\in B(C,r)\big\}
 }
 \exp\left(nF\right)
 \,dP_{n}\\
\allowdisplaybreak 
&=
 c^{|t|}
 \,
 N^{n}
 \,\,
 Q_{n}\Big(\,\Big\{U\in B(C,r)\Big\}\,\Big)
 \,\,
 \int\exp\left(nF\right)\,dP_{n}\,.
 \tag5.4
 \endalign
 $$
It follows immediately from (5.4)
that 
 $$
 \align
 \limsup_{n}
 \frac{1}{n}
 \log
 \sum
   \Sb
   |\bold i|=n\\
   {}\\
   U L_{n}[\bold i]\subseteq B(C, r)
	 \endSb
 s_{\bold i}^{t} 
&\le
 \log N
 \,+\,
 \limsup_{n}
 \frac{1}{n}
 \log Q_{n}\Big(\,\Big\{U\in B(C,r) \Big\}\,\Big)\\
&\qquad\qquad
  \qquad\qquad
 \,+\,
 \limsup_{n}
 \frac{1}{n}
 \log 
 \int
 \exp
 \left(
 nF
 \right)
 \,dP_{n}\,.
 \tag5.5
 \endalign
 $$

Next, we observe that it follows from [El]
that the sequence 
$(P_{n}=P\circ M_{n}^{-1})_{n}
\subseteq
\Cal P\big(\,\Cal P_{S}(\Sigma^{\Bbb N})\,\big)$ 
has the large deviation property with respect to
the sequence
$(n)_{n}$ and rate function 
$I:\Cal P_{S}(\Sigma^{\Bbb N})\to\Bbb R$ given by
$I(\mu)=\log N-h(\mu)$.
We therefore 
conclude from Part (1) of Theorem 5.1 that
 $$
  \limsup_{n}
 \frac{1}{n}
 \log 
 \int
 \exp
 \left(
 nF
 \right)
 \,dP_{n}
 =
 -
  \inf_{\nu\in\Cal P_{S}(\Sigma^{\Bbb N})}(I(\nu)-F(\nu))\,.
  \tag5.6
 $$
Also,
since the sequence 
$(P_{n}=P\circ M_{n}^{-1})_{n}
\subseteq
\Cal P\big(\,\Cal P_{S}(\Sigma^{\Bbb N})\,\big)$ 
has the large deviation property with respect to
the sequence
$(n)_{n}$ and rate function 
$I:\Cal P_{S}(\Sigma^{\Bbb N})\to\Bbb R$ given by
$I(\mu)=\log N-h(\mu)$, 
we conclude from Part (2) of Theorem 5.1
that the sequence $(Q_{n})_{n}$ has the large deviation property with 
respect to 
the sequence
$(n)_{n}$ and rate function
$(I-F)-\inf_{\nu\in\Cal P_{S}(\Sigma^{\Bbb N})}(I(\nu)-F(\nu))$.
As the set
$\{U\in B(C,r)\}
=
U^{-1}( B(C,r))$
is closed,
 it therefore follows from the large deviation property that
 $$
 \align
 \limsup_{n}
 \frac{1}{n}
&\log Q_{n}\Big(\,\Big\{U\in B(C,r) \Big\}\,\Big)\\
&\le
-
 \inf
  \Sb
	\mu\in\Cal P_{S}(\Sigma^{\Bbb N})\\
	{}\\
	U\mu\in B(C,r)
	\endSb
 \Bigg(
 (I(\mu)-F(\mu))
 -
 \inf_{\nu\in\Cal P_{S}(\Sigma^{\Bbb N})}
 (I(\nu)-F(\nu))
 \Bigg)\,.
 \tag5.7
 \endalign
 $$

Combining (5.5). (5.6) and (5.7) now yields
 $$
 \align
 \limsup_{n}
 \frac{1}{n}
 \log
 \sum
   \Sb
   |\bold i|=n\\
   {}\\
   U L_{n}[\bold i]\subseteq B(C, r)
	 \endSb
 s_{\bold i}^{t} 
&\le
 \log N
 \,+\,
 \limsup_{n}
 \frac{1}{n}
 \log Q_{n}\Big(\,\Big\{U\in B(C,r) \Big\}\,\Big)\\
&\qquad\qquad
  \qquad\qquad
 \,+\,
 \limsup_{n}
 \frac{1}{n}
 \log 
 \int
 \exp
 \left(
 nF
 \right)
 \,dP_{n}\\
&\le
 \log N\\
&\qquad
-
 \inf
  \Sb
	\mu\in\Cal P_{S}(\Sigma^{\Bbb N})\\
	{}\\
	U\mu\in B(C,r)
	\endSb
 \Bigg(
 (I(\mu)-F(\mu))
 -
 \inf_{\nu\in\Cal P_{S}(\Sigma^{\Bbb N})}
 (I(\nu)-F(\nu))
 \Bigg)\\
&\qquad\qquad
  \qquad\qquad
 \,-\,
 \inf_{\nu\in\Cal P_{S}(\Sigma^{\Bbb N})}(I(\nu)-F(\nu))\\
&=
 \log N
 \,+\,
 \sup
  \Sb
	\mu\in\Cal P_{S}(\Sigma^{\Bbb N})\\
	{}\\
	U\mu\in B(C,r)
	\endSb
 (F(\mu)-I(\mu))\\
&=
 \sup
  \Sb
	\mu\in\Cal P_{S}(\Sigma^{\Bbb N})\\
	{}\\
	U\mu\in B(C,r)
	\endSb
 \left(t\int\Lambda\,d\mu+h(\mu)\right)\,. 
 \endalign
 $$
This completes the proof.
\hfill$\square$

\bigskip

\noindent
We will now use Theorem 5.2
to prove Theorem 5.5
providing a proof of inequality (2.1).
However, we first prove two small lemmas.

\bigskip

\proclaim{Lemma 5.3}
Let $X$ be  a metric space and 
let $f,g:X\to\Bbb R$ be upper 
semi-continuous functions
with $f,g\ge 0$.
Then $fg$ is upper semi-continuous.
\endproclaim
\noindent{\it Proof}\newline
Since $f$ and $g$ are upper semi-continuous with
$f,g\ge 0$,
this result follows easily from the definition of upper semi-continuity, and the 
proof is therefore omitted.
\hfill$\square$

\bigskip

\proclaim{Lemma 5.4}
Let $X$ be  a metric space and let $\Phi:X\to\Bbb R$ be an upper 
semi-continuous function.
Let $K_{1},K_{2},\ldots\subseteq X$ be non-empty compact subsets of 
$X$ with
$K_{1}\supseteq K_{2}\supseteq\ldots$. Then
 $$
 \inf_{n}\sup_{x\in K_{n}}\Phi(x)
 =
 \sup_{x\in\cap_{n}K_{n}}\Phi(x)\,.
 $$
\endproclaim
\noindent{\it Proof}\newline
First note that it is clear that
$\inf_{n}\sup_{x\in K_{n}}\Phi(x)
 \ge
 \sup_{x\in\cap_{n}K_{n}}\Phi(x)$.
 We will now prove the reverse inequality, namely,
$\inf_{n}\sup_{x\in K_{n}}\Phi(x)
 \le
 \sup_{x\in\cap_{n}K_{n}}\Phi(x)$. 
Let $\varepsilon>0$.
For each $n$, we can choose $x_{n}\in K_{n}$ such that
$\Phi(x_{n})\ge \sup_{x\in K_{n}}\Phi(x)-\varepsilon$.
Next, since $K_{n}$ is compact for all $n$ and
$K_{1}\supseteq K_{2}\supseteq\ldots$, we can find a subsequence $(x_{n_{k}})_{k}$ and a 
point $x_{0}\in\cap_{n}K_{n}$ such that $x_{n_{k}}\to x_{0}$.
Also, since
$K_{n_{1}}\supseteq K_{n_{2}}\supseteq\ldots$, we conclude that
$\sup_{x\in K_{n_{1}}}\Phi(x)
\ge
\sup_{x\in K_{n_{2}}}\Phi(x)
\ge
\ldots$, 
whence
$\inf_{k}\sup_{x\in K_{n_{k}}}\Phi(x)
=
\limsup_{k}\sup_{x\in K_{n_{k}}}\Phi(x)$.
This implies that
$\inf_{n}\sup_{x\in K_{n}}\Phi(x)
\le
\inf_{k}\sup_{x\in K_{n_{k}}}\Phi(x)
=
\limsup_{k}\sup_{x\in K_{n_{k}}}\Phi(x)
\le
\limsup_{k}\Phi(x_{n_{k}})+\varepsilon$.
However, since 
$x_{n_{k}}\to x_{0}$, we deduce from the 
upper semi-continuity of the function $\Phi$, that
$\limsup_{k}\Phi(x_{n_{k}})\le \Phi(x_{0})$.
Consequently
$\inf_{n}\sup_{x\in K_{n}}\Phi(x)
\le
\limsup_{k}\Phi(x_{n_{k}})+\varepsilon
\le
\Phi(x_{0})+\varepsilon
\le
 \sup_{x\in\cap_{n}K_{n}}\Phi(x)+\varepsilon$.
 Finally, letting $\varepsilon\searrow 0$ gives the desired result.
\hfill$\square$

\bigskip

\noindent
We can now state and prove Theorem 5.5.

\bigskip

\proclaim{Theorem 5.5}
Let $X$ be a metric space and let $U:\Cal P(\Sigma^{\Bbb N})\to X$ be 
continuous with respect to the weak topology.
Let $C\subseteq X$ be a closed subset of $X$ and $r>0$.

\roster
\item"(1)"
We have
 $$
 \sigma_{\abs}\big(\,\zeta_{C}^{U,\Lambda}(\cdot;r)\,\big)
\le
 \sup
  \Sb
  \mu\in\Cal P_{S}(\Sigma^{\Bbb N})\\
  {}\\
  U\mu\in B(C,r)
  \endSb
 -\frac{h(\mu)}{\int \Lambda\,d\mu}\,.
 $$
\item"(2)"
We have
 $$
 \limsup_{r\searrow 0}
 \,\,
 \sigma_{\abs}\big(\,\zeta_{C}^{U,\Lambda}(\cdot;r)\,\big)
\le
\,\,\,\,
 \sup
  \Sb
  \mu\in\Cal P_{S}(\Sigma^{\Bbb N})\\
  {}\\
  U\mu\in C
  \endSb
  \,\,\,\,
 -\frac{h(\mu)}{\int \Lambda\,d\mu}\,.
 \qquad\qquad
 $$

\endroster 
\endproclaim

\noindent{\it Proof}\newline

\noindent
(1)
For brevity write
 $$
 u
 =
  \sup
  \Sb
  \mu\in\Cal P_{S}(\Sigma^{\Bbb N})\\
  {}\\
  U\mu\in B(C,r)
  \endSb
 -\frac{h(\mu)}{\int \Lambda\,d\mu}\,.
 $$
We must now prove that if $t>u$, then 
  $$
 \sum
  \Sb
  \bold i\\
  {}\\
  U L_{|\bold i|}[\bold i]\subseteq B(C, r)
    \endSb
 s_{\bold i}^{t}
 <\infty\,.
 $$

Let $t>u$ and write $\varepsilon=\frac{t-u}{3}>0$.
It follows from the 
definition
of $u$ that
 if 
$\mu\in\Cal P_{S}(\Sigma^{\Bbb N})$
with $U\mu\in B(C,r)$,
then
 we have
$-\frac{h(\mu)}{\int\Lambda\,d\mu}<u+\varepsilon=(u+2\varepsilon)-\varepsilon$, 
whence
$-h(\mu)
>
(u+2\varepsilon)\int\Lambda\,d\mu-\varepsilon\int\Lambda\,d\mu$
where
we have used the fact that
$\int\Lambda\,d\mu<0$
because 
$\Lambda<0$. This implies that
 if 
$\mu\in\Cal P_{S}(\Sigma^{\Bbb N})$
with $U\mu\in B(C,r)$,
then
 $$
 \align
 (u+2\varepsilon)\int\Lambda\,d\mu+h(\mu)
&\le
 \varepsilon\int\Lambda\,d\mu\\
&\le
 \varepsilon c_{\max}\\
&=
-\varepsilon\,|c_{\max}|\,.
 \endalign
 $$
We deduce from this inequality and Theorem 5.2 that
 $$
 \align
 \limsup_{n}
\frac{1}{n}
 \log
&\sum
   \Sb
   |\bold i|=n\\
   {}\\
   U L_{n}[\bold i]\subseteq B(C, r)
	 \endSb
 s_{\bold i}^{t}\\
&=
 \limsup_{n}
 \frac{1}{n}
 \log
 \sum
   \Sb
   |\bold i|=n\\
   {}\\
  U L_{n}[\bold i]\subseteq B(C, r)
	 \endSb
 s_{\bold i}^{u+3\varepsilon} \\
&\le
 \limsup_{n}
 \frac{1}{n}
 \log
 \sum
   \Sb
   |\bold i|=n\\
   {}\\
  U L_{n}[\bold i]\subseteq B(C, r)
	 \endSb
 s_{\bold i}^{u+2\varepsilon} \\
&\le
  \sup
  \Sb
  \mu\in\Cal P_{S}(\Sigma^{\Bbb N})\\
  {}\\
  U\mu\in B(C,r)
  \endSb
  \,\left(\,(u+2\varepsilon)\int\Lambda\,d\mu+h(\mu)\,\right)
  \qquad\qquad
  \text{[by Theorem 5.2]}\\
 &\le
	-\varepsilon \,|c_{\max}|\\
&<
 -\tfrac{1}{2}\varepsilon \,|c_{\max}|\,.	
 \tag5.8
 \endalign	
 $$
Inequality (5.8) shows
that
there is an  integer $N_{0}$ such that 
$ \frac{1}{n}
 \log
 \sum_{
   |\bold i|=n\,,\,
   U L_{n}[\bold i]\subseteq B(C, r)
   }
 s_{\bold i}^{t}
 \le
  -\tfrac{1}{2}\varepsilon \,|c_{\max}|$
for all $n\ge N_{0}$, whence
 $$
  \sum
   \Sb
   |\bold i|=n\\
   {}\\
  U L_{n}[\bold i]\subseteq B(C, r)
	 \endSb
 s_{\bold i}^{t}
 \le
 e^{-\frac{1}{2}\varepsilon \,|c_{\max}|\, n}
 \tag5.9
 $$
for all $n\ge N_{0}$. 
Using (5.9) we now conclude that
 $$
 \align
 \sum
  \Sb
  \bold i\\
  {}\\
 U L_{|\bold i|}[\bold i]\subseteq B(C, r)
  \endSb
 s_{\bold i}^{t}
&=
 \sum_{n<N_{0}}
 \sum
   \Sb
   |\bold i|=n\\
   {}\\
   U L_{n}[\bold i]\subseteq B(C, r)
	 \endSb
 s_{\bold i}^{t}
 \,\,
 +
 \,\,
  \sum_{n\ge N_{0}}
 \sum
   \Sb
   |\bold i|=n\\
   {}\\
   U L_{n}[\bold i]\subseteq B(C, r)
	 \endSb
 s_{\bold i}^{t}\\
&{}\\ 
&\le
 \sum_{n<N_{0}}
 \sum
   \Sb
   |\bold i|=n\\
   {}\\
   U L_{n}[\bold i]\subseteq B(C, r)
	 \endSb
 s_{\bold i}^{t}
 \,\,
 +
 \,\,
  \sum_{n\ge N_{0}}
e^{-\frac{1}{2}\varepsilon \,|c_{\max}|\, n}\\
&{}\\
&<
\infty\,.
\endalign
$$
This completes the proof of (1).

\noindent
(2)
It follows immediately from Part (1) that
 $$
 \limsup_{r\searrow 0}
 \,\,
 \sigma_{\abs}\big(\,\zeta_{C}^{U,\Lambda}(\cdot;r)\,\big)
\le
\limsup_{r\searrow 0}
\,\,
 \sup
  \Sb
  \mu\in\Cal P_{S}(\Sigma^{\Bbb N})\\
  {}\\
  U\mu\in B(C,r)
  \endSb
 -\frac{h(\mu)}{\int \Lambda\,d\mu}\,. 
 \tag5.10
 $$
Also,
the function
$r\to
\sup_{
  \mu\in\Cal P_{S}(\Sigma^{\Bbb N})\,,\,
  U\mu\in B(C,r)
  }
 -\frac{h(\mu)}{\int \Lambda\,d\mu}$
 is clearly increasing, and it therefore follows that
 $$
 \limsup_{r\searrow 0}
\,\,
 \sup
  \Sb
  \mu\in\Cal P_{S}(\Sigma^{\Bbb N})\\
  {}\\
  U\mu\in B(C,r)
  \endSb
 -\frac{h(\mu)}{\int \Lambda\,d\mu}
 =
 \inf_{k}
\,\,
 \sup
  \Sb
  \mu\in\Cal P_{S}(\Sigma^{\Bbb N})\\
  {}\\
  U\mu\in B(C,\frac{1}{k}) 
  \endSb
 -\frac{h(\mu)}{\int \Lambda\,d\mu}\,.
 \tag5.11
 $$

Next,
since
the function $U:\Cal P(\Sigma^{\Bbb N})\to X$ is continuous, we conclude that 
the set
$U^{-1}B(C,\frac{1}{k})$ is closed,
and it therefore follows that the 
set
 $K_{k}=\Cal P_{S}(\Sigma^{\Bbb N})\cap U^{-1}B(C,\frac{1}{k})$
is compact.
Also, since the entropy function $h:\Cal P_{S}(\Sigma^{\Bbb N})\to\Bbb R$
is upper semi-continuous (see [Wa, Theorem 8.2])
with $h\ge 0$
and
the 
function
$f:\Cal P_{S}(\Sigma^{\Bbb N})\to\Bbb R$ given by
$f(\mu)=-\frac{1}{\int\Lambda\,d\mu}$ 
is continuous
(because $\Lambda$ is continuous)
with $f\ge 0$,
we 
conclude from
Lemma 5.3
that the function
$\Phi:\Cal P_{S}(\Sigma^{\Bbb N})\to\Bbb R$ given by
$\Phi(\mu)=f(\mu)h(\mu)=-\frac{h(\mu)}{\int\Lambda\,d\mu}$ 
is upper semi-continuous.
Lemma 5.4 applied to $\Phi$ therefore implies that
 $$
 \align
  \inf_{k}
\,\,
 \sup
  \Sb
  \mu\in\Cal P_{S}(\Sigma^{\Bbb N})\\
  {}\\
 U\mu\in B(C,\frac{1}{k}) 
  \endSb
 -\frac{h(\mu)}{\int \Lambda\,d\mu}
&= 
 \inf_{k}\sup_{\mu\in K_{k}}-\frac{h(\mu)}{\int\Lambda\,d\mu}\\
&=
 \sup_{\mu\in\cap_{k}K_{k}}-\frac{h(\mu)}{\int\Lambda\,d\mu}\,.
 \tag5.12
 \endalign
 $$
However, clearly
$\cap_{k}K_{k}
=
\cap_{k}(\Cal P_{S}(\Sigma^{\Bbb N})\cap U^{-1}B(C,\frac{1}{k}))
=
\Cal P_{S}(\Sigma^{\Bbb N})\cap U^{-1}C$, 
whence
 $$
 \sup_{\mu\in\cap_{k}K_{k}}-\frac{h(\mu)}{\int\Lambda\,d\mu}
 =
  \sup
 \Sb
 \mu\in\Cal P_{S}(\Sigma^{\Bbb N})\\
 {}\\
 U\mu\in C
 \endSb
 -
 \frac{h(\mu)}{\int\Lambda\,d\mu}\,.
 \tag5.13
$$
Combining (5.12) and (5.13) gives
 $$
 \align
  \inf_{k}
\,\,
 \sup
  \Sb
  \mu\in\Cal P_{S}(\Sigma^{\Bbb N})\\
  {}\\
 U\mu\in B(C,\frac{1}{k}) 
  \endSb
 -\frac{h(\mu)}{\int \Lambda\,d\mu}
&=
  \sup
 \Sb
 \mu\in\Cal P_{S}(\Sigma^{\Bbb N})\\
 {}\\
 U\mu\in C
 \endSb
 -
 \frac{h(\mu)}{\int\Lambda\,d\mu}\,.
 \tag5.14
 \endalign
 $$

Finally, the desired result follows by combining (5.10), (5.11) and (5.14).
\hfill$\square$

  \bigskip
  \bigskip

 %%%%%%%%%%%%%%%%%%%%%%%%%%%%%%%%%%%%%%%%%
 %%%%%%%%%%%%%%%%%%%%%%%%%%%%%%%%%%%%%%%%%
 %%%%%%%%%%%%%%%%%%%%%%%%%%%%%%%%%%%%%%%%%
 %%%%%%%%%%%%%%%%%%%%%%%%%%%%%%%%%%%%%%%%%
 %%%%%%%%%%%%%%%%%%%%%%%%%%%%%%%%%%%%%%%%% 

\heading{6. Proof of inequality (2.2)
% $$
%  \sup
%  \Sb
%  \mu\in\Cal P_{S}(\Sigma^{\Bbb N})\\
%  {}\\
%  U\mu\in C
%  \endSb
%  \,\,\,\,
% -\frac{h(\mu)}{\int \Lambda\,d\mu}
% \le
%  \underline f^{U,\Lambda}(C)
%  $$
 }
 \endheading

 The purpose of this section is to prove Theorem
 6.6 providing a proof of inequality (2.2).

We first state and prove a number of auxiliary results.
For $\bold i,\bold j\in\Sigma^{\Bbb N}$ with
with $\bold i\not=\bold j$,
we will write
$\bold i\wedge\bold j$ for the longest common 
prefix of 
$\bold i$ and $\bold j$
(i\.e\.
$
\bold i\wedge\bold j
=
\bold u
$
where $\bold u$ is the unique element in $\Sigma^{*}$
for which there 
are $\bold k,\bold l\in\Sigma^{\Bbb N}$
with
$\bold k=k_{1}k_{2}\ldots$
and
$\bold l=l_{1}l_{2}\ldots$
such that
% $$
% \align
% k_{1}
% &\not= l_{1}\,,\\
% \bold i
%&=
% \bold u\bold k\,,\\
% \bold j
%&=
% \bold u\bold l\,.
% \endalign
% $$ 
$k_{1}
 \not= l_{1}$,
 $\bold i
 =
 \bold u\bold k$
and 
$\bold j
 =
 \bold u\bold l$). 
We will always equip $\Sigma^{\Bbb N}$ with the metric
$\distance_{\Sigma^{\Bbb N}}$ defined by
 $$
 \distance_{\Sigma^{\Bbb N}}(\bold i,\bold j)
 =
 \cases
 0
&\quad 
 \text{if $\bold i=\bold j$;}\\
 s_{\bold i\wedge\bold j}
&\quad 
 \text{if $\bold i\not=\bold j$,}
 \endcases
 \tag6.1
 $$
for $\bold i,\bold j\in\Sigma^{\Bbb N}$. 
In the results below,
we will always
compute
the Hausdorff dimension of a subset of $\Sigma^{\Bbb N}$ 
with respect to the metric $ \distance_{\Sigma^{\Bbb N}}$.
Note that when $\Sigma^{\Bbb N}$ 
is equipped with the
metric 
$ \distance_{\Sigma^{\Bbb N}}$, then
 $$
 \diam[\bold i]
 =
 s_{\bold i}
 \tag6.2
 $$
for all $\bold i\in\Sigma^{*}$.

\bigskip

\proclaim{Lemma 6.1}
Let $(X,\distance)$ be a metric space
and let $U:\Cal P(\Sigma^{\Bbb N})\to X$ be continuous with respect to the weak topology.
Let $C$ be a closed subset of $X$
and $r>0$.
\roster
\item"(1)"
There is a positive integer $M_{r}$
such that
if $k\ge M_{r}$,
$\bold u\in\Sigma^{k}$ and $\bold k,\bold l\in\Sigma^{\Bbb N}$, then
 $$
 \distance
 \big(
 \,
 UL_{k}(\bold u\bold k)
 \,,\,
 UL_{k}(\bold u\bold l)
 \,
 \big)
 \le
 \frac{r}{2}\,.
 $$
 \item"(2)"
There is a positive integer $M_{r}$
such that
if $m\ge M_{r}$,
then
 $$
 \align
 \Big\{
 \bold i\in\Sigma^{\Bbb N}
 \,\Big|\,
 \,\,
&UL_{k}\bold i\in B(C,\tfrac{r}{2})
\,\,\,\,
\text{for all $k\ge m$}
\Big\}\\
&\subseteq
\Big\{
 \bold i\in\Sigma^{\Bbb N}
 \,\Big|\,
 \,\,
 UL_{k}[\bold i|k]\subseteq B(C,r)
\,\,\,\,
\text{for all $k\ge m$}
\Big\}\,.
\endalign
$$
\endroster
\endproclaim
\noindent{\it Proof}\newline
\noindent
(1)
For a function $f:\Sigma^{\Bbb N}\to\Bbb R$, let
$\Lip(f)$ denote the Lipschitz constant of $f$, i\.e\.
$\Lip(f)
=
\sup_{
\bold i,\bold j\in\Sigma^{\Bbb N},
\bold i\not=\bold j}
\frac
{|f(\bold i)-f(\bold j)|}
{\distance_{\Sigma^{\Bbb N}}(\bold i,\bold j)}$
and
define the metric $\LDistance$ in $\Cal P(\Sigma^{\Bbb N})$ by
$$
\LDistance(\mu,\nu)
=
\sup
\Sb
f:\Sigma^{\Bbb N}\to\Bbb R\\
\Lip(f)\le 1
\endSb
\Bigg|
\int f\,d\mu-\int f\,d\nu
\Bigg|;
$$
we note that it is well-known that $\LDistance$ is a metric
and  that $\LDistance$ induces the weak topology.
Since 
$U:\Cal P(\Sigma^{\Bbb N})\to X$ is continuous and 
$\Cal P(\Sigma^{\Bbb N})$ is compact, we
conclude that
$U:\Cal P(\Sigma^{\Bbb N})\to X$
is uniformly continuous.
This implies that we can choose $\delta>0$
such that
all measures
$\mu,\nu\in\Cal P(\Sigma^{\Bbb N})$
satisfy the following implication:
 $$
 \LDistance(\mu,\nu)\le\delta
 \,\,\,\,
 \Rightarrow
 \,\,\,\,
 \distance(U\mu,U\nu)\le \tfrac{r}{2}\,.
 \tag6.3
 $$

Next,
choose a positive integer $M_{r}$ such that
 $$
 \frac{1}{M_{r}(1-s_{\max})}
 <
 \delta\,;
 \tag6.4
 $$
recall, that $s_{\max}$ is defined in (4.1).

If $k\ge M_{r}$,
$\bold u\in\Sigma^{k}$ and $\bold k,\bold l\in\Sigma^{\Bbb N}$, then
it follows from 
(6.4) that
 $$
 \align
\LDistance
 \big(
 \,
 L_{k}(\bold u\bold k)
 \,,\,
 L_{k}(\bold u\bold l)
 \,
 \big)
&=
\sup
\Sb
f:\Sigma^{\Bbb N}\to\Bbb R\\
\Lip(f)\le 1
\endSb
\Bigg|
\int f\,d(L_{k}(\bold u\bold k))-\int f\,d(L_{k}(\bold u\bold l))
\Bigg|\\
&=
\sup
\Sb
f:\Sigma^{\Bbb N}\to\Bbb R\\
\Lip(f)\le 1
\endSb
\Bigg|
\frac{1}{k}\sum_{i=0}^{k-1} f(S^{i}(\bold u\bold k))
-
\frac{1}{k}\sum_{i=0}^{k-1} f(S^{i}(\bold u\bold l))
\Bigg|\\
&\le
\sup
\Sb
f:\Sigma^{\Bbb N}\to\Bbb R\\
\Lip(f)\le 1
\endSb
\frac{1}{k}\sum_{i=0}^{k-1} 
|
f(S^{i}(\bold u\bold k))
-
 f(S^{i}(\bold u\bold l))
|\\
&\le
\frac{1}{k}\sum_{i=0}^{k-1} 
\distance_{\Sigma^{\Bbb N}}
\big(
\,
S^{i}(\bold u\bold k)
\,,\,
S^{i}(\bold u\bold l)
\,
\big)\\
&=
\frac{1}{k}\sum_{i=0}^{k-1} 
s_{
S^{i}(\bold u\bold k)
\wedge
S^{i}(\bold u\bold l)
}\\
&\le
\frac{1}{M_{r}}\sum_{i=0}^{k-1} 
s_{\max}^{k-i}\\
&\le
 \frac{1}{M_{r}(1-s_{\max})}\\
&<
\delta\,,
\endalign
$$
and we therefore conclude from (6.3) that
$\distance
 (
 \,
 UL_{k}(\bold u\bold k)
 \,,\,
 UL_{k}(\bold u\bold l)
 \,
 )
 \le
 \tfrac{r}{2}$.

\noindent
(2)
It follows from (1) that
there is a positive integer $M_{r}$
such that
if $k\ge M_{r}$,
$\bold u\in\Sigma^{k}$ and $\bold k,\bold l\in\Sigma^{\Bbb N}$, then
 $\distance
 (
 \,
 UL_{k}(\bold u\bold k)
 \,,\,
 UL_{k}(\bold u\bold l)
 \,
 )
 \le
 \tfrac{r}{2}$.

We now claim that
if $m\ge M_{r}$,
then
 $$
 \align
 \Big\{
 \bold i\in\Sigma^{\Bbb N}
 \,\Big|\,
 \,\,
&UL_{k}\bold i\in B(C,\tfrac{r}{2})
\,\,\,\,
\text{for all $k\ge m$}
\Big\}\\
&\subseteq
\Big\{
 \bold i\in\Sigma^{\Bbb N}
 \,\Big|\,
 \,\,
 UL_{k}[\bold i|k]\subseteq B(C,r)
\,\,\,\,
\text{for all $k\ge m$}
\Big\}\,.
\endalign
$$
In order to prove this inclusion, we fix
$m\ge M_{r}$
and
$\bold i\in\Sigma^{\Bbb N}$
with
$UL_{k}\bold i\in B(C,\tfrac{r}{2})$
for all $k\ge m$.
We must now prove that
$UL_{k}[\bold i|k]\subseteq B(C,r)$
for all $k\ge m$.
We therefore fix $k\ge m$
and
$\bold j\in[\bold i|k]$.
We must now prove that
$UL_{k}\bold j\in B(C,r)$.
For brevity write
$\bold u=\bold i|k$.
Since $\bold j\in[\bold i|k]=[\bold u]$,
we can now find (unique) $\bold k,\bold l\in\Sigma^{\Bbb N}$ such that
$\bold i=\bold u\bold k$ and $\bold j=\bold u\bold l$.
We now have
 $$
 \align
 \dist
 \big(
 \,
 UL_{k}\bold j
 \,,\,
 C
 \,
 \big)
&\le
\distance
\big(
\,
UL_{k}\bold j
\,,\,
UL_{k}\bold i
\,
\big)
+
 \dist
 \big(
 \,
 UL_{k}\bold i
 \,,\,
 C
 \,
 \big)\\
 &=
\distance
\big(
\,
UL_{k}(\bold u\bold l)
\,,\,
UL_{k}(\bold u\bold k)
\,
\big)
+
 \dist
 \big(
 \,
 UL_{k}\bold i
 \,,\,
 C
 \,
 \big)\,.
 \tag6.5
\endalign
$$
However, since
$k\ge m\ge M_{r}$
and $\bold u\in\Sigma^{k}$,
we conclude that
 $\distance
 (
 \,
 UL_{k}(\bold u\bold k)
 \,,\,
 UL_{k}(\bold u\bold l)
 \,
 )
 \le
 \tfrac{r}{2}$.
Also,
since
$k\ge m$, we deduce that
$UL_{k}\bold i\in B(C,\tfrac{r}{2})$, whence
$\dist
 (
 \,
 UL_{k}\bold i
 \,,\,
 C
 \,
 )
\le
\tfrac{r}{2}$.
It therefore follows from (6.5) that
 $$
 \align
 \dist
 \big(
 \,
 UL_{k}\bold j
 \,,\,
 C
 \,
 \big)
 &=
\distance
\big(
\,
UL_{k}(\bold u\bold l)
\,,\,
UL_{k}(\bold u\bold k)
\,
\big)
+
 \dist
 \big(
 \,
 UL_{k}\bold i
 \,,\,
 C
 \,
 \big)\\
&\le
 \frac{r}{2}+\frac{r}{2}\\
&=
 r\,.  
\endalign
$$
This completes the proof.
\hfill$\square$

\bigskip

\proclaim{Lemma 6.2}
Let $X$ be a metric space
and 
let $U:\Cal P(\Sigma^{\Bbb N})\to X$ be continuous with respect to the weak topology.
Let $C\subseteq X$ be a closed subset of $X$.
Then
 $$
 \dim_{\Haus}
 \Big\{
 \bold i\in\Sigma^{\Bbb N}
 \,\Big|\,
 \lim_{m}\dist(\,UL_{m}\bold i,C\,)=0
 \Big\}
 \le
 \underline f^{U,\Lambda}(C)\,;
 $$
recall that $\dim_{\Haus}$ denotes the Hausdorff dimension.
\endproclaim
\noindent{\it Proof}\newline
\noindent
For a subset $\Xi$ of $\Sigma^{\Bbb N}$, we let
$\underline\dim_{\Box}\Xi$ denote the lower box dimension of $\Xi$; the 
reader is referred to [Fa1] for the definition of the lower box 
dimension.
We will use the fact that $\dim_{\Haus} \Xi\le\underline\dim_{\Box}\Xi$
for all $\Xi\subseteq\Sigma^{\Bbb N}$, see, for example, [Ed].

We now introduce the following notation.
For brevity write
 $$
\align
\Gamma
&=
 \Big\{
 \bold i\in\Sigma^{\Bbb N}
 \,\Big|\,
 \lim_{m}\dist(\,UL_{m}\bold i,C\,)=0
 \Big\}\,.\\
 \endalign
 $$
Also, for a positive integer $m$ and a positive real number $r>0$, write
 $$
 \align
\Gamma_{m}(r)
&=
 \Big\{
 \bold i\in\Sigma^{\Bbb N}
 \,\Big|\,
 \text{
 $UL_{k}\bold i\in B(C, r)$
 for all $k\ge m$
 }
 \Big\}\,,\\
\Delta_{m}(r)
&=
 \Big\{
 \bold i\in\Sigma^{\Bbb N}
 \,\Big|\,
 \text{
 $UL_{k}[\bold i|k]\subseteq B(C, r)$
 for all $k\ge m$
 }
 \Big\}\,.
 \endalign
 $$
Observe that if $M$ is any positive integer, then  
we clearly have
 $$
 \Gamma
 \subseteq
 \bigcup_{m\ge M}
 \Gamma_{m}(\tfrac{r}{2})
 \tag6.6
 $$
for all $r>0$.
We also observe that it follows from Lemma 6.1 that for each positive number $r>0$
there is a positive integer $M_{r}$
such that
 $$
 \Gamma_{m}(\tfrac{r}{2})
 \subseteq
 \Delta_{m}(r)
 \tag6.7
 $$
 for all $m\ge M_{r}$.
 It follows from (6.6) and (6.7) that
  $$
  \align
 \Gamma
 &\subseteq
 \bigcup_{m\ge M_{r}}
 \Gamma_{m}(\tfrac{r}{2})\\
&\subseteq
 \bigcup_{m\ge M_{r}}
  \Delta_{m}(r)\,,
\endalign
 $$
whence 
 $$
 \align
 \dim_{\Haus} \Gamma
&\le
 \dim_{\Haus}\Bigg(\bigcup_{m\ge M_{r}}\Delta_{m}(r)\Bigg)\\
&=
\sup_{m\ge M_{r}}\dim_{\Haus} \Delta_{m}(r)\\
&\le
 \sup_{m\ge M_{r}}\underline{\dim}_{\Box}\Delta_{m}(r)
 \tag6.8
 \endalign
 $$
for all $r>0$.

Fix a positive integer $m$.
We now prove that
 $$
 \Delta_{m}(r)
 \subseteq
 \bigcup_{\bold i\in\Pi_{\delta}^{U,\Lambda}(C,r)}[\bold i]
 \tag6.9
 $$
for all $0<\delta<s_{\min}^{m}$ and all $r>0$.
Indeed, fix $\bold j\in \Delta_{m}(r)$.
Now, let $k_{0}$ denote the unique positive integer
such that
if we write $\bold j_{0}=\bold j|k_{0}$, then
$s_{\bold j_{0}}\le \delta<s_{\widehat{\bold j_{0}}}$, i\.e\.
$s_{\bold j_{0}}\approx \delta$.
Since
it follows from Proposition 4.1 that
 $s_{\min}^{k_{0}}
 =
 s_{\min}^{|\bold j_{0}|}
 \le 
 s_{\bold j_{0}}
 \le 
 \delta
 <
  s_{\min}^{m}$,
we conclude that $k_{0}\ge m$, 
and the fact that $\bold j\in \Delta_{m}(r)$ therefore implies that
$UL_{|\bold j_{0}|}[\bold j_{0}]
=
UL_{k_{0}}[\bold j|k_{0}]
\subseteq
B(C,r)$
This shows that $\bold j_{0}\in\Pi_{\delta}^{U,\Lambda}(C,r)$, whence
$\bold j\in[\bold j|k_{0}]
=
[\bold j_{0}]
\subseteq
\cup_{\bold i\in\Pi_{\delta}^{U,\Lambda}(C,r)}[\bold i]$.
This proves (6.9).

Inclusion (6.9) shows that
for all $0<\delta<s_{\min}^{m}$, the family 
$(\,[\bold i]\,)_{\bold i\in\Pi_{\delta}^{U,\Lambda}(C,r)}$
is a covering of $\Delta_{m}(r)$ of sets 
$[\bold i]$
with
$\bold i\in\Pi_{\delta}^{U,\Lambda}(C,r)$
such that
$\diam [\bold i]=s_{\bold i}\le\delta$
for all $\bold i\in\Pi_{\delta}^{U,\Lambda}(C,r)$.
This implies that
 $$
 \underline\dim_{\Box}\Delta_{m}(r)
 \le
 \liminf_{\delta\searrow 0}
 \frac
 {\log|\Pi_{\delta}^{U,\Lambda}(C,r)|}
 {-\log \delta}
 \tag6.10
 $$
for all $r>0$.
Since (6.10) holds for all $m$, we conclude that
 $$
 \sup_{m\ge M_{r}}\underline\dim_{\Box}\Delta_{n}(r)
 \le
 \liminf_{\delta\searrow 0}
 \frac
 {\log|\Pi_{\delta}^{U,\Lambda}(C,r)|}
 {-\log \delta}
 \tag6.11
 $$
for all $r>0$.
 
Combining (6.8) and (6.11) now shows that
 $$
 \dim_{\Haus}\Gamma
 \le
 \liminf_{\delta\searrow 0}
 \frac
 {\log|\Pi_{\delta}^{U,\Lambda}(C,r)|}
 {-\log \delta}
 \tag6.12
 $$
for all $r>0$. Finally, letting $r\searrow 0$ in (6.12) completes the proof.
\hfill$\square$

\bigskip

In order to statement and prove the 
next lemma we introduce the following notation.
Namely,
for a H\"older continuous function
$\varphi:\Sigma^{\Bbb N}\to\Bbb R$, we will write
 $$
 P(\varphi)
 $$
for the topological pressure of $\varphi$.
We can now state and prove Lemma 6.3.

\bigskip

\proclaim{Lemma 6.3}
Let $\mu\in\Cal P_{S}(\Sigma^{\Bbb N})$ with $\supp\mu=\Sigma^{\Bbb N}$.
(Here $\supp\mu$ denotes the topological support of $\mu$.)
Then there exists a sequence $(\mu_{n})_{n}$ of probability measures 
on $\Sigma^{\Bbb N}$ satisfying the following three conditions.
\roster
\item"(1)" We have $\mu_{n}\to\mu$ weakly.
\item"(2)" For each $n$, the measure $\mu_{n}$ is ergodic.
\item"(3)" We have $h(\mu_{n})\to h(\mu)$.
\endroster
\endproclaim
\noindent{\it Proof}\newline
\noindent 
Fix a positive integer $n$.
Since $\supp\mu=\Sigma^{\Bbb N}$, we deduce that
$\mu[\bold i]>0$ for all
$\bold i\in\Sigma^{*}$.
Hence, for $m\in\Bbb N$ and $ i_{1}\ldots i_{m}\in\Sigma^{m}$,
we can define $p_{n, i_{1}\ldots i_{m}}$ by
 $$
 p_{n, i_{1}\ldots i_{m}}
 =
 \cases
 \mu[ i_{1}\ldots i_{m}]
&\qquad
 \text{for $m\le n$,}\\
&{}\\
{\dsize
 \prod_{k=1}^{m-n}
 \,
 \frac
 {\mu[ i_{k} i_{k+1}\ldots i_{k+(n-1)}]}
 {\mu[ i_{k+1}\ldots i_{k+(n-1)}]}
 \,\,\,
 \mu[ i_{(m-n)+1}\ldots i_{m}]
 }
&\qquad
 \text{for $n<m$.}
 \endcases
 \tag6.13
 $$
Since clearly
$\sum_{i}p_{n,i}=1$ and
$\sum_{i}p_{ n,i_{1}\ldots i_{m} i}=
p_{n, i_{1}\ldots i_{m}}$ for all
$m$
and
all
$ i_{1}\ldots i_{m}\in\Sigma^{m}$, 
there exists a (unique) probability measure 
$\mu_{n}$ on $\Sigma^{\Bbb N}$ such that
 $$
 \mu_{n}[ i_{1}\ldots i_{m}]
 =
 p_{ n,i_{1}\ldots i_{m}}
 $$
for all $m$ and all $ i_{1}\ldots i_{m}\in\Sigma^{m}$
(cf\. [Wa, p\. 5]).

\bigskip

{\it Claim 1. We have $\mu_{n}\to\mu$ weakly.}

\noindent{\it Proof of Claim 1.}
It follows from definition (6.13) that
$\mu_{n}[\bold i]=\mu[\bold i]$ for all $\bold i\in\Sigma^{n}$.
This clearly implies that 
$\mu_{n}\to\mu$ weakly.
This completes the proof of Claim 1.

\bigskip

{\it Claim 2. For each $n$, 
there is 
a 
H\"older continuous 
function
$\varphi_{n}:\Sigma^{\Bbb N}\to\Bbb R$
such that
the following conditions hold.

(1) $P(\varphi_{n})=0$\,,

(2) The measure $\mu_{n}$ is a Gibbs state of
$\varphi_{n}$.
}

\noindent{\it Proof of Claim 2.}
We first note that $\mu_{n}$ is shift invariant. Indeed, 
since $\mu$ is shift invariant, a small calculation shows that
$\sum_{i}\mu_{n}[ i \bold i]
=
\mu_{n}[ \bold i]$ for all
$\bold i\in\Sigma^{*}$.
This implies that
$\mu_{n}(S^{-1}[\bold i])=\mu_{n}[\bold i]$ for all 
$\bold i\in\Sigma^{*}$, whence
$\mu_{n}(S^{-1}B)=\mu_{n}(B)$
for all Borel sets $B$.

Next we show that $\mu_{n}$ is a Gibbs state for a H\"older continuous 
function.
Define $\varphi_{n}:\Sigma^{\Bbb N}\to\Bbb R$ by
 $$
 \varphi_{n}( i_{1} i_{2}\ldots)
 =
 \log
 \left(
 \frac
 {\mu[ i_{1} i_{2}\ldots i_{n}]}{\mu[ i_{2}\ldots i_{n}]}
 \right)\,.
 $$
The map $\varphi_{n}$ is clearly H\"older continuous, and it follows 
from the definition of $\mu_{n}$ that 
 $$
 e^{-n\|\varphi_{n}\|_{\infty}}\min_{\bold j\in{\Sigma^{n}}}\mu[\bold j]
 \le
 \frac{\mu_{n}[\bold i|m]}
 {e^{\sum_{k=0}^{m-1}\varphi_{n}(S^{k}\bold i)}}
 \le
 e^{n\|\varphi_{n}\|_{\infty}}\max_{\bold j\in{\Sigma^{n}}}\mu[\bold j]
 $$
for all $\bold i\in\Sigma^{\Bbb N}$ and all $m>n$.
This shows that $\mu_{n}$ is the Gibbs state of $\varphi_{n}$, and 
that the pressure
$P(\varphi_{n})$ of $\varphi_{n}$ equals $0$, i\.e\.
$P(\varphi_{n})=0$; cf\. [Bo].
This completes the proof of Claim 2.

\bigskip

{\it Claim 3. For each $n$, 
the measure $\mu_{n}$ is ergodic.
}

\noindent{\it Proof of Claim 3.}
It follows from  Claim 2 that
$\mu_{n}$ is the 
a Gibbs state of a 
H\"older continuous function. This implies  that  $\mu_{n}$ is ergodic.
This completes the proof of Claim 3.

\bigskip

{\it Claim 4. We have $h(\mu_{n})\to h(\mu)$.}

\noindent{\it Proof of Claim 4.}
For measurable partitions $\Cal A,\Cal B$ of $\Sigma$, let
$h(\mu;\Cal A)$ and $h(\mu;\Cal A|\Cal B)$
denote the entropy of $\Cal A$ with respect to $\mu$,
and the conditional entropy of
$\Cal A$ given $\Cal B$ with respect to $\mu$, respectively.
Write
$\Cal C=\{\,[ i]\,|\, i\in\Sigma\}$ and
$\Cal C_{n}
=
\vee_{k=0}^{n-1}S^{-k}\Cal C
=
\{\,[\bold i]\,|\, \bold i\in\Sigma^{n}\}$.
It follows from Claim 2 that there is a H\"older continuous function 
$\varphi_{n}:\Sigma^{\Bbb N}\to\Bbb R$
with $P(\varphi_{n})=0$ such that
$\mu_{n}$ 
is a Gibbs state of $\varphi_{n}$.
Since
$P(\varphi_{n})=0$ and $\mu_{n}$ is a Gibbs state of $\varphi_{n}$,
the Variational Principle now shows that
$0=P(\varphi_{n})=h(\mu_{n})+\int\varphi_{n}\,d\mu_{n}$ (cf\. [Bo]), whence
 $$
 \align
 h(\mu_{n})
 &=
 -\int\varphi_{n}\,d\mu_{n}\\
&=
 -\sum_{ i_{1}\ldots i_{n}}
 \mu[ i_{1}\ldots i_{n}]
 \log
 \left(
 \frac
 {\mu[ i_{1} i_{2}\ldots i_{n}]}{\mu[ i_{2}\ldots i_{n}]}
 \right)\\
&=
 \,
 h(\mu;\Cal C_{n}|\Cal C_{n-1})\,.
 \tag6.14
 \endalign
 $$
Next, we note that it follows from  [DGS, 11.4]
that
$h(\mu;\Cal C_{n}|\Cal C_{n-1})\to h(\mu;\Cal C)$,
and we therefore conclude from (6.14) that
 $$
 h(\mu_{n})
 \to
 h(\mu;\Cal C)\,.
 \tag6.15
 $$
Finally, it follows immediately from 
the Kolmogoroff-Sinai theorem
that
 $h(\mu;\Cal C)=h(\mu)$.
 This and (6.15) now show that
 $$
 h(\mu_{n})
 \to
 h(\mu)\,.
 $$ 
This completes the proof of Claim 4.

\bigskip

The proof now follows from Claim 1, Claim 3 and Claim 4.
\hfill$\square$

\bigskip

The next auxiliary result
provides a formula for the
upper Hausdorff 
dimension of is a 
probability measure.
If $\mu$ is a 
probability measure on $\Sigma^{\Bbb N}$, we define the upper Hausdorff 
dimension of $\mu$ by
 $$
 \overline\dim_{\Haus}\mu
 =
 \inf
 \Sb
 \Xi\subseteq\Sigma^{\Bbb N}\\
 \mu(\Xi)=1
 \endSb
 \dim_{\Haus} \Xi \,.
 $$
(Recall that $\dim_{\Haus}$ denotes the Hausdorff dimension.) 
The next result provides a formula for the upper Hausdorff 
dimension of an ergodic probability measure on $\Sigma^{\Bbb N}$.
This result is folklore
 and follows  from
the
Shannon-MacMillan-Breiman theorem and the ergodic theorem. 
However, for sake of completeness we have decided to include the 
short proof.

\bigskip

\proclaim{Proposition 6.4}
Let $\mu$ be an ergodic probability measure on $\Sigma^{\Bbb N}$. Then
$\overline\dim_{\Haus}\mu
 =
 -
 \frac{h(\mu)}{\int\Lambda\,d\mu}$.
\endproclaim
\noindent{\it Proof}\newline
\noindent
Since $\mu$ is ergodic, it  follows from
the
Shannon-MacMillan-Breiman theorem that
 $$
 \frac{\log\mu([\bold i|n])}{n}
 \to
 -
h(\mu)
 \quad
 \text{for $\mu$-a\.a\. $\bold i\in\Sigma^{\Bbb N}$.}
 \tag6.16
 $$
 Also, an application of the ergodic theorem shows that
 $\frac{\sum_{k=0}^{n-1}\Lambda S^{k}\bold i}{n}
 \to
\int\Lambda\,d\mu$
for $\mu$-a\.a\. $\bold i\in\Sigma^{\Bbb N}$.
It follows from this and Proposition 4.1 that
$$
 \frac{\log s_{\bold i|n}}{n}
 \to
\int\Lambda\,d\mu
 \quad
 \text{for $\mu$-a\.a\. $\bold i\in\Sigma^{\Bbb N}$.}
 \tag6.17
 $$ 
Combining (6.16) and (6.17) now gives
$$
 \frac{\log\mu([\bold i|n])}{\log s_{\bold i|n}}
 \to
  -
\frac{h(\mu)}{\int\Lambda\,d\mu}
 \quad
 \text{for $\mu$-a\.a\. $\bold i\in\Sigma^{\Bbb N}$.}
 \tag6.18
 $$

Next,
for each $\bold i\in\Sigma^{\Bbb N}$ and $r>0$, let 
$n_{\bold i,r}$ denote the unique integer such that
$s_{\bold i|n_{\bold i,r}}<r\le s_{\widehat {\bold i|n_{\bold i,r}}  } $.
It follows from the definition of the metric $\distance _{\Sigma^{\Bbb N}}$
on $\Sigma^{\Bbb N}$ (see (6.1) and (6.2))
that
$B(\bold i,r)=[\bold i|n_{\bold i,r}]$.
Also, if we let $c$ denote the constant from Condition (C3) in Section 2.1, then
it follows from Proposition 4.1 that
$s_{\bold i|n_{\bold i,r}}
<
r
\le 
s_{\widehat {\bold i|n_{\bold i,r}}  } 
\le
\frac{c}{s_{\min}}s_{\bold i|n_{\bold i,r}}$.
Combining these facts, we now deduce from (6.18) that
 $$
 \align
\lim_{r\searrow 0}
 \frac{\log\mu(B(\bold i,r))}{\log r}
&=
\lim_{r\searrow 0}
 \frac{\log\mu([\bold i|n_{\bold i,r}])}{\log s_{\bold i|n_{\bold i,r}}}\\
&=
\lim_{n}
 \frac{\log\mu([\bold i|n])}{\log s_{\bold i|n}}\\ 
&=
  -
\frac{h(\mu)}{\int\Lambda\,d\mu}
 \quad
 \text{for $\mu$-a\.a\. $\bold i\in\Sigma^{\Bbb N}$,}
 \endalign
 $$
 whence
 $$
 \mu\text{-}\!\esssup_{\bold i}\liminf_{r\searrow 0}\frac{\log\mu(B(\bold i,r))}{\log r}\
 =
   -
\frac{h(\mu)}{\int\Lambda\,d\mu}\,,
\tag6.19
$$
where
$\mu\text{-}\!\esssup$
denotes the 
$\mu$ essential supremum.

Finally, we note that it is well-known that
$\overline\dim_{\Haus}\mu
=
\mu\text{-}\!\esssup_{\bold i}\liminf_{r\searrow 0}\frac{\log\mu(B(\bold i,r))}{\log r}$
(see, for example, [Fa2]),
and it therefore follows immediately from (6.19) that
$\overline\dim_{\Haus}\mu
=
\mu\text{-}\!\esssup_{\bold i}\liminf_{r\searrow 0}\frac{\log\mu(B(\bold i,r))}{\log r}
=
  -
\frac{h(\mu)}{\int\Lambda\,d\mu}$.
\hfill$\square$

\bigskip

 The final auxiliary result says that the map 
 $C\to\underline f^{U,\Lambda}(C)$ is upper 
 semi-continuous.
 In order to state this result we introduce the following
 notation. 
For a 
metric space $X$,
we write
 $$
 \Cal F(X)
 =
 \Big\{
 F\subseteq X\,\Big|\,\text{$F$ is closed and non-empty}
 \Big\}
 \tag6.20
 $$
and we
equip $\Cal F(X)$ with the Hausdorff metric $\Distance$;
recall, that since $X$ may be unbounded, the 
Hausdorff distance $\Distance$ is defined as follows, namely, for 
$E,F\in\Cal F(X)$, write
$$
\Delta(E,F)
=
\min
\Bigg(
\,
\sup_{x\in E}\dist(x,F)
\,,\,
\sup_{y\in F}\dist(y,E) 
\,
\Bigg)
\tag6.21
$$
and define $\Distance$ by
 $$
 \Distance
 =
 \min(1,\Delta)\,.
 \tag6.22
 $$
\bigskip

\proclaim{Lemma 6.5}
Let $X$ be a metric space
and let $U:\Cal P(\Sigma^{\Bbb N})\to X$ be 
continuous with respect to the weak topology.
Equip $\Cal F(X)$ with the Hausdorff metric $\Distance$.
Then the function
$\underline f^{U,\Lambda}:\Cal F(X)\to\Bbb R$
is upper semicontinuous, i\.e\.
for each $C\in\Cal F(X)$ and each $\varepsilon>0$, there exists
a real number $\rho>0$ 
such that
if
$F\in\Cal F(X)$ and $\Distance(F,C)<\rho$, then
 $$
 \underline f^{U,\Lambda}(F)
 \le
 \underline f^{U,\Lambda}(C)+\varepsilon\,.
 $$
\endproclaim
\noindent{\it Proof}\newline
\noindent
Let $C\in\Cal F(X)$ and $\varepsilon>0$. Next, 
it follows from the definition of 
$\underline f^{U,\Lambda}(C)$ that we can choose a real number $r_{0}$
with $0<r_{0}<1$ such 
that
 $$
 \align
  \underline f^{U,\Lambda}(C,r_{0})
&<
 \underline f^{U,\Lambda}(C)+\varepsilon\,.
 \tag6.23
 \endalign
 $$
Let $\rho=\frac{r_{0}}{2}$. 
We now prove the following claim.

\medskip

{\it Claim 1. 
Let $F\in\Cal F(X)$ with $\Distance(F,C)<\rho$.
For all $0<r<\rho$ and all $\delta>0$, we have
 $$
  N_{\delta}^{U,\Lambda}(F,r)
 \le
  N_{\delta}^{U,\Lambda}(C,r_{0})\,.
 $$
 }
\noindent{\it Proof of Claim 1.}
Fix $0<r<\rho$ and $\delta>0$.
Since $\Distance(F,C)<\rho=\frac{r_{0}}{2}$ and $r_{0}<1$,
we first conclude that
$B(F,\frac{r_{0}}{2})\subseteq B(C,r_{0})$.
Hence, if $\bold i\in\Pi_{\delta}^{U,\Lambda}(F,r)$, then
this and the fact that
$0<r<\rho=\frac{r_{0}}{2}$ imply that
$UL_{|\bold i|}[\bold i]
\subseteq
B(F,r)
\subseteq
B(F,\rho)
=
B(F,\frac{r_{0}}{2})
\subseteq B(C,r_{0})$
and so
$\bold i\in \Pi_{\delta}^{U,\Lambda}(C,r_{0})$.
This shows that
$\Pi_{\delta}^{U,\Lambda}(F,r)
 \subseteq
 \Pi_{\delta}^{U,\Lambda}(C,r_{0})$, whence
$ N_{\delta}^{U,\Lambda}(F,r)
 \le
  N_{\delta}^{U,\Lambda}(C,r_{0})$. 
This completes the proof of Claim 1.

\medskip

We now claim that 
if
$F\in\Cal F(X)$ and $\Distance(F,C)<\rho$, then
 $$
 \underline f^{U,\Lambda}(F)
 \le
 \underline f^{U,\Lambda}(C)+\varepsilon\,.
 \tag6.24
 $$
To prove this, let $F\in\Cal F(X)$ with $\Distance(F,C)<\rho$. 
It follows from Claim 1 and (6.23) that if $0<r<\rho$, then
 $$
 \align
 \underline f^{U,\Lambda}(F,r)
&=
 \liminf_{\delta\searrow 0}
 \frac{\log  N_{\delta}^{U,\Lambda}(F,r)}{-\log \delta}\\
&\le
 \liminf_{\delta\searrow 0}
 \frac{\log  N_{\delta}^{U,\Lambda}(C,r_{0})}{-\log \delta}\\
&=
  \underline f^{U,\Lambda}(C,r_{0})\\
&<
 \underline f^{U,\Lambda}(C)+\varepsilon\,.
 \endalign
 $$
Since this inequality holds for all $0<r<\rho$, we finally conclude that 
$ \underline f^{U,\Lambda}(F)
=
\lim_{r\searrow 0}\underline f^{U,\Lambda}(F,r)
 \le
 \underline f^{U,\Lambda}(C)+\varepsilon$.
\hfill$\square$

\bigskip

We can now state and prove the main result in this section, namely,
Theorem 6.6
providing a proof of inequality (2.2).

\bigskip

\proclaim{Theorem 6.6}
Let $X$ be a metric space
and let $U:\Cal P(\Sigma^{\Bbb N})\to X$ be 
continuous with respect to the weak topology.
Let $C\subseteq X$ be a closed subset of $X$. We have
 $$ 
  \sup
   \Sb
   \mu\in\Cal P_{S}(\Sigma^{\Bbb N})\\
   {}\\
   U\mu\in C
   \endSb
 -
 \frac{h(\mu)}{\int\Lambda\,d\mu}
 \le
 \underline f^{U,\Lambda}(C)
 $$
\endproclaim 
\noindent{\it Proof}\newline
\noindent
Let $\varepsilon>0$.
Next,
fix $\mu\in\Cal P_{S}(\Sigma^{\Bbb N})$ with $U\mu\in C$.
We will now prove that
 $$ 
  -
 \frac{h(\mu)}{\int\Lambda\,d\mu}
 \le
 \underline f^{U,\Lambda}(C)
 +
 \varepsilon\,.
 \tag6.25
 $$

Let $\Cal F(X)$ be denied as in (6.20), i\.e\.
$\Cal F(X)
=
\{F\subseteq X\,|\,\text{$F$ is closed and non-empty}
\}$,
and 
and equip $\Cal F(X)$ with the Hausdorff metric $\Distance$, see (6.21) and (6.22).
It follows from Lemma 6.5 that
the function
$ \underline f^{U,\Lambda}:\Cal F(X)\to\Bbb R$ is upper 
semi-continuous, 
and we can therefore choose 
$\rho_{\varepsilon}>0$ such that:
 $$
 \gathered
 \text{
 if $F\in\Cal F(X)$ and $\Distance(F,C)<\rho_{\varepsilon}$,
then}\\
 \underline f^{U,\Lambda}(F)
 \le
 \underline f^{U,\Lambda}(C)+\varepsilon\,.
 \endgathered
 \tag6.26
 $$

Next, observe that we can choose
an $S$-invariant probability measure $\gamma$ on $\Sigma^{\Bbb N}$ such that 
$\supp\gamma=\Sigma^{\Bbb N}$.
For $t\in(0,1)$, we now write
$\mu_{t}=(1-t)\mu+t\gamma\in\Cal P_{S}(\Sigma^{\Bbb N})$.
As $U$ is continuous with $U\mu\in C$ and $\mu_{t}\to\mu$ weakly as $t\searrow 0$,
there exists 
$0<t_{\varepsilon}<1$ such that
for all $0<t<t_{\varepsilon}$, we have
 $$
 \dist(U\mu_{t},C)
 <
 \rho_{\varepsilon}\,.
 \tag6.27
 $$

Fix $0<t<t_{\varepsilon}$.
Since $U$ is continuous and  
$\dist(U\mu_{t},C)
 <
 \rho_{\varepsilon}$ (by (6.27)),
it follows from Lemma 6.3  that we may choose 
a sequence $(\mu_{t,n})_{n}$ of  $S$-invariant probability 
measures on $\Sigma^{\Bbb N}$ such that
 $$
 \align
& \mu_{t,n}\to\mu_{t}\,\,\,\,\text{weakly}\,,
\tag6.28\\
&\text{$\mu_{t,n}$ is ergodic,}
\tag6.29\\
&h(\mu_{t,n})\to h(\mu_{t})
\tag6.30
\endalign
 $$
and
 $$
 \dist(U\mu_{t,n},C)<\rho_{\varepsilon}
 \tag6.31
 $$
for all $n$. 
Observe that it follows from  
(6.31) that
$\Distance(\,C\cup\{\mu_{t,n}\}  \,,\, C  \,)<\rho_{\varepsilon}$,
and we therefore conclude from (6.31) that
 $$
 \underline f^{U,\Lambda}(\,C\cup\{\mu_{t,n}\}  \,)
 \le
 \underline f^{U,\Lambda}(C)
 +
 \varepsilon
 \tag6.32
 $$
for all $n$.
We now prove the following two claims.

\bigskip

{\it Claim 1. For all $0<t<t_{\varepsilon}$, we have
 $$
 \align
 -
 \frac
 {(1-t)h(\mu)+th(\gamma)}
 {(1-t)\int\Lambda\,d\mu+t\int\Lambda\,d\gamma}
&\le
 \lim_{n} \overline\dim_{\Haus}\,\mu_{t,n}\,.
 \endalign
 $$
}
\noindent
{\it Proof of Claim 1.}
Using the fact that the entropy map
$h:\Cal P_{S}(\Sigma)\to\Bbb R$ is affine (cf\. [Wa]) 
we conclude that
 $$
 \align
 -
 \frac
 {(1-t)h(\mu)+th(\gamma)}
 {(1-t)\int\Lambda\,d\mu+t\int\Lambda\,d\gamma}
&\le
 -
 \frac
 {h((1-t)\mu+t\gamma)}
 {\int\Lambda\,d((1-t)\mu+t\gamma)}\\
&=
 -\frac{h(\mu_{t})}{\int\Lambda\,d\mu_{t}}\,.
%&=
% \lim_{n}-\frac{h(\mu_{t,n})}{\int\Lambda\,d\mu_{t,n}}\,.
\tag6.33
 \endalign
 $$
However, since $\Lambda$ is continuous and 
$\mu_{t,n }\to\mu_{t}$ weakly (by (6.28)),
we conclude that
$\int\Lambda\,d\mu_{t,n}\to\int\Lambda\,d\mu_{t}$.
We deduce from this and the fact that
$h(\mu_{t,n})\to h(\mu_{t})$ (by (6.30)) that
 $$
 -\frac{h(\mu_{t})}{\int\Lambda\,d\mu_{t}}
=
 \lim_{n}-\frac{h(\mu_{t,n})}{\int\Lambda\,d\mu_{t,n}}\,.
 \tag6.34
$$
Combining (6.33) and (6.34) now yields
$$
 -
 \frac
 {(1-t)h(\mu)+th(\gamma)}
 {(1-t)\int\Lambda\,d\mu+t\int\Lambda\,d\gamma}
 \le
 \lim_{n}-\frac{h(\mu_{t,n})}{\int\Lambda\,d\mu_{t,n}}\,.
 \tag6.35
 $$

Also, since $\mu_{t,n}$ is ergodic (by (6.29)), it follows from Proposition 6.4 that
$\overline\dim_{\Haus}\,\mu_{t,n}
 =
 -\frac{h(\mu_{t,n})}{\log N}$,
 and we therefore conclude from (6.35) that
 $$
 \align
 -
 \frac
 {(1-t)h(\mu)+th(\gamma)}
 {(1-t)\int\Lambda\,d\mu+t\int\Lambda\,d\gamma}
&\le
 \lim_{n}-\frac{h(\mu_{t,n})}{\int\Lambda\,d\mu_{t,n}}\\
&=
 \lim_{n} \overline\dim_{\Haus}\,\mu_{t,n}\,.
 \endalign
 $$
This completes the proof of Claim 1.

\bigskip

{\it Claim 2. For all $0<t<t_{\varepsilon}$, we have
 $$
 \lim_{n}\overline\dim_{\Haus}\,\mu_{t,n}
 \le
  \underline f^{U,\Lambda}(C)
 +
 \varepsilon\,.
 $$
}
\noindent
{\it Proof of Claim 2.}
It follows immediately 
from the ergodicity of $\mu_{t,n}$ and the ergodic theorem that
$\mu_{t,n}( \{\bold i\in\Sigma^{\Bbb N}\mid\lim_{m}L_{m}\bold i=\mu_{t,n}\} )=1$.
Hence
 $$
 \align
 \overline\dim_{\Haus}\,\mu_{t,n}
&\le
 \dim_{\Haus} 
 \Big\{
 \bold i\in\Sigma^{\Bbb N}
 \,\Big|\,\lim_{m}L_{m}\bold i=\mu_{t,n}
 \Big\}\\
&\le
 \dim_{\Haus}
 \Big\{
 \bold i\in\Sigma^{\Bbb N}
 \,\Big|\,
 \lim_{m}
 UL_{m}\bold i
  =
  U\mu_{t,n}
 \Big\}\\ 
&\le
 \dim_{\Haus}
 \Big\{
 \bold i\in\Sigma^{\Bbb N}
 \,\Big|\,
 \lim_{m}
 \dist\big(
 \,
 UL_{m}\bold i\,,\,
 C\cup\{U\mu_{t,n}\}
 \,
 \big)
 =0
 \Big\}\,.
 \tag6.36
 \endalign
 $$
 Next, it follows from (6.36) using Lemma 6.2 and (6.32) that
 $$
  \align
 \overline\dim_{\Haus}\,\mu_{t,n}
&\le
 \dim_{\Haus}
 \Big\{
 \bold i\in\Sigma^{\Bbb N}
 \,\Big|\,
 \lim_{m}
 \dist\big(
 \,
 UL_{m}\bold i\,,\,
 C\cup\{U\mu_{t,n}\}
 \,
 \big)
 =0
 \Big\}
 \qquad
 \text{[by (6.36)]}\\
&\le 
 \underline f^{U,\Lambda}(\,C\cup\{\mu_{t,n}\}  \,)
 \qquad\qquad
 \qquad\qquad
 \qquad\qquad
 \qquad\qquad
 \,\,\,\,\,
 \text{[by Lemma 6.2]}\\
&\le
 \underline f^{U,\Lambda}(C)
 +
 \varepsilon\,.
  \qquad\qquad
 \qquad\qquad
 \qquad\qquad
 \qquad\qquad
 \qquad\,\,\,\,\,
 \text{[by (6.32)]}
 \endalign
 $$
This completes the proof of Claim 2.

\bigskip

Combining
Claim 1 and Claim 2 shows that for all $0<t<t_{\varepsilon}$, we have
 $$
 \align
 -
 \frac
 {(1-t)h(\mu)+th(\gamma)}
 {(1-t)\int\Lambda\,d\mu+t\int\Lambda\,d\gamma}
&\le
 \underline f^{U,\Lambda}(C)
 +
 \varepsilon\,.
 \tag6.37
 \endalign
 $$ 
Letting $t\searrow 0$ in (6.37) now gives
gives
$-
 \frac
 {h(\mu)}
 {\int\Lambda\,d\mu}
\le
 \underline f^{U,\Lambda}(C)
 +
 \varepsilon$.
This proves (6.25).

Since $\mu\in\Cal P_{S}(X)$ with $U\mu\in C$ was arbitrary, 
it follows immediately from 
(6.25) that
 $$
 \sup
 \Sb
 \mu\in\Cal P(\Sigma^{\Bbb N})\\
 {}\\
 U\mu\in C
 \endSb
 -
 \frac
 {h(\mu)}
 {\int\Lambda\,d\mu}
\le
 \underline f^{U,\Lambda}(C)
 +
 \varepsilon\,.
 $$
Finally, letting $\varepsilon\searrow 0$ gives the desired result.
\hfill$\square$

  \bigskip
  \bigskip

 %%%%%%%%%%%%%%%%%%%%%%%%%%%%%%%%%%%%%%%%%
 %%%%%%%%%%%%%%%%%%%%%%%%%%%%%%%%%%%%%%%%%
 %%%%%%%%%%%%%%%%%%%%%%%%%%%%%%%%%%%%%%%%%
 %%%%%%%%%%%%%%%%%%%%%%%%%%%%%%%%%%%%%%%%%
 %%%%%%%%%%%%%%%%%%%%%%%%%%%%%%%%%%%%%%%%% 

\heading
{
7. Proof of inequality (2.3)
% $$
% \underline f^{U,\Lambda}(C)
% \le
%  \liminf_{r\searrow 0}
% \,\,
% \sigma_{\abs}\big(\,\zeta_{C}^{U,\Lambda}(\cdot;r)\,\big) 
% $$
}
\endheading

The purpose of this section is to prove
Theorem 7.1 providing a proof of inequality (2.3).

\bigskip

\proclaim{Theorem 7.1}
Let $X$ be a metric space and let $U:\Cal P(\Sigma^{\Bbb N})\to X$ be 
continuous with respect to the weak topology.
Let $C\subseteq X$ be a closed subset of $X$ and $r>0$.

\roster
\item"(1)"
We have
 $$
 \underline f^{U,\Lambda}(C,r)
 \le
  \,\,
 \sigma_{\abs}\big(\,\zeta_{C}^{U,\Lambda}(\cdot;r)\,\big) \,.
 \qquad\quad\,\,\,\,
 $$
\item"(2)"
We have
 $$
 \underline f^{U,\Lambda}(C)
 \le
  \liminf_{r\searrow 0}
 \,\,
 \sigma_{\abs}\big(\,\zeta_{C}^{U,\Lambda}(\cdot;r)\,\big) \,.
 $$

\endroster 
\endproclaim

\noindent{\it Proof}\newline
\noindent
(1)
Fix $\varepsilon>0$.
For brevity write
$t=\underline f^{U,\Lambda}(C,r)-\varepsilon$.
Since
$t=\underline f^{U,\Lambda}(C,r)-\varepsilon
<
\underline f^{U,\Lambda}(C,r)
=
 \liminf_{\delta\searrow 0}
 \frac{\log  N_{\delta}^{U,\Lambda}(C,r)}{-\log \delta}$,
we can find
 $\delta_{\varepsilon}$
 with $0<\delta_{\varepsilon}<1$ such that
 $$
 t
 <
 \frac{\log  N_{\delta}^{U,\Lambda}(C,r)}{-\log \delta} 
 $$
for all $0<\delta<\delta_{\varepsilon}$. 
Consequently, for all  $0<\delta<\delta_{\varepsilon}$, we have
 $$
 \delta^{-t}
 \le
  N_{\delta}^{U,\Lambda}(C,r)\,.
  \tag7.1
 $$

 Next, let $c$ denote the constant from Condition (C3) in Section 2.1
 and fix $\rho>0$ with
 $\rho<\min(\, \frac{s_{\min}}{c} \,,\, \delta_{\varepsilon} \,))$.
 We now prove the following two claims.
 
\bigskip

\noindent
{\it
Claim 1.
For $\in\Bbb N$ and $\bold i\in\Sigma^{*}$, the following implication holds:
 $$
 s_{\bold i}\approx \rho^{n}
 \,\,\,\,
 \Rightarrow
 \,\,\,\,
 \rho^{n+1}<s_{\bold i}\le\rho^{n}\,;
 $$
recall, that for $\delta>0$, we write
$s_{\bold i}\approx \delta$
if
$s_{\bold i}\le \delta<s_{\hat\bold i}$, see Section 2.1.
}

\noindent
{\it Proof of Claim 1.} 
Indeed, if $\bold i=i_{1}\ldots i_{m}\in\Sigma^{m}$ with $s_{\bold i}\approx \rho^{n}$, then
$s_{\bold i}\le \rho^{n}<s_{\hat{\bold i}}$, whence
$s_{\bold i}\le \rho^{n}$.
It also follows from Proposition 4.1 that
$s_{\bold i}
=
s_{\hat{\bold i}i_{m}}
\ge
\frac{1}{c}s_{\hat{\bold i}}s_{i_{m}}
>
\frac{1}{c}\rho^{n}s_{\min}
=
\frac{s_{\min}}{c\rho}\rho^{n+1}
\ge
\rho^{n+1}$
where the last inequality is due to the fact that
$\frac{s_{\min}}{c\rho}\ge 1$
because $\rho<\min(\,\frac{s_{\min}}{c}\,,\,\delta_{\varepsilon}\,)\le\frac{s_{\min}}{c}$.
This completes the proof of Claim 1.

\bigskip

\noindent
{\it
Claim 2.
We have
 $$
  \sum
 \Sb
 \bold i\\
 {}\\
 UL_{|\bold i|}[\bold i]\subseteq B(C,r)
 \endSb
 s_{\bold i}^{t}
=
\infty\,.
 $$
}

\noindent
{\it Proof of Claim 2.} 
It is clear that
 $$
 \align
 \sum
 \Sb
 \bold i\\
 {}\\
UL_{|\bold i|}[\bold i]\subseteq B(C,r)
 \endSb
 s_{\bold i}^{t}
&=
 \sum_{n}
 \,
 \sum
 \Sb
 \bold i\\
 {}\\
 \rho^{n+1}<s_{\bold i}\le\rho^{n}\\
 {}\\
 UL_{|\bold i|}[\bold i]\subseteq B(C,r)
 \endSb
 s_{\bold i}^{t}
 \,
 +
 \,
  \sum
 \Sb
 \bold i\\
 {}\\
 \rho<s_{\bold i}\\
 {}\\
UL_{|\bold i|}[\bold i]\subseteq B(C,r)
 \endSb
 s_{\bold i}^{t}\\ 
 &{}\\
&\ge
 \sum_{n}
 \,
 \sum
 \Sb
 \bold i\\
 {}\\
 \rho^{n+1}<s_{\bold i}\le\rho^{n}\\
 {}\\
UL_{|\bold i|}[\bold i]\subseteq B(C,r)
 \endSb
 s_{\bold i}^{t}\,.
 \tag7.2
 \endalign
 $$

Also, for $n\in\Bbb N$ and $\bold i\in\Sigma^{*}$,
the following implication follows from Claim 1:
 $$
 s_{\bold i}\approx \rho^{n}
 \,\,\,\,
 \Rightarrow
 \,\,\,\,
 \rho^{n+1}<s_{\bold i}\le\rho^{n}\,.
 \tag7.3
 $$
We conclude immediately from (7.3) that
 $$
 \sum_{n}
 \,
 \sum
 \Sb
 \bold i\\
 {}\\
 \rho^{n+1}<s_{\bold i}\le\rho^{n}\\
 {}\\
UL_{|\bold i|}[\bold i]\subseteq B(C,r)
 \endSb
 s_{\bold i}^{t}
 \ge
 \sum_{n}
 \,
 \sum
 \Sb
 \bold i\\
 {}\\
 s_{\bold i}\approx \rho^{n}\\
 {}\\
 UL_{|\bold i|}[\bold i]\subseteq B(C,r)
 \endSb
 s_{\bold i}^{t}\,.
 \tag7.4
 $$

Combining (7.2) and (7.4) shows that
 $$
 \align
  \sum
 \Sb
 \bold i\\
 {}\\
 UL_{|\bold i|}[\bold i]\subseteq B(C,r)
 \endSb
 s_{\bold i}^{t}
&\ge
 \sum_{n}
 \,
 \sum
 \Sb
 \bold i\\
 {}\\
 s_{\bold i}\approx \rho^{n}\\
 {}\\
 UL_{|\bold i|}[\bold i]\subseteq B(C,r)
 \endSb
 s_{\bold i}^{t}\\
 &{}\\
&=
 \sum_{n}
 \,
 \sum_{\bold i\in\Pi_{\bold s,\rho^{n}}^{U}(C,r)} 
 s_{\bold i}^{t}\,.
 \tag7.5
 \endalign
 $$
However, if 
${\bold i\in\Pi_{\bold s,\rho^{n}}^{U}(C,r)}$,
then
$s_{\bold i}\approx\rho^{n}$, and it therefore follows from Claim 1 that
$\rho^{n+1}<s_{\bold i}\le \rho^{n}$, whence
$s_{\bold i}\ge \rho^{nt}\rho^{|t|}$.
We conclude from this and (7.5) that
 $$
 \align
  \sum
 \Sb
 \bold i\\
 {}\\
 UL_{|\bold i|}[\bold i]\subseteq B(C,r)
 \endSb
 s_{\bold i}^{t}
&\ge
 \sum_{n}
 \,
 \sum_{\bold i\in\Pi_{\bold s,\rho^{n}}^{U}(C,r)} 
 s_{\bold i}^{t}\\ 
&\ge
\rho^{|t|}
 \sum_{n}
 \,
 \sum_{\bold i\in\Pi_{\bold s,\rho^{n}}^{U}(C,r)} 
 \rho^{nt}\\
&=
\rho^{|t|}
 \sum_{n}
 \,
 \Big|\,\Pi_{\bold s,\rho^{n}}^{U}(C,r)\,\Big|\,\, \rho^{nt}\\ 
&=
\rho^{|t|}
 \sum_{n}
 \,
 N_{\bold s,\rho^{n}}^{U}(C,r)\, \rho^{nt}\,.
 \tag7.6
 \endalign
 $$
Finally, since 
$\rho^{n}
\le
\rho
<\min(\,\frac{s_{\min }}{c}\,,\,\delta_{\varepsilon}\,)
\le
\delta_{\varepsilon}$, we deduce from (7.1) that
$\rho^{-nt}
=
(\rho^{n})^{-t}
\le
N_{\bold s,\rho^{n}}^{U}(C,r)$. This and (7.6) now implies that
 $$
 \align
  \sum
 \Sb
 \bold i\\
 {}\\
UL_{|\bold i|}[\bold i]\subseteq B(C,r)
 \endSb
 s_{\bold i}^{t}
&\ge
\rho^{|t|}
 \sum_{n}
 \,
 \rho^{-nt}\, \rho^{nt}\\  
&=
\rho^{|t|}\sum_{n}1\\
&=
\infty\,.
 \endalign
 $$
This completes the proof of Claim 2.

\bigskip

We conclude
immediately from Claim 2 that
$\underline f^{U,\Lambda}(C,r)-\varepsilon
=
t
\le
\sigma_{\abs}\big(\,\zeta_{C}^{U,\Lambda}(\cdot;r)\,\big)$.
Finally, letting $\varepsilon\searrow 0$ completes the proof.

\noindent
(2)
This follows immediately from (1).
\hfill$\square$

  \bigskip
  \bigskip

 %%%%%%%%%%%%%%%%%%%%%%%%%%%%%%%%%%%%%%%%%
 %%%%%%%%%%%%%%%%%%%%%%%%%%%%%%%%%%%%%%%%%
 %%%%%%%%%%%%%%%%%%%%%%%%%%%%%%%%%%%%%%%%%
 %%%%%%%%%%%%%%%%%%%%%%%%%%%%%%%%%%%%%%%%%
 %%%%%%%%%%%%%%%%%%%%%%%%%%%%%%%%%%%%%%%%% 

\heading
{
8. Proof of Theorem 2.2
}
\endheading

For $x,y\in\Bbb R^{M}$, write
 $$
  [\![x,y]\!]
  =
  \Big\{
 (1-t)x+ty
 \,\Big|\,
 t\in[0,1]
 \Big\}\,,
 $$
 i\.e\.
$[\![x,y]\!]$
denotes the line-segment between $x$ and $y$.

\bigskip

\proclaim{Lemma 8.1}
Let $E\subseteq\Bbb R^{M}$
and
let $x\in E$ and $y\in\Bbb R^{M}\setminus E$.
Then
$ [\![x,y]\!]\cap\partial E\not=\varnothing$.
\endproclaim
\noindent{\it Proof}\newline
\noindent
Let 
$t_{0}
=
\sup\{t\in[0,1]\,|\,
(1-t)x+ty\in E\}$.
Then $(1-t_{0})x+t_{0}y\in [\![x,y]\!]$,
and
since $x\in E$ and $y\in\Bbb R^{M}\setminus E$, it is easily seen that
$(1-t_{0})x+t_{0}y\in\partial E$.
\hfill$\square$

\bigskip

\proclaim{Lemma 8.2}
Let $C\subseteq\Bbb R^{M}$ be a closed subset of $\Bbb R^{M}$
and
let $r,\varepsilon>0$ with $r< \varepsilon$.
Then
$B\big(\,I(C,\varepsilon)\,,r\,\big)\subseteq C$;
recall, that
$I(C,\varepsilon)
=
\{x\in C\,|\,\dist(x,\partial C)\ge\varepsilon\}$,
see Section 2.3.
\endproclaim
\noindent{\it Proof}\newline
\noindent
Let 
$y\in B\big(\,I(C,\varepsilon)\,,r\,\big)$.
We must now prove that $y\in C$.
Assume, in order to reach a contradiction, that $y\not\in C$.
Since 
$I(C,\varepsilon)$ is a closed, 
%subset 
%of the compact set $C$, we deduce that
%$I(C,\varepsilon)$ is compact, and
%we 
%therefore 
it follows that we can
find $x\in I(C,\varepsilon)$
such that
$|y-x|=\dist\big(\,y\,,\,I(C,\varepsilon)\,\big)$.
Also, since $x\in I(C,\varepsilon)\subseteq C$ and $y\not\in C$, it follows from 
Lemma 8.2 that there is $v\in[\![x,y]\!]\cap \partial C$.
We now conclude that
%$r
%\ge
%\dist\big(\,y\,,\,I(C,\varepsilon)\,\big)
%=
%|y-x|
%=
%|y-v|+|v-x|
%>
%|v-x|
%\ge
%\dist\big(\,x\,,\,\partial C\,\big)
%\ge 
%\varepsilon$.
$$
\align
r
&
\ge
\dist\big(\,y\,,\,I(C,\varepsilon)\,\big)
\qquad
\text{[since $y\in B\big(\,I(C,\varepsilon)\,,r\,\big)$]}\\
&=
|y-x|\\
%&=
%|y-v|+|v-x|
%\qquad
%\,\,
%\text{[since $v\in[\![x,y]\!]$]}\\
%&
%>
%|v-x|
%\qquad
%\qquad
%\qquad
%\,\,
%\text{[since $y\not\in C$ and $v\in\partial C\subseteq C$]}\\
&
\ge
|v-x|
\qquad
\qquad
\qquad
\,\,
\text{[since $v\in[\![x,y]\!]$]}\\&
\ge
\dist\big(\,x\,,\,\partial C\,\big)
\qquad
\quad
\,\,\,
\text{[since $v\in\partial C$]}\\
&
\ge 
\varepsilon\,.
\qquad
\qquad
\qquad
\qquad
\,\,
\text{[since $x\in I(C,\varepsilon)$]}
\endalign
$$
However, this inequality contradicts the fact that $r<\varepsilon$.
\hfill$\square$

\bigskip

\noindent{\it Proof of Theorem 2.2}\newline
\noindent
We first note that it follows from Theorem 2.1 that
 $$
 \underline f^{U,\Lambda}(C)\
 =
\sup
  \Sb
  \mu\in\Cal P_{S}(\Sigma^{\Bbb N})\\
  {}\\
  U\mu\in C
  \endSb
  \,\,\,\,
 -\frac{h(\mu)}{\int \Lambda\,d\mu}\,.
 $$
Hence it suffices to prove that
 $$
 \sigma_{\abs}\big(\,\zeta_{C}^{U,\Lambda}\,\big)
=
 \underline f^{U,\Lambda}(C)\,.
 $$
 
We first show that
 $$
 \sigma_{\abs}\big(\,\zeta_{C}^{U,\Lambda}\,\big)
 \le
 \underline f^{U,\Lambda}(C)\,.
 \tag8.1
 $$
Indeed, it follows immediately from the definitions 
of the zeta-functions
$\zeta_{C}^{U,\Lambda}$ and $\zeta_{C}^{U,\Lambda}(\cdot;r)$ that
if $r>0$, then
$\sigma_{\abs}\big(\,\zeta_{C}^{U,\Lambda}\,\big)
\le
\sigma_{\abs}\big(\,\zeta_{B(C,r)}^{U,\Lambda}\,\big)
=
\sigma_{\abs}\big(\,\zeta_{C}^{U,\Lambda}(\cdot;r)\,\big)$, whence
$\sigma_{\abs}\big(\,\zeta_{C}^{U,\Lambda}\,\big)
\le
\liminf_{r\searrow}\sigma_{\abs}\big(\,\zeta_{C}^{U,\Lambda}(\cdot;r)\,\big)$.
We conclude from this  and Theorem 2.1
that
$\sigma_{\abs}\big(\,\zeta_{C}^{U,\Lambda}\,\big)
\le
\liminf_{r\searrow0}\sigma_{\abs}\big(\,\zeta_{C}^{U,\Lambda}(\cdot;r)\,\big)
=
\underline f^{U,\Lambda}(C)$.
This proves (8.1).

Next, we show that
 $$
 \sigma_{\abs}\big(\,\zeta_{C}^{U,\Lambda}\,\big)
 \ge
 \underline f^{U,\Lambda}(C)\,.
 \tag8.2
 $$
Observe that
if $r,\varepsilon>0$ with $r<\varepsilon$, then it follows from Lemma 8.2 that
$B\big(\,I(C,\varepsilon)\,,r\,\big)\subseteq C$, and 
the definitions 
of the zeta-functions
$\zeta_{C}^{U,\Lambda}$ and $\zeta_{C}^{U,\Lambda}(\cdot;r)$ 
therefore imply that
$\sigma_{\abs}\big(\,\zeta_{C}^{U,\Lambda}\,\big)
\ge
\sigma_{\abs}\big(\,\zeta_{B(\,I(C,\varepsilon)\,,r\,)}^{U,\Lambda}\,\big)
=
\sigma_{\abs}\big(\,\zeta_{I(C,\varepsilon)}^{U,\Lambda}(\cdot;r)\,\big)$
for all
$r,\varepsilon>0$ with $r<\varepsilon$.
Hence, for all
$\varepsilon>0$ we have
 $$
 \sigma_{\abs}\big(\,\zeta_{C}^{U,\Lambda}\,\big)
\ge
\liminf_{r\searrow0}
\sigma_{\abs}\big(\,\zeta_{I(C,\varepsilon)}^{U,\Lambda}(\cdot;r)\,\big)\,.
\tag8.3
 $$
Also, since $I(C,\varepsilon)$ is closed, it follows from Theorem 2.1
that
$\liminf_{r\searrow0}
\sigma_{\abs}\big(\,\zeta_{I(C,\varepsilon)}^{U,\Lambda}(\cdot;r)\,\big)
=
\underline f^{U,\Lambda}(\,I(C,\varepsilon)\,)$.
We  conclude
from this and (8.3) that
 $$
 \sigma_{\abs}\big(\,\zeta_{C}^{U,\Lambda}\,\big)
\ge
\underline f^{U,\Lambda}\big(\,I(C,\varepsilon)\,\big)\,.
\tag8.4
 $$
for all $\varepsilon>0$.
Finally, using inner continuity at $C$ and 
letting $\varepsilon\searrow0$, it follows from (8.4)
that
$\sigma_{\abs}\big(\,\zeta_{C}^{U,\Lambda}\,\big)
\ge
\lim_{\varepsilon\searrow0}\underline f^{U,\Lambda}(\,I(C,\varepsilon)\,)
=
\underline f^{U,\Lambda}(C)$.
This proves (8.2).
\hfill$\square$

\bigskip
\bigskip

\heading
{
9. Proof of Theorem 2.3
}
\endheading

The purpose of this section is to prove Theorem 2.3.

\bigskip

\noindent{\it Proof of Theorem 2.3}\newline
\noindent
For brevity write
$G=\{s\in\Bbb C\,|\,\Real(s)>\sigma_{\abs}(\,\zeta_{C}^{U,\Lambda}\,)\}$.
Since 
$\sup_{|\bold i|=n}\frac{1}{\log s_{\bold i}}\to 0$ as $n\to\infty$
(because 
$\sup_{|\bold i|=n}s_{\bold i}\to 0$ as $n\to\infty$),
we conclude that
the series
$Z_{C}^{U,\Lambda}(s)
 =
\sum_{
\bold i\,,\,
UL_{|\bold i|}[\bold i]\subseteq C
}
\frac{1}{\log s_{\bold i}}
s_{\bold i}^{s}$
converges
uniformly in the variable $s$
on all compact subsets of $G$.
%of $s$.
%% for $s$ belonging to compact subsets of $G$.

Since
the series
$Z_{C}^{U,\Lambda}(s)
 =
\sum_{
\bold i\,,\,
UL_{|\bold i|}[\bold i]\subseteq C
}
\frac{1}{\log s_{\bold i}}
s_{\bold i}^{s}$
converges
uniformly in the variable $s$
on all compact subsets of $G$,
we conclude that
the formal calculations below are justified, namely, if $s\in G$, then we have
%This implies that if $s\in G$, 
%then
%the formal calculations below are justified,
%namely
 $$
 \align
 \exp Z_{C}^{U,\Lambda}(s)
 &=
 \exp
\sum
\Sb
\bold i\\
{}\\
UL_{|\bold i|}[\bold i]\subseteq C
\endSb
\frac{1}{\log s_{\bold i}}
s_{\bold i}^{s}\\
{}\\
&=
\exp
\sum
\Sb
\bold i\\
{}\\
\text{$\bold i$ is prime}\\
{}\\
UL_{|\bold i|}[\bold i]\subseteq C
\endSb
\,
\sum_{n}
\,
\frac{1}{\log s_{\!\!\!\!\undersetbrace\text{$n $ times}\to{\ssize{\bold i\ldots\bold i}}}  }
s_{\!\!\!\!\undersetbrace\text{$n $ times}\to{\ssize{\bold i\ldots\bold i}}}^{s}\\
{}\\
&=
\exp
\sum
\Sb
\bold i\\
{}\\
\text{$\bold i$ is prime}\\
{}\\
UL_{|\bold i|}[\bold i]\subseteq C
\endSb
\,
\sum_{n}
\,
\frac{1}{n\log s_{\bold i}}
s_{\bold i}^{sn} 
\\
{}\\
&=
\prod
\Sb
\bold i\\
{}\\
\text{$\bold i$ is prime}\\
{}\\
UL_{|\bold i|}[\bold i]\subseteq C
\endSb
\,
\exp
\,
\Bigg(
\frac{1}{\log s_{\bold i}}
\sum_{n}
\,
\frac{1}{n}
s_{\bold i}^{sn} 
\Bigg)
\\
{}\\
&=
\prod
\Sb
\bold i\\
{}\\
\text{$\bold i$ is prime}\\
{}\\
UL_{|\bold i|}[\bold i]\subseteq C
\endSb
\,
\exp
\,
\Bigg(
\frac{1}{\log s_{\bold i}}
\log
\Bigg(
\frac{1}{1-s_{\bold i}^{s}}
\Bigg)
\Bigg)\\
{}\\
&=
\prod
\Sb
\bold i\\
{}\\
\text{$\bold i$ is prime}\\
{}\\
UL_{|\bold i|}[\bold i]\subseteq C
\endSb
\Bigg(
\frac{1}{1-s_{\bold i}^{s}}
\Bigg)^{\frac{1}{\log s_{\bold i}}}\\
&{}\\
&=
 Q_{C}^{U,\Lambda}(s)\,.
 \tag9.1
\endalign
$$
It follows from 
the calculations involved in establishing
(9.1)
that 
the product
$Q_{C}^{U,\Lambda}(s)$
converges
and that
$Q_{C}^{U,\Lambda}(s)\not=0$
for all $s\in G$.
In addition, we 
deduce from (9.1) that for all $s\in G$, we have
 $\frac{d}{ds}Q_{C}^{U,\Lambda}(s)
 =
  \frac{d}{ds}\exp Z_{C}^{U,\Lambda}(s)
 =
 (\exp Z_{C}^{U,\Lambda}(s))
 \,
\frac{d}{ds}Z_{C}^{U,\Lambda}(s)
 =
  Q_{C}^{U,\Lambda}(s)
 \,
\frac{d}{ds}Z_{C}^{U,\Lambda}(s)$, whence
 $$
 \frac{d}{ds}Z_{C}^{U,\Lambda}(s)
 =
 \frac{\frac{d}{ds}Q_{C}^{U,\Lambda}(s)}{Q_{C}^{U,\Lambda}(s)}
 =
 L\,
 Q_{C}^{U,\Lambda}(s)\,.
 \tag9.2
 $$

Once again using the fact
that
the series
$Z_{C}^{U,\Lambda}(s)
 =
\sum_{
\bold i\,,\,
UL_{|\bold i|}[\bold i]\subseteq C
}
\frac{1}{\log s_{\bold i}}
s_{\bold i}^{s}$
converges
uniformly in the variable $s$
on all compact subsets of $G$.
we deduce that  if $s\in G$, then we have
 $$
 \align
 \frac{d}{ds}
 Z_{C}^{U,\Lambda}(s)
 &=
  \frac{d}{ds}
\sum
\Sb
\bold i\\
{}\\
UL_{|\bold i|}[\bold i]\subseteq C
\endSb
\frac{1}{\log s_{\bold i}}
s_{\bold i}^{s}\\
{}\\
&=
\sum
\Sb
\bold i\\
{}\\
UL_{|\bold i|}[\bold i]\subseteq C
\endSb
\frac{1}{\log s_{\bold i}}
  \frac{d}{ds}
s_{\bold i}^{s}\\
{}\\
&=
\sum
\Sb
\bold i\\
{}\\
UL_{|\bold i|}[\bold i]\subseteq C
\endSb
s_{\bold i}^{s}\\
{}\\
&=
 \zeta_{C}^{U,\Lambda}(s)\,.
 \tag9.3
\endalign
$$

Finally, combining (9.2) and (9.3) gives
$ \zeta_{C}^{U,\Lambda}(s)
  =
   \frac{d}{ds}
 Z_{C}^{U,\Lambda}(s)
 =
   L\, Q_{C}^{U,\Lambda}(s)$
   for all $s\in G$.
\hfill$\square$

\end